\magnification=\magstep 1
\baselineskip=15pt     %normal baselineskip = 12 pt
%\baselineskip=18pt     %normal baselineskip = 12 pt
%\hsize=6.5truein
%\hoffset=.25truein
\parskip=3pt plus1pt minus.5pt
\overfullrule=0pt
\font\hd=cmbx10 scaled\magstep1

\def\Pone{{\bf P}^1}
\def\P{{\bf P}}

\def\O{{\cal O}}
\def\L{{\cal L}}

\def\T{{\cal T}}
\def\H{{\cal H}}

\def\U{{\cal U}}

\def\F{{\cal F}}

\def\M{{\cal M}}

\def\E{{\cal E}}

\def\HH{{\underline H}}

\def\EY{\E(1)|_Y}

\def\PEY{\P(\EY)}

\def\OP{\O_{\PEY}(1)}
\def\OP1{\O_{\Pone}}

\def\Pic{{\rm Pic}}

\def\coker{\mathop{\rm coker}}

\def\re{\mathop{\rm Re}}
\def\im{\mathop{\rm Im}}

\def\mod{\mathop{\rm mod}}
\def\ext{{\rm Ext}}
\def\tor{{\rm Tor}}
\def\hom{{\rm Hom}}

\def\rank{\mathop{\rm rank}}

\def\lhom{{\underline{\hom}}}
\def\lext{{\underline{\ext}}}

\def\boldz{{\bf Z}}

\def\dual#1{{#1}^{\scriptscriptstyle \vee}}

\def\exact#1#2#3{0\rightarrow#1\rightarrow#2\rightarrow#3\rightarrow0}

\def\mapright#1{\smash{
  \mathop{\longrightarrow}\limits^{#1}}}

\def\mapdown#1{\Big\downarrow
   \rlap{$\vcenter{\hbox{$\scriptstyle#1$}}$}}

\def\b{{\bf B}}
\input amssym.tex
\def\relbow{\llcorner}

%\input prepictex.tex
%\input pictex.tex
%\font\thinlinefont=cmr5
\centerline{\hd Special Lagrangian Fibrations II: Geometry.}
\medskip
\centerline{\it Mark Gross\footnote{*}{Supported in part by NSF grant 
DMS-9700761}}
\medskip
\centerline{February 23rd, 1999}
\medskip
\centerline{Mathematics Institute}
\centerline{University of Warwick}
\centerline{Coventry, CV4 7AL}
\centerline{mgross@maths.warwick.ac.uk}
\bigskip
\bigskip
{\hd \S 0. Introduction.}

This paper is a progress report on work surrounding the
Strominger-Yau-Zaslow mirror symmetry conjecture [25]. Roughly put,
this conjecture suggests the following program for attacking an 
appropriate form of the mirror conjecture. Let $X$ be a Calabi-Yau
$n$-fold, with a large complex structure limit point $p$ in a compactification
of the complex moduli space of $X$. One expects mirrors of $X$ to be associated
to such boundary points of the complex moduli space of $X$.
For complex structures on $X$
in an open neighborhood of the boundary point $p$ and suitable choice
of a Ricci flat metric on $X$, one attempts to construct the mirror
of $X$ via the following program:
\item{(1)} There is an $n$-torus representing 
a homology class in $H_n(X,\boldz)$
which is invariant under all monodromy transformations about the
discriminant locus passing through the point $p$. (See [14], \S 3
for details of this representative). The first task is to find an 
homologous $n$-torus which is special
Lagrangian. 
\item{(2)} Having found one special Lagrangian torus, show that it
deforms to yield a fibration $f:X\rightarrow B$ all of whose
fibres are special Lagrangian and whose general fibre is an $n$-torus.
\item{(3)} Construct the dual $n$-torus fibration as follows. Let
$B_0\subseteq B$ be the complement of the discriminant locus of $f$, 
$f_0:X_0\rightarrow B_0$ the restriction of $f$ to $X_0=f^{-1}(B_0)$.
The dual $n$-torus fibration over $B_0$ is 
$\check X_0=R^1f_{0*}{\bf R}/R^1f_{0*}\boldz
\rightarrow B_0$. Find a suitable compactification of $\check X_0$
to a manifold $\check X$ along with a fibration $\check f:\check X\rightarrow
B$.
\item{(4)} Show that $X$ and $\check X$ satisfy a topological form of
mirror symmetry. It is not clear what this means in arbitrary dimension,
but for threefolds, this will be an isomorphism $H^{even}(X,{\bf Q})\cong
H^{odd}(\check X,{\bf Q})$ and $H^{odd}(X,{\bf Q})\cong H^{even}(\check X,
{\bf Q})$. This in particular implies the usual interchange of Hodge numbers
for threefolds. One might also dare to hope these isomorphisms
also hold over $\boldz$; we will see this will often be the case
in Theorem 3.10.
\item{(5)} Put a complex and K\"ahler structure on $\check X$. The choice
of such structures determines the mirror map. One expects that the complex
structure on $X$ should entirely determine the K\"ahler structure on
$\check X$, while the K\"ahler structure on $X$ along with a choice of
the $B$-field,
a cohomology class in $H^2(X,{\bf R}/\boldz)$, or a related cohomology group,
will determine the complex structure on $\check X$. In [15], a somewhat
more precise conjecture was given as to how this interchange of structures
should look on the level of cohomology. 
Specifically, let $\omega$, $\Omega$ be the K\"ahler form and holomorphic
$n$-form on $X$ with $\Omega$ normalised so that $\int_{X_b} \Omega=1$.
In addition, one is given a choice of
$B$-field, which right now we'll take to be a
cohomology class $\b\in H^1(B,R^1f_*{\bf R})$. (The $B$-field will always
be denoted by a bold-face $\b$ to differentiate it typographically from
the base $B$ of the fibration.)
The choice of K\"ahler and complex structures on the mirror, determined
by forms $\check\omega$ and $\check\Omega$, should satisfy the following
relationship: using the identifications $H^1(B,R^1f_*{\bf R})\cong
H^1(B,R^{n-1}\check f_*{\bf R})$ and $H^1(B,R^{n-1}f_*{\bf R})
\cong H^1(B,R^1\check f_*{\bf R})$ which conjecturally hold, 
the following identities should hold in these cohomology groups:
$$\eqalign{[\check\omega]&=[\im\Omega]\cr
[\im\check\Omega]&=[\omega]\cr
[\re\check\Omega]-[\sigma_0]&=\b\cr}$$
where $\sigma_0$ is a chosen zero-section of $\check f:\check X\rightarrow
B$.
\item{(6)} Show that the above procedure yields the correct enumerative
predictions for Gromov-Witten invariants of $X$ and $\check X$.

This program is still a long way from completion, and this paper
represents only one small step in this direction. Not much is known about
items (1) and (2) yet; this may well prove to be the hardest part of
the program. In [15], we gave examples of special Lagrangian $T^3$-fibrations
with a degenerate metric. 
In [29], $T^n$-fibrations are constructed on Calabi-Yau hypersurfaces
in smooth toric varieties. These fibrations are constructed as
deformations of the natural $T^n$-fibration on the large
complex structure limit given by the moment map. Unfortunately these
tori are neither special nor Lagrangian, but this fibration may be sufficient
for a purely topological version of mirror symmetry and provides
evidence for the SYZ conjecture.

We will not address issues (1) and (2) further in this paper.
Instead, throughout this paper,
we assume the existence of a special Lagrangian fibration $f:X\rightarrow B$
on $X$, much as we did in [14]. However, we wish to delve more deeply
into the properties of such fibrations, and to do so, we need to make
reasonable guesses as to what kind of regularity properties such fibrations
will possess. We discuss these issues and assumptions in \S 1. In [14],
we did not make real use of the fact that $f:X\rightarrow B$ was special 
Lagrangian, but rather only used the topological fact that $f$ was a torus
fibration, along with a technical condition we called {\it simplicity}
to control the cohomological contributions of the singular fibres. 
In this paper, we try to make more serious use of the special Lagrangian
condition. We find that we make very heavy use of the Lagrangian condition,
but as yet we do not use the additional special condition
in a profound way. We do however use the fact that special
Lagrangian submanifolds are in fact volume minimizing. 
This allows us to provide moral guidelines as to
what we might expect of special Lagrangian fibrations,
by drawing on the wealth of material known about volume minimizing
rectifiable currents. 
Some of these ideas are discussed in \S 1.

Of course, the study of Lagrangian torus fibrations is a familiar
one in the subject of completely integrable Hamiltonian systems.
In \S 2 of this paper we review and generalise slightly for our
purposes Duistermaat's work on global action-angle coordinates [12]. 
Simplifying a bit, if $f:X\rightarrow B$ is a Lagrangian $T^n$ fibration with
a Lagrangian section
and if $X^{\#}$ is the complement of the critical locus of $f$ in $X$,
then one can canonically write $X^{\#}$ as a quotient of the cotangent
bundle of $B$ by a possibly degenerating family of lattices 
suitably embedded in $\T_B^*$.
Furthermore, the canonical symplectic form on $\T_B^*$ descends to the
symplectic form on $X^{\#}$. Thus if $X$ is a Calabi-Yau manifold,
the existence of a special Lagrangian fibration gives us coordinates
on a large open subset of $X$ on which it is easy to write the
symplectic form. I think of these coordinates as {\it special
Lagrangian coordinates}.
These can be compared with traditional complex coordinates,
where it is easy to write the holomorphic $n$-form ($dz_1\wedge\cdots\wedge
dz_n$), but very difficult to write down the K\"ahler form of a Ricci
flat metric. Thus we should expect that in the coordinates given
by a special Lagrangian fibration, the difficulty will be to write down 
the holomorphic $n$-form, or equivalently, the complex structure.

Taking this analogy further, we note that complex coordinates give
a filtration on de Rham cohomology, namely the Dolbeault cohomology
groups. Similarly, special Lagrangian coordinates can be thought of
as giving rise to a filtration on cohomology, namely that given by
the Leray spectral sequence associated to the special Lagrangian fibration.
Just as the Dolbeault cohomology groups give rise to the Hodge filtration,
we saw in [14] that the Leray filtration should, modulo some
conjectures about monodromy, give rise to the monodromy weight filtration
associated to the large complex structure limit point.

In [14], we studied the Leray spectral sequence of $f$ with coefficients
in ${\bf Q}$. We now show in \S 3, modulo suitable regularity assumptions
stated in \S 1,
that three-dimensional special Lagrangian fibrations satisfy the condition
of $\boldz$-simplicity introduced in [14]. This is a stronger hypothesis
then was used in [14], and as a result, we can analyse the Leray
spectral sequence over $\boldz$. This leads to some interesting
results as to the role torsion in cohomology plays in mirror symmetry.
In particular, it provides an explanation of the phenomenon proposed
in [5] of the role of ``discrete torsion'' (torsion in
$H^3(X,\boldz)$) in mirror symmetry. This
is discussed in \S 3. 
In addition, the analysis of the Leray spectral sequence over $\boldz$
sheds light on item (4) of the program proposed above.

Moving on, we next address the subject of putting symplectic and complex
structures on the mirror $\check X$. This subject has already been discussed
to some extent by Hitchin in [17]. Our approach is inspired by Hitchin's,
and is a generalisation of that approach. In fact, one of the goals
accomplished is to restate Hitchin's constructions in a more
coordinate independent form, so as to allow us to understand
the cohomological ramifications of these constructions.

One should in fact consider any special Lagrangian submanifold
$M\subseteq X$, and consider the moduli space of deformations of this
submanifold. Calling this moduli space $B$, the $D$-brane moduli space
of $M$ is then the set of pairs $(M',\alpha)$ where $\alpha$ is a flat
$U(1)$ connection modulo gauge equivalence on $M'$ a deformation of $M$.
This is a $T^s$-bundle over $B$, where $s=b_1(M)$. Specifically,
if $f:\U\rightarrow B$ is the universal family of special Lagrangian
submanifolds parametrized by $B$, $\U\subseteq B\times X$, then 
the $D$-brane moduli space is $\M=R^1f_*{\bf R}/R^1f_*\boldz\rightarrow B$.
In addition, Mclean [21] gives us a canonical isomorphism between
the tangent bundle of $B$, $\T_B$, and $R^1f_*{\bf R}\otimes C^{\infty}(B)$.
This gives a canonical embedding of $R^1f_*\boldz$ in $\T_B$. However,
$\T_B$ of course does not carry a canonical symplectic form; rather,
it is $\T_B^*$ which does. So to find a symplectic form on $\M$,
we need to reembed $R^1f_*\boldz$ in $\T_B^*$. There are two ways of doing
this. One is to use periods integrals related to $\im\Omega$, the imaginary
part of the holomorphic $n$-form on $X$; the other is to use a canonical
metric introduced by McLean
on $\T_B$ to identify $\T_B$ with $\T_B^*$. In fact these two methods
give the same embedding of $R^1f_*\boldz$ in $\T_B^*$. This allows us
to define a symplectic form on $\M$ by writing $\M$ as
$\T_B^*/R^1f_*\boldz$ and taking the form on $\M$ induced by
the canonical symplectic form on $\T_B^*$. This is the same method
as proposed by Hitchin, recast in a slightly more invariant way,
which makes it easy to see that the cohomology class of the symplectic
form defined in this manner is as predicted by the conjecture of [15].
Thus, if we can understand how to extend this symplectic form
to the compactification of the $D$-brane moduli space, we will solve
the first half of item (5). This circle of ideas is discussed in \S 4.

Because it is easy to write the symplectic form in special Lagrangian
coordinates, we expect it to be very difficult to write the complex
structure. We take up this issue in \S\S 5 and 6. First we explore
what data is necessary to place a complex structure on a Lagrangian
torus fibration in order to make the fibration special Lagrangian
and the induced metric Ricci-flat. A moment's thought shows that given
knowledge of $\omega$, to specify an almost complex structure 
compatible with $\omega$, it is enough to give for 
each point $b\in B$ a metric $g$ on the fibre $X_b$,
and for each point $x\in X_b$ the Lagrangian subspace $J(\T_{X_b,x})
\subseteq \T_{X,x}$. Then $J$ is completely determined by the requirement
that $g(v,w)=\omega(Jv,w)$. The collection of subspaces $J(\T_{X_b,x})$
can be thought of as the horizontal subspaces of an Ehresmann connection
on $f:X^{\#}\rightarrow B$.
We then need to ask when this data determines an integrable complex structure
which induces a Ricci-flat metric in which the fibres of $f$ are
special Lagrangian. 

Following [17], it is easier to determine this by describing the
holomorphic $n$-form $\Omega$. In local coordinates $y_1,\ldots,y_n$
on the base, and canonical coordinates $x_1,\ldots,x_n$ on the fibres
of $\T_B^*$, we can write the general form of $\Omega$ as
$$\Omega=V\bigwedge_{i=1}^n (dx_i+\sum_{j=1}^n \beta_{ij} dy_j).$$
where $\beta_{ij}$ is a complex-valued function on $\T_B^*$ and $V$
is a real-valued function.
The Ehresmann connection is in fact encoded
in $\re\beta$, while $\im\beta$ is just the inverse of the metric on
the fibres. 
The integrability and Ricci-flatness conditions are then easy to
write down. In fact, we show that one needs the conditions:
\item{(1)} The matrix $\beta=(\beta_{ij})$ is symmetric, $\im\beta$
is positive definite, and $V=1/\sqrt{\det\im\beta}$.
\item{(2)} $d\Omega=0$.

The first condition is of course easily achieved, but the second condition
is a quite subtle condition. 
The second condition is what requires real effort, and 
understanding it is really at the heart of the SYZ program. 
It turns out that (2)
is equivalent to the following three conditions:
\item{(1)} The almost complex structure is integrable. (Note that
$d\Omega=0$ is a much stronger condition than integrability. Given any
Lagrangian fibration on a Calabi-Yau manifold with a holomorphic
$n$-form $\Omega$, by replacing $\Omega$
by $e^{i\theta(x)}\Omega$ for a suitable function $\theta(x)$, one can ensure
that $\im\Omega$ restricted to each fibre is zero. However now $\Omega$
is no longer holomorphic.)
\item{(2)} $dx_1,\ldots,dx_n$ are harmonic 1-forms on $X_b$ for each
$b\in B$. (That this is a necessary condition follows from Mclean's results.)
\item{(3)} The volume form $Vdx_1\wedge\cdots\wedge dx_n$ on fibres is parallel
under translation via the Ehresmann connection.

In [17], Hitchin constructed a complex structure on the D-brane moduli
space by specifying the holomorphic $n$-form $\Omega$. The
form he constructed had very special properties: he used $\re\beta=0$
and $\im\beta$ constant along fibres. As a result, the integrability
condition $d\Omega=0$ was much easier to check. In the general case,
it requires a great deal more effort to analyse this equation,
and the results of \S 5 and 6 represent only a beginning. The conditions
(1)-(3) above should perhaps be thought of as mirror equations to the usual
complex Monge-Amp\`ere equations which arise in the study of Ricci-flat
metrics. Understanding their solution should be one of the key
steps in the Strominger-Yau-Zaslow program. We analyse the solution
in several simple cases. For example, we show that if the Ehresmann
connection is in fact flat, then the metric on each fibre must be
flat. This explains why in Hitchin's situation, where
the connection was trivial, it would be impossible to obtain any
solutions which are not flat along the fibres. In particular, the
connection clearly cannot be flat if the fibration possesses singular
fibres. We also give a whole family of solutions related to Hitchin's,
but which take the $B$-field into account. This gives some hint
as to the role the $B$-field plays and leads us to a refined form of
the mirror symmetry conjecture (Conjecture 6.6). We will argue that the
$B$-field should not be thought of as an element of $H^2(X,{\bf R}/\boldz)$
but rather as an element of $H^1(B,R^1f_*{\bf R}/\boldz)$. Thus the group
the $B$-field takes values in should depend not just on $X$ but on the 
fibration $f$. These two groups can in fact be different, so this
suggestion is quite a serious modification of previous interpretations
of the $B$-field (for example that of [5]).

Finally, in \S 7 we apply much of what we have done to the situation
of K3 surfaces, by way of an extended example. Here we know that
special Lagrangian fibrations do exist, and we know the precise construction
of the mirror map. We are able to show that the various recipes
and constructions given in earlier sections of the paper can
be carried out completely for K3 surfaces. In particular, we
obtain an almost purely differential geometric description of mirror
symmetry for K3 surfaces.

I would like to thank P. Aspinwall, D. Calderbank,
N. Hitchin, M. Micallef, A. Todorov,
P. Wilson, S.-T. Yau,  and E. Zaslow
for useful discussions. In addition, I would like to thank S.-T. Yau
for his hospitality at Harvard, where some of this work was carried out.
\medskip

{\it Convention.} For a $p$-form $\alpha$ and tangent vectors $v_1,\ldots,
v_q$, $\iota(v_1,\ldots,v_q)\alpha$ denotes the $p-q$-form 
$\alpha(v_1,\ldots,v_q,\cdot)$.
\bigskip

{\hd \S 1. Special Lagrangian Fibrations.}

In what follows, $X$ will denote a complex Calabi-Yau manifold with a 
Ricci-flat K\"ahler metric $g$, K\"ahler form $\omega$, and holomorphic
$n$-form $\Omega$. The form
$\Omega$ is always normalised to be of unit length, i.e.
$$\omega^n/n!=(-1)^{n(n-1)/2}(i/2)^n\Omega\wedge\bar\Omega.$$
This only fixes $\Omega$ up to a phase factor.
As the metric is Ricci-flat, we can take the metric to be real analytic.

Recall

\proclaim Definition-Proposition 1.1. [16] 
$\re\Omega$ is a calibration, called the {\it special Lagrangian}
calibration. An $n$-dimensional real submanifold $M\subseteq X$
is {\it special Lagrangian} if $\re\Omega|_M=Vol(M)$.
Modulo orientation,
a submanifold $M\subseteq X$ of real dimension $n$ is special
Lagrangian if and only if $\omega|_M=0$ and $\im\Omega|_M=0$.

It is natural to extend the notion of special Lagrangian
submanifolds to special Lagrangian integral currents, as done
in [16]. For an introduction to integral and rectifiable currents
and geometric measure theory, see [22], and for more rigourous
treatments see [13] and [24]. We will not make much
use of the language of geometric measure theory here except to justify
some of the assumptions made on special Lagrangian fibrations. The
reader unfamiliar with this language should just keep in mind that passing
from submanifolds to integral
currents means extending the class of submanifolds to subsets with integer
multiplicities attached, over which
one can still integrate forms.
There are natural compactness theorems
for integral currents, enabling one to easily construct volume minimizing
currents. One can then try to control the singularities of such currents.
The biggest regularity result of this nature is Almgren's monumental

\proclaim Theorem 1.2. [1,2] Suppose $N$ is an $m+l$ dimensional submanifold
of ${\bf R}^{m+n}$ of class $k+2$ and that $T$ is an $m$-dimensional
rectifiable current in ${\bf R}^{m+n}$ which is absolutely
area minimizing with respect to $N$. Then there is an open subset
$U$ of ${\bf R}^{m+n}$ such that $Supp(T)\cap U$ is an $m$-dimensional
minimal submanifold of $N$ of class $k$ and the Hausdorff dimension
of $Supp(T)-(U\cup Supp(\partial T))$ does not exceed $m-2$.

The proof [2] is unpublished and is over 1500 pages long in preprint form.
Thus an urgent question is

\proclaim Question 1.3. Are there nice regularity theorems for special
Lagrangian currents?

For example, [18] shows that holomorphic integral currents are obtained
by integration over complex analytic subvarieties, and of course
the singular locus of such varieties is well-behaved. This is
much stronger than Almgren's result, which doesn't even guarantee finite
$m-2$ dimensional Hausdorff measure of the singular set. We do not want
the theory of special Lagrangian currents to have to depend on Almgren's
result.

Next we should consider what is a reasonable definition of a special
Lagrangian fibration. It might be reasonable to say that

\proclaim Definition 1.4. If $B$ is a topological space and
$f:X\rightarrow B$ is a continuous map, we say $f$ is a {\it
special Lagrangian fibration} if for all $b\in B$, $X_b:=f^{-1}(b)$
is the support of a special Lagrangian integral current $T$ with
$\partial T=0$.

Even this might be too strong; one might insist that the fibres be
special Lagrangian
on a dense subset of $B$. This would allow some fibres to jump
dimension. Nevertheless, we do not expect fibres to
decrease in dimension as this would suggest the cohomology class of the general
fibre $X_b$ was trivial, which contradicts $\int_{X_b}\Omega\not=0$.
I do not want to give a rigorous argument of this sort here as it
requires being clearer about concepts such as the dimension of the fibre.
However it is clear that special Lagrangian fibrations cannot behave
as many completely integrable Hamiltonian systems do, in which
some fibres are tori of smaller dimension. On the other hand, it
is not as clear that we want to rule out the possibility of fibres
jumping up in dimension, something which often happens in algebro-geometric
contexts. Nevertheless, we will stick to Definition 1.4.

We are interested in very specific sorts of special Lagrangian $T^n$
fibrations. As argued in [14], \S 3, we are looking for special
Lagrangian fibrations on Calabi-Yau manifolds near a specific large
complex structure limit point in the boundary of complex moduli space, and the
homology class of a fibre should be represented by a specific vanishing
cycle associated to the boundary point. It was argued in [14], Observation
3.4 that, in the three dimensional case, if this homology class is
primitive, then for {\it general} choice of complex moduli near the large
complex structure limit point, all fibres of $f:X\rightarrow B$ must
be irreducible. Let us be more precise. There is a notion of {\it indecomposable
integral current} ([13], 4.2.25) which is analagous to the notion 
of irreducibility in algebraic geometry, and we will say a fibre
$f^{-1}(b)$ is irreducible if the current $T_b$ obtained by integrating
forms over $f^{-1}(b)$ (with orientation induced by $\Omega$) is 
indecomposable.

\proclaim Definition 1.5. A special Lagrangian fibration is {\it integral}
if each fibre is irreducible and the currents $\{T_b|b\in B\}$ all
represent the same integral homology class.

This last condition rules out fibres which need to be thought of as multiple
fibres. However, it is difficult to define multiplicity of a fibre
in the context of a continuous map. Here the term ``integral''
is being used in its algebro-geometric
sense: each fibre is irreducible and ``reduced''.

Now a special Lagrangian fibration need not be integral any more than
an elliptic fibration need only have integral fibres. But this is the generic 
behaviour of special Lagrangian fibrations in dimensions $\le 3$, and
the assumption of integrality vastly decreases the range of
potential singular fibres. Without assumptions of integrality, it would
be very difficult to relate the cohomology of a $T^n$-fibration and its dual.
But one should keep in mind that even in dimension $3$, integrality
should fail for special values of the complex structure. Furthermore,
if the homology class of a fibre is not primitive, then there is still
the possibility that integrality will fail, and perhaps
even multiple fibres appear. I have
no argument to rule this out, but as we shall see in later sections, it is
special Lagrangian fibrations with sections which are most important, and 
for these the homology class of a fibre is primitive.

The next natural question is whether we expect $B$ itself to be a manifold. I
would like to give a rough argument
that this is a natural expectation if the fibration
$f:X\rightarrow B$ is integral. Indeed, given a fibre $X_b$,
$X_b$ will be smooth at a general point $x\in X_b$, and in a neighborhood
of $x$, using the exponential map, the deformations of $X_b$ can be
identified with deformations of $X_b$ inside its normal bundle near $x$.
Thus locally the fibre of the normal bundle of $X_b$ at $x$ yields a 
natural local section for $f:X\rightarrow B$, and hence gives a manifold
structure on $B$. This construction hopefully in addtion yields a
$C^{\infty}$ structure on $B$.
Thus we will feel justified in making the following assumption, to be in
force throughout the remainder of the paper.

\proclaim Assumption 1.6. All special Lagrangian fibrations $f:X\rightarrow
B$ will be $C^{\infty}$ maps of $C^{\infty}$ manifolds. Furthermore, $f$ will
be assumed to have a local section at each point $b\in B$. In addition,
if $f$ is assumed to be integral, we will assume that for any point
$x\in X_b-Sing(X_b)$, there is a local $C^{\infty}$ section passing through
$x$.

Again, if the fibration is not integral, I would not be surprised
if singular $B$ can arise naturally. Also, since the metric on $X$
is real analytic, we can hope that with suitable coordinates on
$B$, $f$ will in fact be real analytic. We will not assume that here, however.

Next, we pass to the nature of the discriminant locus. Given $f$
integral as in Assumption 1.6, we can now consider the rank of
the differential $f_*:\T_{X,x}\rightarrow \T_{B,f(x)}$ at various points
$x\in X$. At points with a local $C^{\infty}$ section, $\rank f_*=n$.
It will be shown in Proposition 2.2 that if $\rank f_*=\nu$, then the fibre
in fact contains a submanifold of dimension $\nu$ on which $\rank f_*=\nu$.
Thus the existence of points $x\in X$ for which $\rank f_*=n-1$ contradicts
Almgren's theorem. Even if one doesn't accept Almgren's theorem,
it is not difficult to rule out the existence of such points
in low dimension using regularity results for minimizing hypersurfaces.
Thus it is quite safe to assume that $\rank f_*\not= n-1$ for any
point $x\in X$.
This gives a stratification of the discriminant locus
$$\Delta_{\nu}=f(\{x\in X|\rank f_*:\T_{X,x}\rightarrow \T_{B,f(x)}\le\nu\},$$
with 
$\Delta_0\subseteq\cdots\subseteq \Delta_{n-2}=\Delta$, the discriminant 
locus. Then Federer's generalization of Sard's Theorem, [13], 3.4.3, states
that $\H^{\nu+(2n-\nu)/k}(\Delta_{\nu})=0$, where $\H^n$ denotes
the $n$-dimensional Hausdorff measure. In particular, if $k=\infty$,
$\H^{\nu+\epsilon}(\Delta_{\nu})=0$ for all $\epsilon>0$, so $\Delta_{\nu}$
is of Hausdorff dimension $\le \nu$. We need however stronger information
than Hausdorff dimension to reach any cohomological conclusions.
Thus the following assumptions will be used at various points
of this paper; unlike Assumption 1.6, we will assume these only
when we need them, mentioning them specifically.

\proclaim Assumption 1.7. For an integral special Lagrangian fibration
$f:X\rightarrow B$, $b\in \Delta_{\nu}-\Delta_{\nu-1}$, there is
a dense subset $\L$ of the real Grassmannian of $n-\nu$ dimensional
subspaces of $\T_{B,b}$ such that one can find locally, for each
$L\in\L$, a submanifold $B'$ of $B$ passing through $b$ such that
$\T_{B',b}=L$ and $B'\cap\Delta_{\nu}=\{b\}$.

This requires reasonable regularity results about the discriminant,
which certainly hold if $f:X\rightarrow B$ is real analytic, so that
$\Delta_{\nu}$ is a sub-analytic set. A stronger form of this assumption
which we will need in \S 3 and will comment on there is

\proclaim Assumption 1.7$'$. In addition in Assumption 1.7, $B'$ can
be found so that $f^{-1}(B')$ is a submanifold in a neighborhood
of the fibre $X_b$.

Finally, we will require some assumptions on $Sing(X_b)$. Almgren's
theorem that the Hausdorff dimension  of $Sing(X_b)$ is no more
than $n-2$ is not sufficiently strong for most purposes,
and at the very least, we will frequently need to use

\proclaim Assumption 1.8. If $f:X\rightarrow B$ is an integral special
Lagrangian fibration then $$H^i(Sing(X_b),\boldz)=0$$ for $i>n-2$, for all
$b\in B$.

This is a restriction on homological dimension. If $f$ is real analytic, then
this should hold given Almgren's result.
\bigskip

{\hd \S 2. Action-angle Coordinates.}

There is a standard theory of global action-angle coordinates
due to Duistermaat [12]. We will extend this slightly so as to include
information about the smooth part of the singular fibres.
In this section, $f:X\rightarrow B$ denotes any special Lagrangian
fibration satisfying Assumption 1.6. However, many of the results
in this section apply to $C^k$ Lagrangian fibrations $f:X\rightarrow
B$ with $B$ a manifold and $f$ having a local
section in a neighborhood of each point $b\in B$.

The first observation is that there is an action of $\T_B^*$, the
cotangent bundle of $B$, on $X$.

\proclaim Proposition-Definition 2.1. There is a $C^{\infty}$ 
action of $\Gamma(U,\T_B^*)$ on
$f^{-1}(U)$ for any $U\subseteq B$, which we write for any $\alpha\in
\Gamma(U,\T_B^*)$ as
a map $T_{\alpha}:f^{-1}(U)\rightarrow f^{-1}(U)$.
This satisfies the following properties:
\item{(1)} If $d\alpha=0$, then $T_{\alpha}$ is a symplectomorphism.
\item{(2)} $T_{\alpha}$ acts fibrewise, and $T_{\alpha}|_{f^{-1}(b)}:f^{-1}(b)
\rightarrow f^{-1}(b)$ only depends on the value of $\alpha$ at $b$.
\item{(3)} $T_{\alpha}\circ T_{\beta}=T_{\alpha+\beta}$.

Proof. This is standard: see for example [12], \S 1, or [3], Chap.
10. We review the definition of these maps however.
If $\alpha$ is a compactly supported 1-form on $Y$, then $f^*\alpha$
is a compactly supported 1-form on $f^{-1}(U)$. There is then a 
vector field $v_{\alpha}$ on $f^{-1}(U)$ with $\iota(v_{\alpha})\omega
=f^*\alpha$. 
This generates a flow $\phi_t:f^{-1}(U)\rightarrow f^{-1}(U)$
for all $t$, and we take $T_{\alpha}=\phi_1$. It is then standard that if
$d\alpha=0$, $v$ is locally Hamiltonian and $\phi_t$ then is  a 1-parameter
family of symplectomorphisms.

If locally $d\alpha=0$, we can write $\alpha=dH$ for some function $H$
on $B$, and if $G$ is any other function on $B$,
$\{G\circ f, H\circ f\}=0$ since $f:X\rightarrow B$ is a Lagrangian fibration.
But then $\phi_t$ is the Hamiltonian flow associated to $H$, and 
$$0=\{G\circ f, H\circ f\}(x)={d\over dt}\bigg |_{t=0}(G\circ f)(\phi_t(x))$$
so $G\circ f$ is a constant on $\phi_t(x)$ for any $x$. Thus $\phi_t$
acts on fibres. Clearly $v_{\alpha}|_{f^{-1}(b)}$ depends only on the value of
$\alpha$ at $b$, so the action of $\phi_t$ on $f^{-1}(b)$ depends only on the
value of $\alpha$ at $b$. In particular $\phi_t$ acts on fibres for arbitrary
$\alpha$, not just compactly supported $\alpha$. $\bullet$

\bigskip
Next, a standard analysis of the orbits of this action.

\proclaim Proposition 2.2. If $b\in B$, $x\in f^{-1}(b)$, then the orbit
$\{T_{\alpha}(x)|\alpha\in T_{B,b}^*\}$ is diffeomorphic to ${\bf R}^l\times
T^s$ and $l+s$ coincides with $rank(f^*:\T_{B,b}^*\rightarrow \T_{X,x}^*)
=rank(f_*:\T_{X,x}\rightarrow\T_{B,b})$.

Proof. If $\alpha\in\T_{B,b}^*$ is in the kernel of $f^*$, then $v_{\alpha}$
is zero at $x$, so $T_{\alpha}(x)=x$. Thus for any $\alpha\in\T_{B,b}^*$,
$T_{\alpha}(x)$ depends only on $\alpha$ modulo $\ker f^*$, and the orbit
of $x$ is homeomorphic to a quotient of $V=\T_{B,b}^*/\ker f^*$ by a subgroup
$\Gamma$, via the map $\alpha\in V\mapsto T_{\alpha}(x)$. The differential of
this map at $0\in V$ is injective, and hence the map $V\rightarrow X$ is a
diffeomorphism of an open neighborhood of $0\in V$ with its image in $X$. Thus
$\Gamma$ is a discrete subgroup of $V$, and the orbit of $x$ is diffeomorphic
to $V/\Gamma$. $\bullet$
\bigskip

Suppose first that $f:X\rightarrow B$ is equipped with a $C^{\infty}$ section
$\sigma_0:B\rightarrow X$. We define a $C^{\infty}$
map 
$\pi:\T_B^*\rightarrow X$ by, for $\alpha\in \T_{B,b}^*$, $\pi(\alpha)
=T_{\alpha}(\sigma_0(b))$. The image of the zero section of $\T_B^*$ is
$\sigma_0(B)$. Let $\Lambda\subseteq\T_B^*$ be given by
$\Lambda=\pi^{-1}(\sigma_0(B))$.
Then $\Lambda_b\subseteq \T_{B,b}^*$ is the discrete subgroup
of the vector space $\T_{B,b}^*$ given by $\Lambda_b=\pi^{-1}(\sigma_0(b))$.
Let $X_0^{\#}$ be the image of the map $\pi$, $f^{\#}:X_0^{\#}\rightarrow
B$ the projection. Then clearly $\Lambda_b$ is canonically isomorphic to
$H_1((f^{\#})^{-1}(b),\boldz)$, which is isomorphic to
$H^{n-1}_c((f^{\#})^{-1}(b),\boldz)$ by Poincar\'e duality. Here we use
the fact that $f$ is special Lagrangian to give a canonical orientation
on the fibres of $f^{\#}$.
Also, $\Lambda$, as a subset of the total space $\T_B^*$, is closed
as $\pi$ is continuous.
Since the map $\pi$ is a local isomorphism,
$\Lambda$ is also \'etale over $B$. Thus we can think of $\Lambda$ as the
\'espace \'etal\'e of $R^{n-1}_cf^{\#}_*\boldz$, and in particular, we obtain
an exact sequence of sheaves of abelian groups
$$\exact{R^{n-1}_cf_*^{\#}\boldz}{\T_B^*}{X_0^{\#}},$$
where this now defines the group structure on $X_0^{\#}$. Observe
also that since
$\Lambda$ is \'etale over $B$, $R^{n-1}_cf_*^{\#}\boldz$
has no sections with support in proper closed subsets of $B$. In particular,
$\HH^0_{\Delta}(B,R^{n-1}_cf_*^{\#}\boldz)=0$.

We recall the notion of canonical coordinates on the total space of the
bundle $\T_B^*$. Given $U\subseteq B$ an open set with coordinates
$y_1,\ldots, y_n$, {\it canonical coordinates} on $\T_U^*$
are $y_1,\ldots,y_n,x_1,\ldots,x_n$, where $(y_1,\ldots,y_n,x_1,\ldots,x_n)$
is the coordinate representation of the differential form
$\sum x_idy_i\in \T_{U,(y_1,\ldots,y_n)}^*$. We will use canonical 
coordinates consistently throughout this paper, whenever we use local
coordinates to perform calculations. Note that $\T_B^*$ always carries
a standard symplectic form, which in canonical coordinates can be written
as $\sum_{i=1}^n dx_i\wedge dy_i$. (This is the opposite of some sign
conventions).

\proclaim Proposition 2.3. With notation as above,
let $y_1,\ldots,y_n$ be local coordinates on
a neighborhood $U\subseteq B$. Then on $\T_U^*$,
$$\pi^*\omega=\sum_i dx_i\wedge dy_i +\sum_{i,j} a_{ij} dy_i\wedge dy_j$$
where the $a_{ij}$ are functions depending only on $y_1,\ldots,y_n$.
Furthermore, if $\sigma_0:B\rightarrow X$ is a Lagrangian section, then
$$\pi^*\omega=\sum_i dx_i\wedge dy_i$$
on $U$, and thus $\pi^*\omega$ is the standard symplectic form on
$\T_B^*$. Finally, if $H^2(B,{\bf R})=0$, then every section of $f$ is
homotopic to a Lagrangian section.

Proof. Since the fibres of $\T_B^*\rightarrow B$ are Lagrangian with respect to
$\pi^*\omega$, we can locally write
$$\pi^*\omega=\sum_i dx_i\wedge \theta_i+\sum_{i,j} a_{ij}dy_i\wedge dy_j$$
with the $\theta_i$ 1-forms not involving the $dx_i$'s and $a_{ij}$ functions
on $\T_B^*$. Now the function $y_i$ induces on $X$ a Hamiltonian vector
field which, by definition of the map $\pi$, must be $\pi_*\partial/\partial
x_i$. Thus $\iota(\partial/\partial x_i)(\pi^*\omega)=dy_i$, from which we see
that $\theta_i=dy_i$.

The condition $d(\pi^*\omega)=0$ then implies that the functions $a_{ij}$
are independent of $x_1,\ldots,x_n$.

If $\sigma_0$ is Lagrangian with respect to $\omega$, then the zero section
of $\T_B^*$ is Lagrangian, from which we see that $a_{ij}=0$.

If $\sigma_0$ is not Lagrangian, let $\omega'=\sum_i dx_i\wedge dy_i$
(locally)
be the standard symplectic form on $\T_B^*$. Then $\pi^*\omega-\omega'$
is a closed 2-form locally given by $\sum_{i,j} a_{ij}dy_i\wedge dy_j$,
and hence is the pull-back of a closed 2-form on $B$. Thus if $H^2(B,{\bf
R})=0$, there exists a one-form $\theta$ on $B$ with
$d\theta=\pi^*\omega-\omega'$. Then $-\theta$ defines a section of $\T_B^*$
which is Lagrangian with respect to $\pi^*\omega$, and this maps to
a Lagrangian section of $f$ homotopic to $\sigma_0$. 
$\bullet$

This allows us to prove a result stated in [14]. Recall that $X^{\#}$
is the complement of the critical locus of $f$ in $X$.

\proclaim Theorem 2.4. (Theorem 3.6 of [14])
Let $X$ be a Calabi-Yau $n$-fold, $B$
a smooth real $n$-dimensional manifold, with $f:X\rightarrow B$ a
$C^{\infty}$ special Lagrangian torus fibration such that $R^nf_*{\bf Q}
={\bf Q}_B$ and such that the singular locus of each singular fibre has
cohomological dimension $\le n-2$.
Suppose furthermore that $f$ has a $C^{\infty}$ section $\sigma_0$.
Then $X^{\#}$ has the structure of a fibre space
of groups with $\sigma_0$ the zero
section. In fact there is an exact sequence of sheaves of abelian groups
$$\exact{R^{n-1}f_*\boldz}{\T_B^*}{X^{\#}}.$$
Given a section $\sigma\in \Gamma(U,X^{\#})$, one obtains a
$C^{\infty}$ diffeomorphism $T_{\sigma}:f^{-1}(U)\cap X^{\#}
\rightarrow f^{-1}(U)\cap X^{\#}$ given by $x\mapsto x+\sigma(f(x))$,
and this diffeomorphism extends to a diffeomorphism $T_{\sigma}:f^{-1}(U)
\rightarrow f^{-1}(U)$.

Proof. This follows from the previous discussion
if we can show that the hypotheses imply two
things: 
\item{(1)} $X^{\#}=X_0^{\#}$
\item{(2)} $R^{n-1}_cf^{\#}_*\boldz\cong R^{n-1}f_*\boldz$.

To show (1) and (2),
let $X_b$ be a fibre, and let $Z\subseteq X_b$ be the singular
locus of $Supp(X_b)$. 
We then have an exact sequence
$$\eqalign{H^{n-2}(Z,{\bf Q})&\rightarrow H^{n-1}_c(X_b-Z,{\bf Q})
\rightarrow H^{n-1}(X_b,{\bf Q})\cr\rightarrow &H^{n-1}(Z,{\bf Q})=0
\rightarrow H^n_c(X_b-Z,{\bf Q})\rightarrow H^n(X_b,{\bf Q})
\rightarrow 0.\cr}$$
Here $H^{n-1}(Z,{\bf Q})=0$ by the assumption on the cohomological
dimension of $Z$.
Thus, since $H^n(X_b,{\bf Q})={\bf Q}$ by assumption, $X_b-Z$ can have only
one connected component. Since the $\T_{B,b}^*$-orbit
of $\sigma_0(b)$ is already one connected component of 
$X_b-Z$, we see that $X_b-Z=(X_0^{\#})_b$. But as $X_b^{\#}\subseteq
X_b-Z$, we must have $X_0^{\#}=X^{\#}$. 

We also see from the above exact sequence that there is a surjection
$$R^{n-1}_cf_*^{\#}\boldz\rightarrow R^{n-1}f_*\boldz\rightarrow 0.$$
This is an isomorphism outside of $\Delta$. But since
$\HH^0_{\Delta}(B,R^{n-1}_cf_*^{\#}\boldz)=0$, we conclude this surjection is 
in fact an isomorphism. $\bullet$

Note that the hypotheses of Theorem 2.4 hold if $f:X\rightarrow B$ is integral
and satisfies Assumption 1.8.
Another useful observation about the topology of singular fibres:

\proclaim Lemma 2.5. If $f:X\rightarrow B$ is integral and
satisfies Assumption 1.8, and $b\in B$ with
$Z=Sing(X_b)$, then $Z$ is connected.

Proof. Let $U=X_b-Z$. As argued above, $U\cong {\bf R}^k\times T^{n-k}$
for some $k$. Assuming $Z$ is non-empty, we have an exact sequence
$$0=H^0_c(U,\boldz)\rightarrow H^0(X_b,\boldz)\rightarrow H^0(Z,\boldz)
\rightarrow H^1_c(U,\boldz)\rightarrow H^1(X_b,\boldz)$$
and we wish to show $H^0(Z,\boldz)=\boldz$, which is equivalent to 
the injectivity of $H^1_c(U,\boldz)\rightarrow H^1(X_b,\boldz)$. 
If $H^1_c(U,\boldz)=0$, then there is no problem, so the only possiblity
is that $k=1$ and $H^1_c(U,\boldz)=\boldz$. In this case, $H^{n-1}_c(U,\boldz)
=\boldz^{n-1}$ and thus in a small neighborhood $V\subseteq B$ of $b$,
the sheaf $R^{n-1}_cf_*^{\#}\boldz$ contains $\boldz^{n-1}$ as a subsheaf.
Thus over $V$, there is a $T^{n-1}$ bundle $T\rightarrow B$,
$T\subseteq X^{\#}$. The fundamental class of each fibre yields a non-zero
section $\sigma$ of $R^1_cf_*^{\#}\boldz$ over $B$, and under the map
$R^1_cf_*^{\#}\boldz\rightarrow R^1f_*\boldz$, $\sigma$ maps to a section
of $R^1f_*\boldz$ which is non-zero on $V-\Delta$. But the map $H^1_c(U,\boldz)
\rightarrow H^1(X_b,\boldz)$ is the induced map on stalks, and hence
is non-zero, thus injective. $\bullet$
\bigskip

We now address the situation when $f:X\rightarrow B$ does not have
a Lagrangian section but we assume $f$ is integral. 
We recall though that Assumption 1.6 is always in force, giving the existence
of a local
section. This theory was developed in Duistermaat's paper
[12], and is completely analagous to Kodaira's theory of elliptic
surfaces, or to Ogg-Shafarevich theory.

Let $X^{\#}=X-Crit(f)$ as usual. We now obtain a map 
$$\psi:R^{n-1}_cf_*^{\#}\boldz\rightarrow \T_B^*$$
as follows. Given a fibre $X_b^{\#}$, $\gamma\in H_1(X_b^{\#},\boldz)
\cong H^{n-1}_c(X_b^{\#},\boldz)$,
map $\gamma$ to the differential
$$v\mapsto -\int_{\gamma}\iota(v)\omega$$
where we choose any lifting of $v\in\T_{B,b}$ to $X^{\#}$.

Now in the case $f$ did have a section, we previously 
constructed an embedding $R^{n-1}_cf_*^{\#}\boldz\hookrightarrow\T_B^*$.
We compare these two constructions.
Let $U\subseteq B$ be an open set where $f^{-1}(U)\rightarrow U$
possesses a section, which we can take to be Lagrangian. 
Using this section as the zero section
we obtain an exact sequence
$$0\rightarrow \Lambda=R^{n-1}_cf_*^{\#}\boldz|_U\mapright{\psi'}
{\T_U^*}\mapright{\pi} f^{-1}(U)^{\#}\rightarrow 0$$
as before. Now if $\lambda\in \Lambda_b
=H^{n-1}_c(X_b^{\#},\boldz)=H_1(X_b^{\#},\boldz)\subseteq \T_{B,b}^*$ via
the map $\psi'$, then
in local coordinates $(y_1,\ldots,y_n)$ on $B$, $\lambda
=\sum\lambda_i dy_i$. Now
$$\eqalign{-\int_{\lambda}\iota(\partial/\partial y_i)\omega
&=-\int_{(0,\ldots,0)}^{(\lambda_1,\ldots,\lambda_n)}
\iota(\partial/\partial y_i)\pi^*\omega\cr
&=
\int_{(0,\ldots,0)}^{(\lambda_1,\ldots,\lambda_n)}
dx_i\cr
&=\lambda_i\cr}$$
so the two maps $\psi,\psi':
R^{n-1}_cf_*^{\#}\boldz|_U\rightarrow\T_U^*$ coincide.
Since $\Lambda\subseteq\T_U^*$ is Lagrangian, we conclude that the
image of $R^{n-1}_cf_*^{\#}\boldz$ in $\T_B^*$ under $\psi$ is
Lagrangian.
Thus there is an exact sequence
$$\exact{R^{n-1}_cf^{\#}_*\boldz}{\T_B^*}{J^{\#}}$$
defining $J^{\#}$ in which
$J^{\#}$ inherits the standard symplectic form from $\T_B^*$. 
Since
$f:X^{\#}\rightarrow B$ is locally isomorphic to $J^{\#}\rightarrow B$,
a standard argument [12] shows that one can obtain $X^{\#}\rightarrow B$ from
$J^{\#}\rightarrow B$ by regluing using a \v Cech 1-cocycle 
$\{(U_i,\sigma_i)\}$ where $\sigma_i$ is a Lagrangian section of 
$J^{\#}\rightarrow B$ over $U_i$. We call $j:J^{\#}\rightarrow B$
the {\it Jacobian fibration} of $f:X^{\#}\rightarrow B$, in analogy
with the theory of elliptic curves. This gives a one-to-one correspondence
between the group $H^1(B,\Lambda(J^{\#}))$, where $\Lambda(J^{\#})$
is the sheaf of Lagrangian sections of $J^{\#}$, and the set
$$\{\hbox{$f:Y^{\#}\rightarrow B$ a Lagrangian fibration with local
section and Jacobian $j:J^{\#}\rightarrow B$}\}/\cong.$$
In fact, this can be extended to the compactifications. We phrase
this more generally for Lagrangian fibrations.

\proclaim Theorem 2.6. Let $f:X\rightarrow B$ be a proper Lagrangian 
fibration with connected fibres, with a local section everywhere.
Then there is a symplectic manifold $J$, called the Jacobian of $f$,
and a (proper) Lagrangian fibration $j:J\rightarrow B$ which is locally
isomorphic to $f:X\rightarrow B$, and which has a Lagrangian
section. Furthermore, there is a one-to-one correspondence between
the sets 
$$\{\hbox{$f:Y\rightarrow B$ a Lagrangian fibration
with local section and Jacobian $j:J\rightarrow B$}\}/
\cong$$
and $H^1(B,\Lambda(J^{\#}))$.

Proof. To construct $J$, choose an open covering $\{U_i\}$ such that
$f^{-1}(U_i)\rightarrow U_i$ has a Lagrangian section $\sigma_i$
for each $i$. Then $f^{-1}(U_i\cap U_j)$ has a symplectomorphism we will write
as $T_{i,j}$ obtained by treating $\sigma_i$ as the zero-section
and then translating by $\sigma_j$, so that $T_{i,j}$ takes
$\sigma_i$ into $\sigma_j$. We construct $J$ by identifying
$f^{-1}(U_i)$ and $f^{-1}(U_j)$ along $f^{-1}(U_i\cap U_j)$,
using $T_{i,j}$ to identify $f^{-1}(U_i\cap U_j)\subseteq f^{-1}(U_i)$
and $f^{-1}(U_i\cap U_j)\subseteq f^{-1}(U_j)$. These identifications
are compatible over $U_i\cap U_j\cap U_k$ because $T_{j,k}\circ T_{i,j}
=T_{i,k}$. Thus we obtain a Lagrangian fibration $j:J\rightarrow B$ with a 
section, as desired. The usual regluing construction gives the 1-1 
correspondence.
$\bullet$

For $j:J\rightarrow B$ special Lagrangian and integral,
one computes $H^1(B,\Lambda(J^{\#}))$ by using the exact sequence
$$\exact{R^{n-1}_cj_*^{\#}\boldz}{\Lambda(\T_B^*)}{\Lambda(J^{\#})}.$$
$\Lambda(\T_B^*)$ is just the kernel of exterior differentiation
acting on $\T_B^*$, so $H^i(B,\Lambda(\T_B^*))=H^{i+1}(B,{\bf R})$ for
$i\ge 1$. From this we obtain the sequence
$$H^2(B,{\bf R})\rightarrow H^1(B,\Lambda(J^{\#}))\rightarrow H^2(B,R^{n-1}_c
j_*^{\#}\boldz)\rightarrow H^3(B,{\bf R}).$$
In any dimension, if $H^2(B,{\bf R})=0$, then 
$H^2(B,R_c^{n-1}j_*^{\#}\boldz)_{tors}\subseteq H^1(B,\Lambda(J^{\#}))$.
Duistermaat [12] observed that if an element $\alpha\in H^1(B,\Lambda(J^{\#}))$
comes from an element $[\alpha']\in H^2(B,{\bf R})$, then the corresponding
$f:X\rightarrow B$ can be obtained by choosing a $2$-form $\alpha'$
on $B$ representing $[\alpha']$, and taking $X=J$ with symplectic form
$\omega+f^*\alpha'$. Any two choices of $\alpha'$ can be related
by translation by a section.

{\it Example 2.7.} If $n=3$ and $H^2(B,{\bf R})=0$, then we have
$$H^2(B,R^2_cf^{\#}_*\boldz)_{tors}\subseteq H^1(B,\Lambda(J^{\#}))
\subseteq H^2(B,R^2_cf_*^{\#}\boldz).$$
If $n=2$ and $J$ is a K3 surface, then $H^2(B,R^1_cf_*^{\#}\boldz)=0$
and so there is a sequence
$$H^1(B,R^1_cf^{\#}_*\boldz)\rightarrow H^2(B,{\bf R})\rightarrow
H^1(B,\Lambda(J^{\#}))\rightarrow 0.$$

{\it Remark 2.8.} If $U\subseteq B-\Delta$ is a 
simply connected set, then there is no monodromy in the local system
$(R^{n-1}f_*\boldz)|_U\subseteq\T_U^*$. Thus, if $\lambda_1,\ldots,\lambda_n$
are sections of $\T_U^*$ generating $(R^{n-1}f_*\boldz)|_U$, the fact
that $d\lambda_i=0$ shows there are functions $u_i$ such that
$du_i=\lambda_i$ on $U$. The $u_1,\ldots,u_n$ form local coordinates
on $U$, since $du_1,\ldots,du_n$ are independent. This is the standard
construction of action coordinates on $U$. We will say that $u_1,\ldots,
u_n$ are {\it action coordinates} for $f:X\rightarrow B$ on $U$. Canonical
coordinates $u_1,\ldots,u_n,x_1,\ldots,x_n$ are called action-angle 
coordinates.
These are also the coordinates that Hitchin introduces after [17], Prop. 1.
The advantage of working in this coordinate system is that the
periods are now just the constant periods $du_1,\ldots,du_n$. 
\bigskip

{\hd \S 3. Simplicity and the Leray spectral sequence revisited.}

Recall from [14] that a special Lagrangian $T^n$-fibration
$f:X\rightarrow B$ was said to be {\it $G$-simple} if
$$i_*R^pf_{0*}G=R^pf_*G$$
for all $p$, where $i:B-\Delta\rightarrow B$ is the inclusion,
$f_0=f|_{f^{-1}(B-\Delta)}$, and $G$ is an abelian group.
This condition was crucial for getting a handle on the topology
of $X$ and the relationship between the topology of $X$ and its dual.
In [14], we only made use of ${\bf Q}$-simplicity, while here, we
will go further. We are interested in a broader range
of groups $G$. In particular, $G={\bf R},{\bf Q},\boldz,\boldz/m\boldz$,
or ${\bf R}/\boldz$ will be of relevance for us. Clearly
$\boldz$-simplicity implies ${\bf Q}$-simplicity or ${\bf R}$-simplicity,
but $\boldz/m\boldz$ or ${\bf R}/\boldz$-simplicity provides extra information
about monodromy modulo $m$ which may be valuable in trying to classify
possible monodromy transformations about the discriminant locus.

If $f$ is integral, then $f$ has connected fibres, and $G
=f_*G=i_*f_{0*}G$.
We have seen in \S 2 that if Assumption 1.8 holds then
$R^nf_*\boldz=R^n_cf_*^{\#}\boldz=\boldz$,
and thus by the universal coefficient theorem, $R^nf_*G=G$. Thus
for integral fibrations satisfying Assumption 1.8,
the simplicity condition holds for $p=0$ and $n$ and any abelian group $G$.

We now prove it holds for $p=1$, under the additional Assumption 1.7$'$.

\proclaim Definition 3.1.
We say a point $b\in B$ is a {\it rank $k$} point if 
$$k=\min_{x\in X_b} \rank f_*:\T_{X,x}\rightarrow \T_{B,b}.$$

In particular $b$ is rank $n$ if and only if $f$ is smooth over $b$.

We first comment on when Assumption 1.7$'$ might hold.

\proclaim Lemma 3.2. If $b$ is a rank $k$ point, $k\not=0$,
then Assumption 1.7 implies Assumption $1.7'$ at $b$ if $X_b$
contains only a finite number of orbits of the action of $\T_{B,b}^*$.
Note Assumption 1.7 automatically implies Assumption $1.7'$ at a rank 0
point.

Proof. Consider the set $\L'=\{\im f_*:\T_{X,x}\rightarrow\T_{B,b}|x\in X_b\}$
of subspaces of $\T_{B,b}$. Note that if $x,y\in X_b$ with $T_{\sigma}(x)
=y$ for some $\sigma$, then $f_{*,x}=f_{*,y}\circ (T_{\sigma})_*$ where
$f_{*,x}$ denotes the pushforward of tangent vectors at $x$. Thus
$\im f_{*,x}=\im f_{*,y}$. By the assumption that $X_b$ is a union
of a finite number of orbits, we see then that $\L'$ is finite,
and $\min_{L\in\L'} \dim L=k$. Because $\L'$ is finite we can choose a 
subspace $T\subseteq \T_{B,b}$ of dimension $n-k$ in the set
$\L$ given by Assumption 1.7 which intersects
every element of $\L'$ in the expected dimension. It then follows
that $\dim f_{*,x}^{-1}(T)=2n-k$ for all $x\in X_b$. Now if $B'$ is
taken to be a submanifold of $B$ with tangent space $T$ at $b$, we see
that the implicit function theorem implies that $f^{-1}(B')$ is
a manifold in a neighborhood of $X_b$. $\bullet$

This makes Assumption 1.7$'$ appear quite reasonable, at least in the 
three-dimensional case, as the finiteness of the number of orbits
for rank 1 fibres would follow from finiteness of the $\H^1$-measure
of the singular locus. (As mentioned in \S 1 however,
such a result is not actually known.) 

\proclaim Theorem 3.3. If $f:X\rightarrow B$ is an integral
special Lagrangian fibration satisfying Assumption $1.7'$, 
then $i_*R^1f_{0*}G=R^1f_*G$.

Proof. Let
$$\Delta_0\subseteq\Delta_1\subseteq\cdots
\subseteq \Delta_{n-2}=\Delta\subseteq B$$ be the stratification 
of $\Delta$ given in \S 1.
Let $i_k:B-\Delta_k\rightarrow B-\Delta_{k-1}$ be the inclusion.
We will show using descending induction that
$$i_{k*}(R^1f_*G)|_{B-\Delta_k}=(R^1f_*G)|_{B-\Delta_{k-1}}$$
for each $k$. One always has a functorial map
$$(R^1f_*G)|_{B-\Delta_{k-1}}\rightarrow i_{k*}i_k^*(R^1f_*G)|_{B-
\Delta_{k-1}}=i_{k*}(R^1f_*G)|_{B-\Delta_k},$$
and we just need to show this is an isomorphism on the level of stalks
at each point $b\in\Delta_k$, $b\not\in\Delta_{k-1}$.

We can choose a $B'$ through the point $b$ using Assumption 1.7$'$,
so that $\Delta_k\cap B'=\{b\}$ and $X'=f^{-1}(B')$ is a manifold.
This gives a diagram
$$\matrix{B'-\{b\}&\mapright{i'}&B'\cr
\mapdown{j'}&&\mapdown{j}\cr
B-\Delta_k&\mapright{i_k}&B-\Delta_{k-1}\cr}$$
$i',j',j$ the inclusions. Let $f':X'\rightarrow B'$
be the restriction of $f$. Then $j^*R^1f_*G=R^1f'_*G$
and in particular the stalks of $R^1f_*G$ and
$R^1f'_*G$ at $b$ are the same. On the other hand, there is
a natural map $(i_{k*}i_k^*R^1f_*G)_b\rightarrow (i_*'i'^*R^1f'_*G)_b$.
Indeed, an element of $(i_{k*}i_k^*R^1f_*G)_b$ represented by a germ
$(U,\alpha)$, $\alpha\in\Gamma(U-\Delta_k,R^1f_*G)$, is mapped
to $(U\cap B',\alpha|_{(U\cap B')-\{b\}})$, $\alpha|_{(U\cap B')-\{b\}}
\in\Gamma((U\cap B')-\{b\},R^1f'_*G)$. This map is in fact injective.
Indeed, by descending induction, $(R^1f_*G)|_{B-\Delta_k}$
has no sections over any open subset of $B-\Delta_k$ supported on
a proper closed subset. Thus the restriction maps of the sheaf 
$(R^1f_*G)|_{B-\Delta_k}$ are injective, and it follows that the map
$(i_{k*}i_k^*R^1f_*G)_b\rightarrow (i_*'i'^*R^1f'_*G)_b$ is injective.
We then have a diagram
$$\matrix{(R^1f_*G)_b&\mapright{\cong}&(R^1f'_*G)_b\cr
\mapdown{}&&\mapdown{}\cr
(i_{k*}i_k^*R^1f_*G)_b&\hookrightarrow&(i_*'i'^*R^1f_*G)_b\cr}$$
so $(R^1f_*G)_b\rightarrow (i_{k*}i_k^*R^1f_*G)_b$ is an isomorphism
if $(R^1f'_*G)_b\rightarrow (i_*'i'^*R^1f_*'G)_b$ is.
Thus we need to show that
$$\lim_{\rightarrow\atop b\in U\subseteq B'} H^1(f'^{-1}(U),G)
\rightarrow \lim_{\rightarrow\atop b\in U\subseteq B'} \Gamma
(U-\{b\},R^1f'_*G)$$
is an isomorphism. We can take the direct limit over contractible $U$.

We will need to know what $H^2_{X_b}(X',G)$ is. By the universal
coefficient theorem, there is an exact sequence
$$\exact{H^i_{X_b}(X',\boldz)\otimes_{\boldz} G}{H^i_{X_b}(X',G)}
{\tor^1_{\boldz}(H^{i+1}_{X_b}(X',\boldz),G)}.$$
Now by a suitable form of Poincar\'e duality ([10], V 9.3),
$H^i_{X_b}(X',\boldz)\cong H_{2n-k-i}(X_b,\boldz)$, where the latter
group is Borel-Moore homology. This is computed via the exact sequence
([10] V, \S 3, (9))
$$\exact{\ext^1(H_c^{2n-k-i+1}(X_b,\boldz),\boldz)}{H_{2n-k-i}(X_b,\boldz)}
{\hom(H_c^{2n-k-i}(X_b,\boldz),\boldz)}.$$
Given that $H^n_c(X_b,\boldz)=\boldz$ and $H^{n-1}_c(X_b,\boldz)=
H^{n-1}_c(X_b^{\#},\boldz)$ is free, we see that
$$H^2_{X_b}(X',\boldz)=
\cases{\boldz&if $k=n-2$;\cr
0&if $k<n-2$,\cr}$$
and $H^3_{X_b}(X',\boldz)$ is free, so that
$$H^2_{X_b}(X',G)=
\cases{G&if $k=n-2$;\cr
0&if $k<n-2$.\cr}$$
Of course $H^i_{X_b}(X',G)=0$ for $i<2$. 

There are two cases:

{\it Case 1:} $k<n-2$. 
For $b\in U\subseteq B'$,
$H^1(f'^{-1}(U),G)\cong H^1(f'^{-1}(U-\{b\}),G)$,
by the relative cohomology
long exact sequence. On the other hand, the Leray spectral sequence
for $f':X^*=f'^{-1}(U-\{b\})\rightarrow B^*=U-\{b\}$ yields the exact sequence
$$0=H^1(B^*,G)\rightarrow H^1(X^*,G)\rightarrow 
H^0(B^*,R^1f'_*G)\rightarrow H^2(B^*,G)
\rightarrow H^2(X^*,G).$$
This last map is injective as $X^*\rightarrow B^*$ can be assumed to have
a section. Thus we obtain the isomorphism
$$H^1(f'^{-1}(U),G)\rightarrow H^0(B^*,R^1f'_*G),$$
and taking direct limts gives the desired isomorphism.

{\it Case 2.} $k=n-2$. In this case,
the relative cohomology exact sequence gives for $b\in U\subseteq B'$,
$X^*=f'^{-1}(U-\{b\}), B^*=U-\{b\}$,
$$0\rightarrow H^1(f'^{-1}(U),G)\rightarrow H^1(X^*,G)
\rightarrow H^2_{X_b}(f'^{-1}(U),G)\rightarrow H^2(f'^{-1}(U),G).$$
In this case $H^2_{X_b}(f'^{-1}(U),G)=G$.
We have a commutative square
$$\matrix{H^1(B^*,G)&\mapright{\cong}& H^2_{\{b\}}(U,G)\cr
\mapdown{}&&\mapdown{\cong}\cr
H^1(X^*,G)&\mapright{}&H^2_{X_b}(f'^{-1}(U),G)\cr}$$
showing the map $H^1(X^*,G)\rightarrow H^2_{X_b}(f'^{-1}(U),G)$
is surjective, yielding
$$\exact{H^1(f'^{-1}(U),G)}{H^1(X^*,G)}{G}.$$
On the other hand, the Leray spectral sequence for $f':X^*\rightarrow B^*$
gives
$$\exact{H^1(B^*,G)\cong G}{H^1(X^*,G)}{H^0(B^*,R^1f_*'G)}.
$$
Putting these two sequences together one finds that the map
$$H^1(f'^{-1}(U),G)\rightarrow H^0(B^*, R^1f_*'G)$$
is an isomorphism, and hence
$$\lim_{\rightarrow} H^1(f'^{-1}(U),G)\rightarrow \lim_{\rightarrow}
\Gamma(U-\{b\},R^1f'_*G)$$
is an isomorphism. $\bullet$

We next try to understand $R^{n-1}f_*\boldz$. In any event, in \S 2 we
have seen that if $f:X\rightarrow B$ is integral and satisfies Assumption
1.8 then
$R^{n-1}_cf_*^{\#}\boldz\cong
R^{n-1}f_*\boldz$ and $\HH^0_{\Delta}(B,R^{n-1}_cf_*^{\#}\boldz)=0$.
Since there is an exact sequence
$$0\rightarrow \HH^0_{\Delta}(B,R^{n-1}f_*\boldz)\rightarrow
R^{n-1}f_*\boldz\rightarrow i_*R^{n-1}f_{0*}\boldz
\rightarrow \HH^1_{\Delta}(B,R^{n-1}f_*\boldz)\rightarrow 0,$$
we have already shown at least that the natural map $R^{n-1}f_*\boldz
\rightarrow i_*R^{n-1}f_{0*}\boldz$ is injective.
However, to show surjectivity, we need a more delicate understanding
of the inductive structure of the singular locus. If $b\in B$ is
a rank $k$ point, it would, for example, be sufficient to show that
there is a fixed-point-free
Hamiltonian $T^k$ action in a neighborhood of $X_b$ in order to 
achieve a sufficiently strong inductive description of the singular fibres. 
However failing to prove such a result, we will make an ad hoc argument
in the $n=3$ case. Nevertheless, in any dimension we have

\proclaim Lemma 3.4. If $f$ is integral and satisfies Assumptions
1.7 and 1.8, $\Delta_0$ the set of
rank 0 points of $B$, $i_0:B-\Delta_0\hookrightarrow B$ the inclusion, then
$$i_{0*}i_0^*R^{n-1}f_*\boldz=R^{n-1}f_*\boldz.$$

Proof. If $U$ is a contractible open neighborhood of $b\in\Delta_0$, we have
an exact sequence
$$\exact{H^{n-1}(f^{-1}(U),\boldz)}{H^{n-1}(f^{-1}(U-\{b\}),\boldz)}
{H^n_{X_b}(f^{-1}(U),\boldz)\cong\boldz}.$$
In addition, the map
$$\lim_{\rightarrow\atop b\in U} H^{n-1}(f^{-1}(U),\boldz)
\rightarrow\lim_{\rightarrow\atop b\in U}\Gamma(U-\{b\},R^{n-1}f_*\boldz)$$
is injective.
We just need to show this map is surjective.
From the Leray spectral sequence for $f:X^*=f^{-1}(U-\{b\})\rightarrow B^*
=U-\{b\}$,
we obtain a map $\varphi:H^{n-1}(X^*,\boldz)
\rightarrow\Gamma(B^*,R^{n-1}f_*\boldz)$. We first show this map is
surjective. Indeed, given a section $\sigma\in \Gamma(B^*,R^{n-1}f_*\boldz)
\subseteq\Gamma(B^*,\T_B^*)$, let $M\subseteq X^*$ be the
circle bundle over $B^*$ whose fibre at $b\in B^*$ is
${\bf R}\sigma(b)/\boldz\sigma(b)\subseteq X_b^{\#}=\T_{B,b}^*/(R^{n-1}f_*
\boldz)_b$. Let $[M]\in H^{n-1}(X^*,\boldz)$ be the image of 
$1\in H^{n-1}_M(X^*,\boldz)\cong H^0(M,\boldz)$ in $H^{n-1}(X^*,\boldz)$.
It is then clear that $\varphi([M])=\sigma$. Thus $\varphi$ is surjective.

We now have a diagram
$$\matrix{&&0&&&&&&\cr
&&\mapdown{}&&&&&&\cr
0&\mapright{}&\lim_{\rightarrow}
H^{n-1}(f^{-1}(U),\boldz)&\mapright{}&\lim_{\rightarrow} H^{n-1}(X^*,\boldz)&
\mapright{}&\boldz&\mapright{}&0\cr
&&\mapdown{\varphi'}&&\mapdown{\varphi}&&\mapdown{\varphi''}&&\cr
0&\mapright{}&\im\varphi&\mapright{\alpha}&
\lim_{\rightarrow} 
\Gamma(B^*,R^{n-1}f_*\boldz)&\mapright{}&\coker\alpha&\mapright{}&0\cr
&&\mapdown{}&&\mapdown{}&&&&\cr
&&0&&0&&&&\cr}$$
and we wish to show that $\coker\alpha=0$, so that
$\alpha\circ\varphi'$ is surjective as desired. By the snake lemma,
$\ker\varphi\cong\ker\varphi''$. By the Leray spectral sequence
for $f:X^*\rightarrow B^*$, the image of $H^{n-1}(B^*,\boldz)\rightarrow
H^{n-1}(X^*,\boldz)$ is contained in $\ker\varphi$. On the other hand, by
the commutativity of
$$\matrix{H^{n-1}(B^*,\boldz)&\mapright{\cong}& H^n_{\{b\}}(U,\boldz)\cr
\mapdown{}&&\mapdown{\cong}\cr
H^{n-1}(X^*,\boldz)&\mapright{}&H^n_{X_b}(f^{-1}(U),\boldz)\cr}$$
it is then clear that $\ker\varphi$ surjects onto $H^n_{X_b}(f^{-1}(U),\boldz)
=\boldz$, so $\ker\varphi''=\boldz$ and $\coker\alpha=0$. $\bullet$

\proclaim Theorem 3.5. If $\dim X=3$, $f:X\rightarrow B$ an integral
special Lagrangian fibration satisfying 
Assumptions 1.7 and 1.8, and if the fibres $X_b$ for rank 1 points
$b$ are a union of a finite number of $\T_{B,b}^*$-orbits, 
then $f$ is $\boldz$-simple.

Proof. The hypothesis of this theorem implies the hypothesis of
Theorem 3.3, so the simplicity condition holds for
$p=0,1$ and $3$.
In the notation of the proof of Theorem 3.3, we need to show
$i_{1*}i_1^*(R^2f_*\boldz)|_{B-\Delta_0}=(R^2f_*\boldz)|_{B-\Delta_0}$, 
as Lemma 3.4 then allows
us to complete the proof of simplicity. Choose $b\in\Delta_1$,
$b\not\in\Delta_0$, and choose a 2-dimensional disk $B'$ passing through
$b$ as in the proof of Theorem 3.3, $f':X'=f^{-1}(B')\rightarrow B'$.
As in that proof, we just need to show that
$$H^2(X_b,\boldz)\cong \lim_{\rightarrow\atop
b\in U\subseteq B'} H^2(f'^{-1}(U),\boldz)
\rightarrow \lim_{\rightarrow\atop
b\in U\subseteq B'} \Gamma(U-\{b\},R^2f_*'\boldz)$$
is an isomorphism, and we already know this map is
injective. Note also that the cokernel of this map is torsion-free:
since $R^2f'_*\boldz\subseteq \T_B^*|_{B'}$, if an integer
multiple of a section of $R^2f'_*\boldz$ over $U-\{b\}$
extends to a section over $U$, the section itself extends.
Now $B^*=B'-\{b\}$ is a punctured disk,
and hence $(R^2f_*'\boldz)|_{B^*}$ is a local system determined
by a single monodromy transformation $T$. By Poincar\'e duality,
$(R^1f_*'\boldz)|_{B^*}$ has monodromy ${}^tT$. Since
$\Gamma(B^*,R^2f_*'\boldz)=\ker(T-I)$,
we see that $\rank \Gamma(B^*,R^2f_*'\boldz)
=\rank \Gamma(B^*,R^1f_*'\boldz)=\rank H^1(X_b,\boldz)$
by simplicity for $p=1$.
Thus it will be sufficient to show that $\rank
H^1(X_b,\boldz)=\rank H^2(X_b,\boldz)$ to show the above map
is an isomorphism.

Since $b$ is a rank 1 point, the singular locus of $X_b$ is a union
of circles, each being a closed orbit of the action of $\T_{B,b}^*$ on
$X_b$. Since there are assumed to be only a finite number
of such orbits, each $S^1$ is a connected component of $Sing(X_b)$.
But since $Sing(X_b)$ is connected by Corollary 2.5,
$Z:=Sing(X_b)\cong S^1$. Now we also have exact sequences
$$\exact{H^i_c(X_b-\boldz,\boldz)}{H^i(X_b,\boldz)}{H^i(Z,\boldz)}$$
for all $i$, the exactness for $i=1$ shown in the proof of Corollary 2.5
and for $i=2$ shown in the proof of Theorem 2.4. From this we see
that if $H^2_c(X_b-Z,\boldz)=\boldz$, then
$H^1(X_b,\boldz)=H^2(X_b,\boldz)=\boldz$, while if $H^2_c(X_b-Z,\boldz)
=\boldz^2$, then $H^1(X_b,\boldz)=H^2(X_b,\boldz)=\boldz^2$,
completing the proof in these cases.

Finally, suppose $H^2_c(X_b-Z,\boldz)=0$. 
Then $H^2_c(X_b,\boldz)=0$ and it follows as in the proof of Theorem 3.3
that $H^3_{X_b}(f^{-1}(U),\boldz)=0$, so one obtains from the relative
cohomology sequence a surjection
$$H^2(f'^{-1}(U),\boldz)\rightarrow H^2(f'^{-1}(U-\{b\}),\boldz)
\rightarrow 0.$$
The argument of the proof of Lemma 3.4 shows that $H^2(f'^{-1}(U-\{b\}),
\boldz)\rightarrow \Gamma(U-\{b\},R^2f_*'\boldz)$ is surjective,
and hence so is $H^2(f'^{-1}(U),\boldz)\rightarrow\Gamma(U-\{b\},
R^2f_*'\boldz)$. Taking direct limits, one concludes that
$H^2(X_b,\boldz)=\ker(T-I)$ as desired, showing simplicity. 
But notice in fact this
case can't occur, since as $H^2(X_b,\boldz)=0$, we also have
$H^1(X_b,\boldz)=\ker({}^tT-I)=0$, 
contradicting $H^1(X_b,\boldz)=H^1(Z,\boldz)=\boldz$. $\bullet$

{\it Remark 3.6.} 
If $f:X\rightarrow B$ is a $\boldz$-simple special Lagrangian
$T^3$ fibration, we obtain some restrictions on the cohomology 
of a singular fibre $X_b$. Clearly $H^0(X_b,\boldz)=H^3(X_b,\boldz)=\boldz$,
so if $b_i=\rank H^i(X_b,\boldz)$, we will say for the duration of this 
remark that $X_b$ is of type $(b_1,b_2)$. Clearly $b_1,b_2\le 3$,
and if $b_2=3$, then $X_b$ is non-singular. If $b_1=3$, then $b_2=3$,
since $\bigwedge^2 H^1(X_b,\boldz)\subseteq H^2(X_b,\boldz)$, and so $X_b$
is non-singular. This also shows that if $b_1=2$ then $b_2\ge 1$, and a similar
argument shows that if $b_2=2$, then $b_1\ge 1$. Thus the
possible values for $(b_1,b_2)$ are
$(2,2)$, $(2,1)$, $(1,2)$, $(1,1)$, $(0,1)$, $(1,0)$ or $(0,0)$. We describe
the probable topology of an integral singular fibre with each of the
above possible cohomology groups.
\item{$(2,2)$} Such cohomology is realised by a fibre of the form
$I_1\times S^1$, where $I_1$ denotes a Kodaira type $I_1$-fibre.
\item{$(1,2)$} $S^1\times T^2/\{pt\}\times T^2$. This was seen in the example
in \S 1 of [14].
\item{$(2,1)$} This was described in [14], Remark 1.4. Identify
$T^3$ with the solid cube $[0,1]^3$ with opposite sides identified. Then
take 
$T^3/\sim$, where $(x_1,x_2,x_3)\sim(x_1',x_2',x_3')$ if and only if
$(x_1,x_2,x_3)=(x_1',x_2',x_3')$ or $(x_1,x_2)=(x_1',x_2')\in 
\partial([0,1])^2$. This fibre is singular along a figure eight.
\item{$(1,1)$} There are two possibilities here. The fibre
could be $II\times S^1$, where $II$ denotes a
Kodaira type II fibre, or it could be
$T^3/\sim$, where $(x_1,x_2,x_3)\sim(x_1',x_2',x_3')$ if and only if
$(x_1,x_2,x_3)=(x_1',x_2',x_3')$ or $(x_1,x_2)=(x_1',x_2')\in 
\{0,1\}\times [0,1]$ or $(x_1,x_2),(x_1',x_2')\in [0,1]\times\{0,1\}$.
(This is the same as contracting one loop of the singular figure eight in the
$(2,1)$ case to a point).
\item{$(0,1)$} 
$T^3/\sim$, where $(x_1,x_2,x_3)\sim(x_1',x_2',x_3')$ if and only if
$(x_1,x_2,x_3)=(x_1',x_2',x_3')$ or $(x_1,x_2),(x_1',x_2')\in \partial
[0,1]^2$. (This is equivalent to contracting the singular figure eight of the
$(1,2)$ case to a point).
\item{$(1,0)$}
$T^3/\sim$, where $(x_1,x_2,x_3)\sim(x_1',x_2',x_3')$ if and only if
$(x_1,x_2,x_3)=(x_1',x_2',x_3')$ or $(x_1,x_2),(x_1',x_2')\in 
\partial[0,1]^2$ and $x_3=x_3'$ or $x_3,x_3'\in \{0,1\}$. 
\item{$(0,0)$} 
$T^3/\sim$, where $(x_1,x_2,x_3)\sim(x_1',x_2',x_3')$ if and only if
$(x_1,x_2,x_3)=(x_1',x_2',x_3')$ or $(x_1,x_2,x_3),(x_1',x_2',x_3')\in
\partial[0,1]^3$.
(We are contracting the boundary
of $[0,1]^3$ to a point. This is topologically a sphere).

One notes that for each fibre of type $(m,n)$, there is a fibre of
type $(n,m)$ which should then be its dual. (In particular, the fibres
of type $(2,2)$ and $(1,1)$ should be self-dual).

I cannot prove yet that this provides a complete classfication of integral
three-dimensional singular fibres, but it seems to be a reasonable conjecture.
In addition, these types of examples extend to higher dimensions, and one
finds a much wider range of possible topologies, which nonetheless
exhibit the desired duality.

\bigskip
Having proven $\boldz$-simplicity in some cases for special Lagrangian
fibrations, we wish now to return to a more careful
study of the Leray spectral sequence for special Lagrangian fibrations,
with special consideration of the role torsion plays. First we make
some observations on the Leray spectral sequence in any dimension.

\proclaim Lemma 3.7. If $f:X\rightarrow B$ is a $\boldz$-simple special
Lagrangian $T^n$-fibration with a section, then in the Leray spectral
sequence for $f$, $E^2_{1,1}=E^{\infty}_{1,1}$ and 
$E^2_{1,n-1}=E^{\infty}_{1,n-1}$. In addition, the Leray filtration
yields a surjection $H^n(X,\boldz)_{tors}\rightarrow (E^2_{1,n-1})_{tors}$.

Proof. The only possible non-zero differential to or from $E^2_{1,1}$ is
$d_2:E^2_{1,1}\rightarrow E^2_{3,0}=H^3(B,\boldz)$. But since $f$ has
a section, the map $H^3(B,\boldz)\rightarrow H^3(X,\boldz)$ is injective,
and thus $d_2=0$. Thus $E^2_{1,1}=E^{\infty}_{1,1}$.

Next, $E^2_{1,n-1}=H^1(B,R^{n-1}f_*\boldz)$, and recall from [14] that
that $H^1(B,R^{n-1}f_*\boldz)$ is the group
of sections of $f$ modulo homotopy, with a fixed section,
say $\sigma_0$, the zero section. Then for any section $\sigma$,
the cohomology class $[\sigma]-[\sigma_0]\in H^n(X,\boldz)$ and
the element $[\sigma]\in H^1(B,R^{n-1}f_*\boldz)$ representing the section
coincide up to sign 
in $E^{\infty}_{1,n-1}$ by [14], Theorem 4.1 (which
holds with $\boldz$ coefficients if $f$ is $\boldz$-simple). 
Thus $E^{\infty}_{1,n-1}=E^2_{1,n-1}$, and if $H^n(X,\boldz)
=F^0\supseteq F^1\supseteq\cdots$ is the Leray filtration on
$H^n(X,\boldz)$, then $F^1/F^2\cong E^2_{1,n-1}$. Since $F^0/F^1\cong
\boldz$, $H^n(X,\boldz)_{tors}=F^1_{tors}$ and thus there is a map 
$H^n(X,\boldz)_{tors}\rightarrow (E^2_{1,n-1})_{tors}$. To see this map is 
surjective,
suppose $\sigma$ is a torsion section of $f$. Then $\sigma$ must
be disjoint from $\sigma_0$. Indeed, if $x\in \sigma\cap \sigma_0$,
let $U\subseteq B$ be a small open neighborhood of $f(x)\in B$
in which $\sigma$ is represented by a section $\tilde\sigma\in
\Gamma(U,\T_B^*)$, such that $\tilde\sigma(f(x))=0$. Then if $m$ is
the order of the torsion section $\sigma$, $m\tilde\sigma$ is
a non-zero section of $R_c^{n-1}f^{\#}_*\boldz$ which is zero 
for at least one point, which is impossible. Thus $\sigma$ and
$\sigma_0$ are disjoint.

By [14], Theorem 4.1,
$T_{\sigma}^*:H^*(X,{\bf Q})\rightarrow H^*(X,{\bf Q})$ is a unipotent
operator, but on the other hand $T_{\sigma}^m=I$ since $\sigma$ is $m$-torsion.
Thus $T_{\sigma}^*=I$, so $[\sigma]=[\sigma_0]$ in $H^n(X,{\bf Q})$ and
$[\sigma]-[\sigma_0]$ is in fact a torsion element of $H^n(X,\boldz)$. This
shows that the map $H^n(X,\boldz)_{tors}\rightarrow (E^2_{1,n-1})_{tors}$
is surjective.
$\bullet$

Even if in general the Leray spectral sequence for $f$
does not degenerate, the above result might be sufficient for many
applications; as we will see in later sections, $H^1(B,R^1f_*{\bf R})$ and
$H^1(B,R^{n-1}f_*{\bf R})$ play important roles in mirror symmetry.

Note that if $\boldz$-simplicity fails because $f$ has reducible fibres,
we expect the second part of the above result to fail.
\bigskip

We now focus on the three dimensional case. 

\proclaim Proposition 3.8. 
Let $f:X\rightarrow B$ and $\check f:\check X\rightarrow B$
be dual $\boldz$-simple special Lagrangian $T^3$
fibrations, and suppose that $H^1(B,\boldz)=0$.
Then $H^1(X,\boldz)=0$ if and only if $H^1(\check X,\boldz)=0$.

Proof. Since $H^1(B,\boldz)=0$, $H^2(B,\boldz)$ is torsion, so
$$rank(H^1(X,\boldz))=rank(H^0(B,R^1f_*\boldz)).$$ But since
$X$ is K\"ahler, the first betti number of $X$ is even, so
if $H^1(X,\boldz)\not=0$ then $rank(H^0(B,R^1f_*\boldz))\ge 2$.
The wedge of two independent sections of $R^1f_*\boldz$ yields
a section of $R^2f_*\boldz$, and $R^2f_*\boldz\cong R^1\check f_*\boldz$,
so $H^0(B,R^1\check f_*\boldz)\not=0$, hence $H^1(\check X,\boldz)\not=0$.
Repeating the same argument interchanging $X$ and $\check X$ gives the
result. $\bullet$

\proclaim Theorem 3.9. Let $f:X\rightarrow B$, $\check f:\check X\rightarrow B$
be dual $\boldz$-simple
special Lagrangian $T^3$-fibrations with
sections, and assume $H^1(X,\boldz)=0$.
Then the Leray spectral sequences for $f$ and $\check f$ with coefficients
in $\boldz$ degenerate at the $E_2$-term, and
$$rank_{\boldz} H^i(B,R^jf_*\boldz)=rank_{\boldz} H^{3-i}(B,R^{3-j}f_*\boldz).$$
If in addition $f$ and $\check f$ are ${\bf R}/\boldz$-simple, then
$$Tors(H^i(B,R^jf_*\boldz))\cong Tors(H^{4-i}(B,R^{3-j}f_*\boldz)).$$

Note that the additional assumption of ${\bf R}/\boldz$-simplicity is not
a particularly strong one. Indeed, we have seen that given
suitable regularity hypotheses, we have only failed to show the
$G$-simplicity condition for $p=2$. But $\boldz$-simplicity implies
$R^2f_*\boldz\cong R^1\check f_*\boldz$ and $R^2f_*{\bf R}\cong
R^1\check f_*{\bf R}$, whence $R^2f_*{\bf R}/\boldz\cong
R^1\check f_*{\bf R}/\boldz$. Hence the existence of the dual
fibration $\check f$ and the ${\bf R}/\boldz$-simplicity condition
for $p=1$ implies it for $p=2$.

Proof. 
The $E_2$-term
of the Leray spectral sequence, by the arguments of [14], 
Lemma 2.4, looks like
$$\matrix{
\boldz&0&T^{2,3}&\boldz\cr
0&\boldz^{h^{1,2}}\oplus T^{1,2}&\boldz^{h^{1,1}}\oplus T^{2,2}&T^{3,2}\cr
0&\boldz^{h^{1,1}}\oplus T^{1,1}&\boldz^{h^{1,2}}\oplus T^{2,1}&T^{3,1}\cr
\boldz&0&T^{2,0}&\boldz\cr}\leqno{(3.1)}$$
with a similar diagram for $\check f$. Here $T^{i,j}=
H^i(B,R^jf_*\boldz)_{tors}$. We are using 
$H^1(X,\boldz)=H^1(\check X,\boldz)=0$ to obtain the zeroes
on the left column and the top and bottom rows. Clearly the desired
statement on ranks follows.

Recalling from [14], \S 2 that $\boldz$-simplicity
implies $R^if_*\boldz\cong R^{3-i}\check f_*\boldz$, it follows that
if $\check T^{i,j}=H^i(B,R^j\check f_*\boldz)_{tors}$, then $\check T^{i,j}
\cong T^{i,3-j}$.

Since $f$ has a section, the argument of [14], Lemma
2.4, combined with the
degeneration statement of Proposition 3.7,
shows that the above spectral sequence degenerates. The same holds for
$\check f$. 

As for the statements about torsion, clearly $T^{2,3}
\cong T^{2,0}\cong H^2(B,\boldz)$. To show the rest,
we use Poincar\'e-Verdier duality (see for example
[9]). In any dimension, applying duality to the map 
$s:B\rightarrow pt$, we obtain isomorphisms
$${\bf R}\hom({\bf R}{\Gamma}(R^if_*\boldz),\boldz)\cong {\bf R}\Gamma
{\bf R}\lhom(R^if_*\boldz,\boldz[n]).$$
Applying $H^{-j}$ to both sides, we obtain
$$H^{-j}(
{\bf R}\hom({\bf R}{\Gamma}(R^if_*\boldz),\boldz))\cong \ext^{n-j}(R^if_*
\boldz,\boldz).\leqno{(3.2)}$$
The left hand side is easily computed by choosing a complex of
projective $\boldz$-modules quasi-isomorphic to ${\bf R}\Gamma(R^if_*\boldz)$
and applying the algebraic universal coefficient theorem, which yields
exact sequences
$$\exact{\ext^1(H^{j+1}(B,R^if_*\boldz),\boldz)}{H^{-j}({\bf R}\hom(
{\bf R}\Gamma(R^if_*\boldz),\boldz))}{\hom(H^j(B,R^if_*\boldz),\boldz)}.$$
The difficulty in applying Poincar\'e-Verdier
duality for non-locally constant sheaves is the difficulty of comparing
the $\ext$'s and the cohomology groups.
We will do this for $n=3,i=2$. 

First on $B_0$ one has $\lhom(R^2f_{0*}\boldz,\boldz)\cong
R^1f_{0*}\boldz$ by Poincar\'e duality, and if $i:B_0\hookrightarrow
B$ is the inclusion, the natural map
$$i_*\lhom(R^2f_{0*}\boldz,\boldz)
\rightarrow 
\lhom(i_*R^2f_{0*}\boldz,i_*\boldz)$$
is an isomorphism. Thus by $\boldz$-simplicity, $\lhom(R^2f_*\boldz,\boldz)
\cong R^1f_*\boldz$. Also, ${\bf R}/\boldz$-simplicity implies
$\boldz/m\boldz$-simplicity for any $m$, and so a similar argument
shows $\lhom(R^2f_*\boldz,\boldz/m\boldz)\cong R^1f_*\boldz/m\boldz$.

So by the local-global Ext spectral sequence one has a 
five-term sequence 
$$\eqalign{0&\rightarrow H^1(B,R^1f_*\boldz)\rightarrow \ext^1(R^2f_*\boldz,
\boldz)\rightarrow H^0(B,\lext^1(R^2f_*\boldz,\boldz))\cr
&\rightarrow H^2(B,R^1f_*\boldz)\rightarrow \ext^2(R^2f_*\boldz,\boldz).\cr}
\leqno{(3.3)}$$
I claim that $\lext^1(R^2f_*\boldz,\boldz)$ is a torsion-free sheaf.
Indeed, apply $\lhom(R^2f_*\boldz,\cdot)$ to the exact sequence
$$\exact{\boldz}{\boldz}{\boldz/m\boldz}.$$
We obtain an exact sequence
$$0\mapright{} R^1f_*\boldz\mapright{\cdot m}R^1f_*\boldz
\mapright{} R^1f_*\boldz/m\boldz
\mapright{}\lext^1(R^2f_*\boldz,\boldz)\mapright{\cdot m}\lext^1(R^2f_*\boldz,
\boldz).$$
But in fact $R^1f_*(\boldz/m\boldz)\cong R^1f_*\boldz/mR^1f_*\boldz$ 
since $R^2f_*\boldz$ is torsion-free,
so we see the multiplication by $m$ map is injective on
$\lext^1(R^2f_*\boldz,\boldz)$. Thus this sheaf is torsion free.

Now the left hand side of (3.2) is $\boldz^{h^{1,2}}\oplus T^{2,2}$
for $j=1$ and is $\boldz^{h^{1,1}}\oplus T^{3,2}$ for $j=2$. Thus by (3.2),
$$rank_{\boldz} \ext^1(R^2f_*\boldz,\boldz)=h^{1,1}=rank_{\boldz}
H^1(B,R^1f_*\boldz),$$
so in (3.3) the fact that $H^0(B,\lext^1(R^2f_*\boldz,\boldz))$ is
torsion-free shows that the map $$\ext^1(R^2f_*\boldz,\boldz)
\rightarrow H^0(B,\lext^1(R^2f_*\boldz,\boldz))$$ is zero.
Thus we have
$$H^1(B,R^1f_*\boldz)\cong \ext^1(R^2f_*\boldz,\boldz)$$
and
$$0\rightarrow H^0(B,\lext^1(R^2f_*\boldz,\boldz))
\rightarrow H^2(B,R^1f_*\boldz)\rightarrow \ext^2(R^2f_*\boldz,\boldz)
$$
exact. Thus (3.2) implies $T^{1,1}\cong T^{3,2}$, and by using the
same argument for $\check f$, $T^{1,2}\cong T^{3,1}$. In addition,
$T^{2,1}=H^2(B,R^1f_*\boldz)_{tors}\subseteq\ext^2(R^2f_*\boldz,\boldz)_{tors}
=T^{2,2}$. On the other hand, from (3.1) and Proposition 3.7
there are exact sequences (see the beginning of the proof of Theorem 3.10
for details of the first sequence)
$$\exact{T^{2,1}}{H^3(X,\boldz)_{tors}}{T^{1,2}},$$
$$\exact{T^{3,1}}{H^4(X,\boldz)_{tors}}{T^{2,2}},$$
and in addition for any oriented 6-manifold $H^3(X,\boldz)_{tors}
\cong H^4(X,\boldz)_{tors}$ by Poincar\'e duality and the universal
coefficient theorem.
So $\# T^{2,2}=\# T^{2,1}$ and this implies $T^{2,2}\cong T^{2,1}$.
$\bullet$

\proclaim Theorem 3.10. Let $f:X\rightarrow B$, $\check f:\check X
\rightarrow B$ be as in Theorem 3.9, and assume in addition that
$B$ is simply connected. Then there are non-canonical
isomorphisms
$$\eqalign{H^{even}(X,\boldz[1/2])&\cong H^{odd}(\check X,\boldz[1/2])\cr
H^{odd}(X,\boldz[1/2])&\cong H^{even}(\check X,\boldz[1/2]).\cr}$$
In general, there are short exact sequences
$$\exact{H^2(B,R^1f_*\boldz)_{tors}}{H^3(X,\boldz)_{tors}}
{H^1(B,R^2f_*\boldz)_{tors}}$$
$$\exact{H^2(B,R^1\check f_*\boldz)_{tors}}{H^3(\check X,\boldz)_{tors}}
{H^1(B,R^2\check f_*\boldz)_{tors}}$$
and if they split, the above isomorphisms hold over $\boldz$. This happens,
for example, if both $X$ and $\check X$ are simply connected. In any
event,
$$\# H^{even}(X,\boldz)_{tors}=\# H^{odd}(\check X,\boldz)_{tors}$$
$$\# H^{odd}(X,\boldz)_{tors}=\# H^{even}(\check X,\boldz)_{tors}.$$

Proof. The Leray filtration on $H^3(X,\boldz)$ is
$$0\subseteq F_0=\boldz[T^3]\subseteq F_1\subseteq F_2\subseteq F_3=
H^3(X,\boldz).$$
Since the cohomology class $[T^3]$ is primitive in $H^3(X,\boldz)$,
$F_0\subseteq F_1$ is a primitive embedding and
$(F_1)_{tors}=(F_1/F_0)_{tors}=T^{2,1}$ in the notation of the proof of
Theorem 3.9.
It then follows from Proposition 3.7 that there is an exact sequence
$$\exact{T^{2,1}}{H^3(X,\boldz)_{tors}}{T^{1,2}}.$$
First assume this sequence for $X$ and $\check X$ splits,
so $H^3(X,\boldz)_{tors}=T^{2,1}\oplus T^{1,2}$.
Now $H^3(X,\boldz)_{tors}\cong H^4(X,\boldz)_{tors}$, so $H^4(X,\boldz)_{tors}
=T^{2,1}\oplus T^{1,2}$ also. Putting this together we see that
$$H^{even}(X,\boldz)_{tors}=T^{1,1}\oplus T^{2,1}\oplus T^{1,2}$$
and
$$\eqalign{H^{odd}(X,\boldz)_{tors}=&T^{3,2}\oplus T^{2,1}\oplus T^{1,2}\cr
=&T^{1,1}\oplus T^{2,1}\oplus T^{1,2}.\cr}$$
On the other hand
$$\eqalign{H^{odd}(\check X,\boldz)_{tors}=&\check T^{1,1}\oplus \check T^{2,1}
\oplus\check T^{1,2}\cr
=& T^{1,2}\oplus T^{2,2}\oplus T^{1,1}\cr
=&H^{even}(\check X,\boldz)_{tors}
\cr}.$$
Since $T^{2,2}\cong T^{2,1}$ by Theorem 3.9, we are done.

Note that if $X$ and $\check X$ are simply connected, $0=T^{1,1}
\cong \check T^{1,2}$ and $0=\check T^{1,1}\cong T^{1,2}$, so
the sequences trivially split. If the sequences don't split,
then it is still clear that the numerical equalities hold.

Finally, we finish the proof of the theorem by showing the sequences
do split over $\boldz[1/2]$. We define a map $\phi:T^{1,2}\otimes_{\boldz}
\boldz[1/2]\rightarrow H^3(X,\boldz[1/2])_{tors}$.
Indeed, for $\sigma$ a torsion section
of $f:X\rightarrow B$, given by $\sigma\in H^1(B,R^2f_*\boldz)_{tors}$,
we can set $\phi(\sigma)=
(\log T_{\sigma}^*)([\sigma_0])\in H^3(X,\boldz[1/6])$ where
$[\sigma_0]\in H^3(X,\boldz)$ is the cohomology class of the zero section.
Here $\log T_{\sigma}^*=(T_{\sigma}^*-I)-{1\over 2}(T_{\sigma}^*-I)^2
+{1\over 3}(T_{\sigma}^*-I)^3$, as $(T_{\sigma}^*-I)^4=0$ 
by [14] Theorem 4.1. As observed in the proof of Theorem 3.7,
$(T_{\sigma}^*-I)([\sigma_0])=[\sigma]-[\sigma_0]$ is torsion
and represents the class $\sigma\in H^1(B,R^2f_*\boldz)_{tors}$.
Thus $(T_{\sigma}-I)^3([\sigma_0])\in H^3(B,\boldz)\subseteq H^3(X,\boldz)$
must be zero as this element is also torsion. So $(\log T_{\sigma}^*)
([\sigma_0])=((T_{\sigma}^*-I)-{1\over 2}(T_{\sigma}^*-I))([\sigma_0])
\in H^3(X,\boldz[1/2])$. Furthermore, 
$$\eqalign{\phi(\sigma+\tau)&=(\log T^*_{\sigma+\tau})([\sigma_0])\cr
&=(\log T^*_{\sigma}\circ T^*_{\tau})([\sigma_0])\cr
&=(\log T_{\sigma}^*+\log T_{\tau}^*)([\sigma_0])\cr
&=\phi(\sigma)+\phi(\tau),\cr}$$
so $\phi$ is a group homomorphism. Thus $\phi$ gives the desired
splitting over $\boldz[1/2]$. $\bullet$

The problem that arises with two-torsion in the above theorem
seems at the moment to be unavoidable, and does not make the statement
very aesthetically pleasing. The heart of this issue is the following:
given a two-torsion element in 
$H^2(\check X,\boldz)$,
is the square of this element non-zero in $H^4(\check X,\boldz)$?
If it is non-zero, then it follows from [14], Theorem 4.1, that if $\sigma$
is the corresponding torsion section of $f$, then $[\sigma]-[\sigma_0]$
is not two-torsion, making it unlikely that the exact sequences in
Theorem 3.10 split over $\boldz$.  

{\it Example 3.11.} The only example personally known to me of a Calabi-Yau
threefold $X$ with $H^3(X,\boldz)_{tors}$ non-zero is the ``Enriques 
threefold'', obtained by dividing out $K3\times E$ with the involution
$(\iota,-1)$, where $\iota$ is the Enriques involution on the K3 surface.
(See [4] for calculations of the cohomology of this threefold.)
This possesses a special Lagrangian $T^3$ fibration $f$ in much the same way
as the examples in [15]. In fact $H^3(X,\boldz)_{tors}=\boldz/2\boldz$,
and the fibration $f:X\rightarrow B$ is seen to have a torsion section,
so $H^1(B,R^2f_*\boldz)_{tors}=\boldz/2\boldz$. 

{\it Remark 3.12.}
If $f:X\rightarrow B$ does not have a section, then $\check f:\check X
\rightarrow B$ does, as well as the Jacobian $j:J\rightarrow B$ of
$f$. Then $j$ and $\check f$ are dual, and Theorem 3.11 applies to this pair.
In addition, $R^if_*\boldz=R^ij_*\boldz$, so the $E^2$ term of the
Leray spectral sequence for $f$ is still given by (3.1), and
$H^2(B,R^1f_*\boldz)_{tors}\cong H^2(B,R^2f_*\boldz)_{tors}$.
However now the spectral sequence won't degenerate.
Since there is no class in $H^3(X,\boldz)$ restricting to the
generator of $H^3(X_b,\boldz)$, one of the differentials from
$H^0(B,R^3f_*\boldz)$ must be non-zero.

This now gives an explanation for the speculations of [5] of the role
that $H^3(X,\boldz)_{tors}$ should play in mirror symmetry. There it
was argued on physical grounds that the K\"ahler moduli space
of a Calabi-yau threefold was in fact $H^2(X,{\bf C}/\boldz)$, so in
fact it had one component for each element of $H^3(X,\boldz)_{tors}$. Thus the 
complex moduli space of the mirror $\check X$ should have a similar
number of components. It was not clear what this meant. But in
our current context it is clear: if $H^2(B,R^1f_*\boldz)_{tors}\not=0$,
$f:X\rightarrow B$ will have many dual fibrations, only one of which
will have a section. All the other duals are obtained by
twisting the one with a section. The set of such dual fibrations
is classified by $H^2(B,R^1f_*\boldz)_{tors}=H^2(B,R^2\check f_*\boldz)_{tors}$
and we have seen that these groups are related to $H^3(X,\boldz)_{tors}$.
However, they do not necessarily coincide with $H^3(X,\boldz)_{tors}$,
and this will lead us to modify the definition of $B$-field in
Conjecture 6.6.
\bigskip

We end this section with some comments concerning the de Rham realisation
of the Leray spectral sequence, which we will need later.
In general, let $f:X\rightarrow S$ be a smooth map of differentiable manifolds.
Then one has an exact sequence of vector bundles
$$\exact{f^*\Omega^1_S}{\Omega^1_X}{\Omega^1_{X/S}}.$$
This gives rise to a filtration $F^p\Omega^{\cdot}_X$
on the de Rham complex,
$$\Omega_X^{\cdot}=F^0\Omega_X^{\cdot}\supseteq F^1\Omega_X^{\cdot}
\supseteq\cdots$$
such that
$$F^p\Omega_X^{p+q}/F^{p+1}\Omega_X^{p+q}=
\bigwedge^pf^*\Omega^1_S\otimes \bigwedge^q\Omega^1_{X/S}.$$
This filtration gives rise to a spectral sequence with
$$E^0_{pq}=\Gamma(X,\bigwedge^pf^*\Omega^1_S\otimes\bigwedge^q\Omega^1_{X/S})$$
with differential $d^0$ being exterior differentiation along fibres of the
map $f$. Then
$$E^1_{pq}=\Gamma(S,\Omega^p_S\otimes R^qf_*{\bf R}),$$
and $d^1$ is the Gauss-Manin connection $\nabla_{GM}$. Here, given a form
$\alpha\in F^p\Omega^{p+q}_X$ with $d^0\alpha=0$, the element represented
in $E^1_{pq}$ is given as follows. For $v_1,\ldots,v_p\in \T_{B,b}$,
the form $(\iota(v_1,\ldots,v_p)\alpha)|_{X_b}$ defines a well-defined
cohomology class in $H^q(X_b,{\bf R})$. This yields an $R^qf_*{\bf R}$-valued
$p$-form in $E^1_{pq}$. Next 
$$E^2_{pq}=H^p(S,R^qf_*{\bf R}),$$
which coincides with the $E^2$-term of the Leray spectral
sequence for $f$.

Now let us specialise to the case that $f:X\rightarrow B$ is a 
special Lagrangian $T^n$-fibration, $f_0:X_0\rightarrow B_0=B-\Delta$
the smooth part of the fibration. On $B_0$, McLean's result
gives a natural isomorphism
$$\T_{B_0}\cong R^1f_{0*}{\bf R}\otimes C^{\infty}(B_0).$$
Mclean also defines an $n$-form on the base. This is given by
$$\Theta(v_1,\ldots,v_n)
=\int_{X_b} (-\iota(v_1)\omega)\wedge\cdots\wedge(-\iota(v_n)\omega).$$
In canonical coordinates, this is
$$\Theta(\partial/\partial y_1,\ldots,\partial/\partial y_n)
=\int_{X_b} dx_1\wedge\cdots\wedge dx_n.$$
Here $X_b$ is oriented canonically by $\Omega$ as it is special Lagrangian.
Of course, this form goes to infinity
at singular fibres.
Another way to think about this form is via integration along
fibres of $\omega^n/n!$. Since 
$$\omega^n/n!=(-1)^{n(n+1)/2}dy_1\wedge\cdots\wedge dy_n\wedge dx_1
\wedge\cdots\wedge dx_n,$$
$f_*(\omega^n/n!)=(-1)^{n(n+1)/2}\Theta$.
Thus in particular,
$$\int_X \omega^n/n!=(-1)^{n(n+1)/2}\int_B\Theta.$$

We can identify 
$\bigwedge^{n-q}\T_{B_0}$ with $\Omega^q_{B_0}$ by contracting with
$\Theta$, and so obtain isomorphisms
$$\Omega_{B_0}^q\cong R^{n-q}f_{0*}{\bf R}\otimes C^{\infty}(B_0).$$
Thus the $E^1_{p,n-q}$ term of the de Rham realisation of the Leray spectral
sequence for $f_0$ is $\Gamma(B_0,\Omega_{B_0}^p\otimes \Omega_{B_0}^q)$.

\proclaim Definition 3.13. A cohomology class $[\alpha_0]\in H^1(B_0,
R^{n-1}f_{0*}{\bf R})$ 
is {\it symmetric} if there is a representative $\alpha_0\in E^1_{1,n-1}=
\Gamma(B_0,\Omega^1_{B_0}\otimes\Omega^1_{B_0})$ of $[\alpha_0]$
which is invariant under the involution of $\Omega^1_{B_0}\otimes\Omega^1_{B_0}$
given by $a\otimes b\mapsto b\otimes a$. In other words, 
$\alpha_0\in \Gamma(B_0,S^2\Omega^1_{B_0})\subseteq \Gamma(B_0,
\Omega^1_{B_0}\otimes \Omega^1_{B_0})$. A cohomology class
$\alpha\in H^1(B,R^{n-1}f_*{\bf R})$ is said to be {\it symmetric} if
its restriction to $H^1(B_0, R^{n-1}f_{0*}{\bf R})$ is symmetric.

The following result, which identifies the symmetric cohomology classes,
will be useful in \S 6 in studying the role of the $B$-field.

\proclaim Theorem 3.14. A cohomology class $[\alpha_0]\in H^1(B_0,R^{n-1}f_{0*}
{\bf R})$ is symmetric if and only if $[\alpha_0]\wedge [\omega]
\in H^2(B_0,R^nf_{0*}\bf R)$ is zero. In particular, if 
$H^2(B,{\bf R})=0$, then all elements of $H^1(B,R^{n-1}f_*{\bf R})$
are symmetric.

Proof. We work in action-angle coordinates, on a neighborhood $U$,
so that the Gauss-Manin connection is the trivial connection.
Thus if $y_1,\ldots,y_n$ are the action coordinates, 
$$\nabla_{GM}(\sum f_{IJ}dy_I\otimes dy_J)=
\sum d(f_{IJ}dy_I)\otimes dy_J.$$
In addition, in suitably ordered
action-angle coordinates $\Theta=dy_1\wedge\cdots\wedge dy_n$,
from which one sees easily that the cohomology class of 
$(-1)^{i-1} dx_1\wedge\cdots\wedge \hat{dx_i}\wedge\cdots\wedge dx_n$
on a fibre $X_b$ is identified with $dy_i\in \Omega^1_{B,b}$.

Choose a representative $\alpha_0\in \Gamma(B_0,\Omega^1_{B_0}\otimes
\Omega^1_{B_0})$ for $[\alpha_0]$, so on $U$ we write
$$\alpha_0=\sum_{i,j} \alpha_{ij} dy_i\otimes dy_j,\leqno{(3.4)}$$ with
$\alpha_0$ symmetric if and only if $\alpha_{ij}=\alpha_{ji}$.
Locally, this class can also be represented by the $n$-form
on $f^{-1}(U)$
$$\alpha_0=\sum_{i,j} (-1)^{j-1} \alpha_{ij} dy_i\wedge
dx_1\wedge\cdots\wedge \hat{dx_j}\wedge\cdots\wedge dx_n,$$
from which we see that
$$\eqalign{\alpha_0\wedge\omega&=
\sum_{i,j} (-1)^{j-1}\alpha_{ij}
dy_i\wedge dx_1\wedge\cdots\wedge\hat{dx_j}\wedge
\cdots\wedge dx_n\wedge dx_j\wedge dy_j\cr
&=\sum_{i<j} (\alpha_{ji}-\alpha_{ij})dy_i\wedge dy_j\wedge
dx_1\wedge\cdots\wedge dx_n.\cr}\leqno{(3.5)}$$
Thus we see $\alpha_0$ is a symmetric representative for $[\alpha_0]$
if and only if $\alpha_0\wedge\omega=0$ in $\Gamma(B_0,\Omega^2_{B_0}
\otimes \Omega_{B_0}^0)$. Now $\alpha_0\wedge\omega$ represents
the cohomology class $[\alpha_0]\wedge[\omega]\in H^2(B_0,R^nf_{0*}{\bf R})$.
This is zero if $\alpha_0$ is symmetric. Conversely, suppose
$[\alpha_0]\wedge[\omega]$ is the zero cohomology class. Then
there exists a $\beta\in \Gamma(B_0,\Omega^1_{B_0}\otimes
\Omega^0_{B_0})$ such that $\nabla_{GM} \beta=\alpha_0\wedge\omega$.
$\beta$ also gives rise to an element $\beta'\in \Gamma(B_0,\Omega_{B_0}^0
\otimes \Omega^1_{B_0})$ by using the map $a\otimes b\mapsto b\otimes a$.
The following claim then proves the theorem.

{\it Claim.} $\alpha_0+\nabla_{GM}\beta'$ is a symmetric representative
for $[\alpha_0]$.

Proof. Locally write $\beta=\sum_i \beta_i dy_i\otimes 1$, so
$\beta'=\sum_j \beta_j 1\otimes dy_j$. Now $\nabla_{GM}(\beta)
=\alpha_0\wedge\omega$ means that
$${\partial\beta_j\over \partial y_i}-{\partial\beta_i\over\partial y_j}
=\alpha_{ji}-\alpha_{ij}.$$
Also,
$$\alpha_0+\nabla_{GM}\beta'=\sum_{i,j}(\alpha_{ij}+{\partial \beta_j\over
\partial y_i}) dy_i\otimes dy_j.$$
But
$$\alpha_{ij}+{\partial\beta_j\over\partial y_i}-(\alpha_{ji}+{\partial\beta_i
\over\partial y_j})=0,$$
so $\alpha_0+\nabla_{GM}\beta'$ is a symmetric representative for
$[\alpha_0]$. $\bullet$

{\it Remark 3.15.} There is a related map which partly explains the interest
in symmetric cohomology classes. If $f:X\rightarrow B$ is an
integral, $\boldz$-simple special Lagrangian fibration with a section,
then there is a natural inclusion $R^{n-1}f_*{\bf R}/\boldz\hookrightarrow
\Lambda(X^{\#})$. This is induced from the inclusion $R^{n-1}f_*\boldz
\hookrightarrow \Lambda(\T_B^*)$ and the inclusion obtained from this
by tensoring with ${\bf R}$: $R^{n-1}f_*{\bf R}\hookrightarrow
\Lambda(\T_B^*)$. Thus one obtains a natural map $H^1(B,R^{n-1}f_*
{\bf R}/\boldz)\rightarrow H^1(B,\Lambda(X^{\#}))$. To analyse
this map, consider instead the map $H^1(B,R^{n-1}f_*{\bf R})
\rightarrow H^1(B,\Lambda(\T_B^*))\cong H^2(B,{\bf R})$. This in fact
coincides with the map $H^1(B,R^{n-1}f_*{\bf R})\rightarrow H^2(B,R^nf_*
{\bf R})\cong H^2(B,{\bf R})$ obtained by cup-product with $-\omega$.
It is easiest to see this over $B_0$: in this case,
we have resolutions $0\rightarrow R^{n-1}f_{0*}{\bf R}
\rightarrow \Omega^{\cdot}_{B_0}\otimes\Omega^1_{B_0}$ and
$0\rightarrow\Lambda(\T_{B_0}^*)\rightarrow\Omega_{B_0}^{\cdot+1}$.
The map $R^{n-1}f_{0*}{\bf R}\rightarrow\Lambda(\T_{B_0}^*)$ extends
to a map of complexes $\Omega_{B_0}^{\cdot}\otimes\Omega^1_{B_0}
\rightarrow \Omega_{B_0}^{\cdot+1}$ given by $\alpha\otimes\beta
\mapsto\alpha\wedge\beta$. Thus $\alpha_0$ as in (3.4) is mapped to
$\sum_{i<j} (\alpha_{ij}-\alpha_{ji})dy_i\wedge dy_j$, which by
(3.5) can be identified with $-\alpha_0\wedge\omega$.
Thus the map $H^1(B_0,
R^{n-1}f_{0*}{\bf R})\rightarrow H^1(B_0,\Lambda(\T_{B_0}^*))$ coincides
with the map $\wedge(-\omega)$. It is not difficult to show this holds
over $B$ also, but we omit the cohomological argument.

Thus we see that the symmetric cohomology classes are those that
map to zero in $H^1(B,\Lambda(X^{\#}))$.
\bigskip

{\hd \S 4. The symplectic form on $D$-brane moduli space.}

Let $B$ be a moduli space of special Lagrangian submanifolds in a 
Calabi-Yau manifold $X$ of dimension $n$, along with a universal family
$$\matrix{\U&\subseteq&X\times B\cr
\mapdown{f}&\cr
B&&\cr}$$
Let $p:\U\rightarrow X$ be the projection. For the moment we
assume $f$ is smooth, so that points in $B$ are
parametrizing only smooth special Lagrangian submanifolds. We do not assume
these submanifolds are tori. Here $\dim B=\dim H^1(\U_b,{\bf R})=:s$.
The $D$-brane moduli space is the space of special
Lagrangian submanifolds along with a choice of flat $U(1)$ connection modulo
gauge equivalence,
i.e. an element of $H^1(\U_b,{\bf R}/\boldz)$. Thus the $D$-brane
moduli space $\M$ should be $R^1f_*({\bf R}/\boldz)$. The prediction
from string theory
is that $\M$ should be a complex K\"ahler manifold, so we need to
understand how to put these structures on $\M$. As long as there
are no singular fibres to deal with, Hitchin has described how to
put a complex and K\"ahler structure on $\M$. Here, we will
describe a more coordinate independent way of describing the same
K\"ahler form (i.e. a symplectic form in the absence of a complex
structure). This both allows us to compute the cohomology class
represented by this symplectic form and in principle should allow one
to extend this construction to singular fibres.

In what follows, we assume the fibres of $f$ are special Lagrangian
with respect to the symplectic form $\omega$ and holomorphic
$n$-form $\Omega$, with the standard normalization $\omega^n/n!
=(-1)^{n(n-1)/2}(i/2)^n\Omega\wedge\bar\Omega$. We also use a 
holomorphic $n$-form $\Omega_n$ normalised by
$$\Omega_n={\Omega\over\int_{\U_b} p^*\Omega}.$$
We use $\Omega_n$ instead of $\Omega$ for several reasons. First, we
do not want the symplectic structure on $\M$
to depend on $\omega$, but only on $\Omega$. If we multiply
$\omega$ by a constant, we must also rescale $\Omega$. If 
we rescale $\Omega$ and use $\Omega$ instead
of $\Omega_n$ in the construction below, 
then the symplectic form $\check\omega$ we
construct on $\M$ also changes. Secondly, this normalization fits
with the usual form of the mirror map as described in item (5) of the
introduction.

To obtain a symplectic form on $\M$, we define a map
$$R^1f_*\boldz\rightarrow \T_B^*$$
in such a way that the canonical symplectic form on $\T_B^*$
descends to a symplectic form on $\T_B^*/R^1f_*\boldz$. We follow
Hitchin's suggestion of computing the periods of
$\im p^*\Omega_n$. Now
$(R^1f_*\boldz)_b\cong H^1(\U_b,\boldz)\cong H_{n-1}(\U_b,\boldz)$,
so for $\gamma\in H_{n-1}(\U_b,\boldz)$, map $\gamma$ to the differential
$$v\mapsto -\int_{\gamma} \iota(v)\im p^*\Omega_n$$
where again we choose an arbitrary lifting of $v$ to $\U$.

\proclaim Lemma 4.1. The image of $R^1f_*\boldz$ in $\T_B^*$ is Lagrangian 
with respect to the standard symplectic form on $\T_B^*$.
Thus $\M=\T_B^*/R^1f_*\boldz$ inherits this symplectic form, which we
will call $\check\omega$.

Proof. See Hitchin's paper [17]. His proof is as follows: in a small open set 
of $B$, choose $\Gamma\subseteq \U$ a family of $n-1$-dimensional submanifolds
representing a section of $R^1f_*\boldz$ over $U$, $\pi:\Gamma\rightarrow
B$ the projection. Then the section of $\T_B^*$ obtained by
taking periods with respect to $\Gamma$ is just the $1$-form
$-\pi_*((p^*\im\Omega)|_{\Gamma})$. Since $\Omega$ is a closed form, so is
this push-down, and hence $-\pi_*((p^*\im\Omega|_{\Gamma})$ is a
Lagrangian section of $\T_B^*$ with respect to the standard symplectic
form. $\bullet$

We will now clarify what the cohomology class of $\check\omega$ is. To do
so, we will compare the Leray spectral sequences for $f$ and
$\check f:\M\rightarrow B$, but will use the de Rham realisation of 
these spectral sequences discussed in \S 3, 
which we can do as $f$ and $\check f$ are smooth.
Our construction yields a canonical isomorphism $H_{n-1}(\U_b,
\boldz)\cong H_1(\M_b,\boldz)$ and hence a canonical isomorphism
$H^{n-1}(\U_b,\boldz)\cong H^1(\M_b,\boldz)$, which yields a canonical
isomorphism
$$R^{n-1}f_*{\bf R}\cong R^1\check f_*{\bf R}.\leqno{(4.1)}$$
Since $\im p^*\Omega_n$ restricts to zero on the fibres of $f$,
$\im p^*\Omega_n\in F^1\Omega^n_{\U},$ and since $d\Omega_n=0$, $\im p^*\Omega_n$
in particular gives rise to a class $[\im p^*\Omega_n]$ in $E^1_{1,n-1}
=\Gamma(B,\Omega^1_B\otimes R^{n-1}f_*{\bf R})$. 
Now on $\M$, $\check\omega\in F^1(\Omega^2_{\M})$, as $\check\omega$ restricts
to zero on the fibres of $\check f$, and thus $\check\omega$ determines
an element $[\check\omega]$ of $\check E^1_{1,1}=\Gamma(B,\Omega^1_B\otimes
R^1\check f_*{\bf R})$.
Via the isomorphism (4.1) we can identify $E^1_{1,n-1}$ and $\check E^1_{1,1}$.

\proclaim Proposition 4.2. Under this identification, $[\check\omega]=
[\im p^*\Omega_n]$. Thus in particular, they represent the same class in
$H^1(B,R^1\check f_*{\bf R})\cong H^1(B,R^{n-1}f_*{\bf R})$.

Proof. A section of $\Gamma(B,\Omega^1_B\otimes R^{n-1}f_*{\bf R})$
associates, to any tangent vector $v\in \T_{B,b}$, an element of
$H^{n-1}(\U_b,{\bf R})$. Specifically, $[\im p^*\Omega_n]$ associates
to a tangent vector $v\in\T_{B,b}$ the cohomology class represented
by $\iota(v)(\im p^*\Omega_n)$. To determine what cohomology class this is,
we choose a basis $\gamma_1,\ldots,\gamma_s$ of $H_{n-1}(\U_b,\boldz)$
and calculate the periods
$$\int_{\gamma_i} \iota(v)\im p^*\Omega_n.$$
On $\M$, $[\check\omega]\in \Gamma(B,\Omega^1_B\otimes R^1\check f_*{\bf R})$
is similarly represented by 
$$v\mapsto \iota(v)\check\omega,$$
and $\gamma_1,\ldots,\gamma_s$ also form a basis for $H_1(\M_b,\boldz)$
by construction. Recall that we embedded $H_{n-1}(\U_b,\boldz)
=H_1(\M_b,\boldz)$ in $\T_{B,b}^*$ by mapping $\gamma_i$ to the
differential 
$$v\mapsto -\int_{\gamma_i}\iota(v)\im p^*\Omega_n.$$
Now choosing local coordinates $y_1,\ldots,y_s$ on the base,
$x_1,\ldots,x_s,y_1,\ldots,y_s$ canonical coordinates on $\T_B^*$,
$$\int_{\gamma_i} \iota(\partial/\partial y_j)\check\omega
=-\int_{\gamma_i} dx_j=\int_{\gamma_i} \iota(\partial/\partial y_j)
\im p^*\Omega_n$$
by construction. Thus $[\check\omega]=[\im p^*\Omega_n]$. $\bullet$
\bigskip

We recall here for future use:

\proclaim Observation 4.3. (McLean [21]) Because $\T_{B,b}$ is naturally
isomorphic to the space of harmonic 1-forms on $\U_b$, there is a
metric $h$ on $\T_B$ coming from the Hodge metric. Precisely,
for $v\in\T_{B,b}$, $-(\iota(v)p^*\omega)|_{X_b}$ is the corresponding
harmonic one-form, and $-*(\iota(v)p^*(\omega))|_{X_b}
=(\iota(v)\im p^*(\Omega))|_{X_b}$, and we define, for
$v,w\in\T_{B,b}$,
$$h(v,w)=-\int_{X_b} (\iota(v)p^*\omega)\wedge(\iota(v)\im
p^*\Omega).$$
This is a Riemannian metric on $B$.

Specialising down to the case that
$f:X\rightarrow B$ is a special Lagrangian
$T^n$-fibration with possible degenerate fibres,
the above method gives us a way of constructing an open subset of the dual
fibration along with a symplectic form on that open set. On
$B_0=B-\Delta$, we have defined an embedding $R^1f_{0*}\boldz
\hookrightarrow \T_{B_0}^*$. This allows us to define $\check X_0$
via the exact sequence
$$\exact{R^1f_{0*}\boldz}{\T_{B_0}^*}{\check X_0}$$
and $\check X_0$ acquires a symplectic form $\check \omega$
inherited from the canonical symplectic form on
$\T_{B_0}^*$. 

Next we need to prove

\proclaim Conjecture 4.4. The embedding $R^1f_{0*}\boldz\hookrightarrow
\T_{B_0}^*$ extends to an embedding $R^1f_*\boldz\hookrightarrow \T_B^*$.
If $\check X^{\#}$ is defined as $\T_B^*/R^1f_*\boldz$, then $\check X^{\#}$
is a manifold with symplectic form $\check\omega$ inherited from
the standard symplectic form on $\T_B^*$. Furthermore, $\check X^{\#}$
can be compactified to a manifold $\check X$ with a map
$\check f:\check X\rightarrow B$ extending $\check f^{\#}$,
on which $\check\omega$
extends to a symplectic form on $\check X$.

This involves first understanding the asymptotic behaviour of
the periods as one approaches $\Delta$, as well as understanding the issue
of compactification.

If this conjecture holds and 
$f$ and $\check f$ are both ${\bf R}$-simple, then it is easy to
see that $[\check \omega]\in
H^1(B,R^1\check f_*{\bf R})$ coincides with $[\im\Omega_n]
\in H^1(B,R^{n-1}f_*{\bf R})$. Indeed, by Proposition 4.2,
these classes agree in $H^1(B_0,R^1\check f_*{\bf R})
\cong H^1(B_0,R^{n-1}f_*{\bf R})$. However, since $f$ and 
$\check f$ were assumed to be simple, $H^i_{\Delta}(B, R^jf_*{\bf R})=0$
for $i=0$ and $1$, and ditto for $\check f$. Thus there is an injection
$H^1(B,R^if_*{\bf R})\hookrightarrow H^1(B-\Delta, R^if_*{\bf R})$,
and so the classes $[\check \omega]$ and $[\im\Omega]$ agree also
in $H^1(B,R^1\check f_*{\bf R})\cong H^1(B,R^{n-1}f_*{\bf R})$.

$\check f:\check X\rightarrow B$ is not the only possible Lagrangian fibration
we might construct. This fibration possesses a Lagrangian section by 
construction. By Theorem 2.6, any element of $H^1(B,\Lambda(\check X^{\#}))$
gives rise to another, locally isomorphic Lagrangian fibration $\check g:
\check Y\rightarrow B$.  See also Remark 3.12 and Conjecture 6.6.

{\it Remark 4.5.} Having constructed a symplectic form $\check\omega$
on $\check X$, or on an open subset of $\check X$,
$\check\omega^n/n!$ defines an orientation on $\check X$. Thus we can
check that this agrees with the choice of orientation on $\check X$
made in [14], Convention 4.3. First, note that in having fixed
$\Omega$, we have fixed an orientation on the fibres of $f:X\rightarrow B$.
If we have fixed canonical coordinates $y_1,\ldots,y_n,x_1,\ldots,x_n$,
then $\Omega|_{X_b}=Vdx_1\wedge\cdots\wedge dx_n$ with $V$ a real
function and either $V>0$ or $V<0$. By changing the order of the variables
$y_i$, we can ensure $V>0$, and then $dy_1\wedge\cdots\wedge dy_n$
yields a canonical orientation on $B$. We will always assume our coordinates
are so oriented. Note that this orientation on $B$ is the same as that
induced by the $n$-form $\Theta$ on $B_0$.

Now let us check Convention 4.3 of [14] is correct. Recall
that the convention of [14] for the cohomology class of a fibre,
$[X_b]$, was that
$$\int_{X_b}\alpha=\int_X\alpha\wedge [X_b].$$
With $\alpha=[\Omega]$, we then have
$$0<\int_{X_b}\Omega=\int_X \Omega\wedge [X_b].$$
We can take $[X_b]$ to be the pull-back of a nowhere zero $n$-form
on $B$, locally $fdy_1\wedge\cdots\wedge dy_n$. Then $\Omega\wedge
[X_b]=Vf dx_1\wedge\cdots\wedge dx_n\wedge dy_1\wedge\cdots\wedge dy_n$,
while $\omega^n/n!=(-1)^{n(n-1)/2}dx_1\wedge\cdots\wedge dx_n\wedge
dy_1\wedge\cdots\wedge dy_n$.
Since $V>0$, we need $sign(f)=(-1)^{n(n-1)/2}$. Thus we take $[X_b]$
locally to be of the form $(-1)^{n(n-1)/2}|f|dy_1\wedge\cdots\wedge
dy_n$. Now the dual class of $(-1)^n[X_b]$ in $H^n(B,R^n\check f_*{\bf R})$
will be locally represented by something like $(-1)^{n(n+1)/2}g
dy_1\wedge\cdots\wedge dy_n\wedge dx_1\wedge\cdots\wedge dx_n$, $g>0$.
This is the same sign as $\check\omega^n/n!=(-1)^{n(n+1)/2}dy_1
\wedge\cdots\wedge dy_n\wedge dx_1\wedge\cdots\wedge dx_n$, so the orientations
agree.
\bigskip

In the case of torus fibrations, we now describe an alternative
way of putting a symplectic form on $\check X_0$. We do this by
providing an alternative description of the embedding $R^1f_{0*}\boldz
\hookrightarrow \T_{B_0}^*$, using the Riemannian metric $h$ on the base
defined in Observation 4.3.
Since we obtained in \S 2 an embedding $R^{n-1}f_{0*}\boldz\rightarrow
\T_{B_0}^*$, we obtain dually an embedding $R^1f_{0*}\boldz\rightarrow
\T_{B_0}$. We will then use a normalised form of $h$ to identify
$\T_{B_0}$ and $\T_{B_0}^*$.

First, given that there is an identification between $H^1(X_b,{\bf R})$
and $\T_{B,b}$, via $v\in \T_{B,b}\mapsto -\iota(v)\omega$,
we have also a canonical identification of $H_{n-1}(X_b,{\bf R})$ with
$\T_{B,b}$,
via Poincar\'e duality. In fact, write $X_b=V/\Lambda$,
$V=\T_{B,b}^*$, $\Lambda=H_1(X_b,\boldz)\subseteq V$,
$\dual{V}=\T_{B,b}$, $\dual{\Lambda}=\{\varphi\in \dual{V}|
\varphi(\Lambda)\subseteq\boldz\}\subseteq \dual{V}$. 
There is a canonical identification of $\dual{\Lambda}$
with $H^1(X_b,\boldz)$. On the other hand,
the identification $\bigwedge^n \Lambda\cong\boldz$ determined by the 
orientation on $X_b$ gives us a natural identification of $H_{n-1}
(X_b,\boldz)=\bigwedge^{n-1}\Lambda$ with $\dual{\Lambda}$ via
the perfect pairing
$$\Lambda\times \bigwedge^{n-1}\Lambda\rightarrow\bigwedge^n\Lambda
\mapright{\cong}\boldz.$$
(Note: Whenever we use Poincar\'e duality, there is an arbitrary
choice of order in this pairing which may affect the signs of the
isomorphisms. This was seen in [Slag I], where certain conventions
were chosen. Here we also make a choice, and keep in mind that
we could just as well have chosen the pairing $\bigwedge^{n-1}\Lambda
\times\Lambda\rightarrow\boldz$.)

\proclaim Proposition 4.6. If $\alpha\in H^{n-1}(X_b,{\bf R})$, $\gamma
\in H_{n-1}(X_b,\boldz)\cong \dual{\Lambda}\subseteq \T_{B,b}$
via the above identification, then
$$\int_{\gamma}\alpha=-\int_{X_b}\iota(\gamma)\omega\wedge\alpha.$$

Proof. We compute both sides using local action-angle coordinates.
Let $y_1,\cdots,y_n$ be action coordinates as in Remark 2.8.
Then the lattice $\Lambda\subseteq V=\T_{B,b}^*$ is generated by 
$e_1,\ldots,e_n$, $e_i=dy_i$, and then $e_1^*,\ldots,e_n^*$
form a dual basis for $\dual{\Lambda}\subseteq \T_{B,b}$, $e_i^*=\partial/
\partial y_i$. Suppose $\alpha=
dx_1\wedge\cdots\wedge \hat{dx_j}\wedge\cdots\wedge dx_n$.
Then 
$$\eqalign{-\int_{X_b} \iota(e_i^*)\omega\wedge\alpha&=
\int_{X_b} dx_i\wedge\alpha\cr
&=(-1)^{j-1}\delta_{ij}.\cr}$$

On the other hand, the isomorphism between $\dual{\Lambda}$ and
$\bigwedge^{n-1}\Lambda$ identifies $e_i^*$ with
$(-1)^{i-1}e_1\wedge\cdots\wedge \hat e_i\wedge\cdots \wedge e_n$.
The latter defines an oriented $n-1$-torus in $X_b$, namely
the quotient of the subspace of $V$ spanned by $e_1,\ldots,e_{i-1},e_{i+1},
\ldots,e_n$ by the lattice generated by these vectors, and
$\int_{e_i^*}\alpha$ is just the integral of $\alpha$ over this torus,
which is clearly $(-1)^{j-1}\delta_{ij}$. Thus
$$\int_{e_i^*}\alpha=-\int_{X_b} \iota(e_i^*)\omega\wedge\alpha$$
and the result follows from linearity. $\bullet$
\bigskip

This now gives us the opportunity to rephrase the embedding
$H^1(X_b,\boldz)\hookrightarrow \T_{B,b}^*$.
Let $h_n$ be the normalised metric on $B_0$ given by
$$h_n(v,w)={h(v,w)\over\int_{X_b} \Omega}.$$
Then we have

\proclaim Proposition 4.7. For
$\gamma\in H^1(X_b,\boldz)\cong H_{n-1}(X_b,\boldz)=\dual{\Lambda}
\subseteq \T_{B,b}$, the $1$-form
$$v\mapsto -\int_{\gamma} \iota(v)\im\Omega_n$$
coincides with the $1$-form $-h_n(\gamma,\cdot)$. Thus the embedding
$R^1f_{0*}\boldz\rightarrow\T_{B_0}^*$ previously defined coincides 
up to sign with the embedding $R^1f_{0*}\boldz\rightarrow \T_{B_0}$
composed with the isomorphism $\T_{B_0}\cong \T_{B_0}^*$ induced by
the Riemannian metric $h_n$.

Proof. By Proposition 4.6,
$$\eqalign{-\int_{\gamma}\iota(v)\im\Omega_n&=
\int_{X_b}\iota(\gamma)\omega\wedge\iota(v)\im\Omega_n\cr
&={1\over \int_{X_b}\Omega}\int_{X_b}\iota(\gamma)\omega
\wedge\iota(v)\im\Omega\cr
&=-{h(\gamma,v)\over\int_{X_b}\Omega}.\bullet\cr}$$

In fact, we see that $h_n$ also describes the class
$[\im \Omega]\in E^1_{1,n-1}=\Gamma(B_0,\Omega^1_{B_0}\otimes 
R^{n-1}f_{0*}{\bf R})$:

\proclaim Proposition 4.8.  Under the isomorphism
$$E^1_{1,n-1}=\Gamma(B_0,\Omega^1_{B_0}\otimes R^{n-1}f_{0*}{\bf R})
\cong\Gamma(B_0,\Omega^1_{B_0}\otimes\Omega^1_{B_0})$$
given in \S 3, $[\im\Omega_n]$ coincides with $h_n\in
\Gamma(B_0,S^2\Omega^1_{B_0})\subseteq \Gamma(B_0,\Omega^1_{B_0}
\otimes\Omega^1_{B_0})$.

Proof. First note that for a point $b\in B$,
$e^*_1,\ldots,e^*_n$ a basis of $H_{n-1}(X_b,\boldz)\cong
\dual{\Lambda}\subseteq\T_{B,b}$, $[\im\Omega_n]$ associates to a vector
$v\in\T_{B,b}$ the class of $\iota(v)\im\Omega_n\in H^{n-1}(X_b,{\bf R})$,
and in terms of the periods,
$$\eqalign{\int_{e^*_i}\iota(v)\im\Omega_n&=
-\int_{X_b}\iota(e^*_i)\omega\wedge\iota(v)\im\Omega_n\cr
&=h_n(e^*_i,v).\cr}$$
Thus, in action-angle coordinates as used in the proof of Theorem 3.14,
$[\im\Omega_n]$ corresponds to the element of $\Gamma(B_0,
\Omega^1_{B_0}\otimes R^{n-1}f_{0*}{\bf R})$ given as
$$\sum_{i,j} h_n(\partial/\partial y_i,\partial/\partial y_j)
dy_i\otimes (-1)^{j-1}dx_1\wedge\cdots\wedge \hat dx_j\wedge\cdots
\wedge dx_n$$
which coincides with
$$\sum_{i,j} h_n(\partial/\partial y_i,\partial/\partial y_j)
dy_i\otimes dy_j\in \Gamma(B_0,S^2\Omega^1_{B_0})$$
as desired. $\bullet$
\bigskip

{\hd \S 5. Complex structures on special Lagrangian torus fibrations.}

Recall from [17] the following: (for K3 surfaces, this was noticed in [28];
see also [27].)

\parskip=0pt
\proclaim Theorem 5.1. Let $X$ be a real $2n$-dimensional
manifold. If
$\Omega$ is a complex-valued $C^{\infty}$
$n$-form on $X$ satisfying the three properties
\item{(1)} $d\Omega=0$;
\item{(2)} $\Omega$ is locally decomposable (i.e. can be written locally
as $\theta_1\wedge\cdots\wedge \theta_n$ where $\theta_1,\ldots,\theta_n$
are 1-forms);
\item{(3)} $(-1)^{n(n-1)/2}(i/2)^n\Omega\wedge\bar\Omega>0$ everywhere on $X$,

{\noindent\it
then $\Omega$ determines a complex structure on $X$ for which $\Omega$ is
a holomorphic $n$-form.} 

\parskip=3pt plus1pt minus.5pt

\parskip=0pt
\proclaim Theorem 5.2. Let $X$ be a real $2n$-dimensional
manifold. Suppose $\omega$ is a symplectic
form on $X$ and
$\Omega$ is a complex-valued $n$-form on $X$ such that
\item{(1)} $\Omega$ satisfies the conditions of Theorem 5.1;
\item{(2)} $\omega$ is a positive $(1,1)$ form in the complex structure of
Theorem 5.1;
\item{(3)} $(-1)^{n(n-1)/2}(i/2)^n\Omega\wedge\bar\Omega=\omega^n/n!$.

{\noindent\it
Then $\Omega$ induces a complex structure on $X$ such that $\omega$
is a K\"ahler form on $X$ whose corresponding metric is Ricci-flat.}

\parskip=3pt plus1pt minus.5pt
Now let $f:X\rightarrow B$ be an integral special Lagrangian fibration
with a Lagrangian section,
so that in local coordinates, $\omega$ takes
the standard form. Let $\Omega$ be the holomorphic
$n$-form on $X$ normalised so that $\im\Omega|_{X_b}=0$ for all $b\in B$,
$\omega^n/n!=(-1)^{n(n-1)/2}(i/2)^n\Omega\wedge
\bar\Omega$, and $\int_{X_b}\Omega=Vol(X_b)$ with respect to
the metric induced by the K\"ahler form $\omega$. As before, 
we set $\Omega_n=\Omega/\int_{X_b} \Omega$.

Following the suggestion of [17], 
since $\Omega$ is locally decomposable, we can write, for local
coordinates $y_1,\ldots,y_n$ on $B$  as usual,
$$\Omega=V\bigwedge_i (dx_i+\sum_j \beta_{ij} dy_j),$$
where $V$ is a real function of $y_1,\ldots,x_n$ and $V|_{X_b}dx_1\wedge\cdots
\wedge dx_n$ is the volume form on $X_b$, while $\beta_{ij}$ is
a complex valued function. Using Remark 4.5, we will always
assume that $V>0$. The forms of type $(1,0)$ are spanned by the
1-forms $\theta_i:=dx_i+\sum_{i,j} \beta_{ij} dy_j$. Thus the entire complex
structure is encoded in the matrix $(\beta_{ij})$. 
We now look to see how the conditions of Theorems 5.1 and 5.2
translate into conditions for the functions
$V$ and $\beta_{ij}$:

\proclaim Calculation 5.3. 
$$(-1)^{n(n-1)/2} (i/2)^n\Omega\wedge\bar\Omega
=V^2\det(\im\beta)\omega^n/n!.$$

Proof.
Write
$$\Omega=V\theta_1\wedge\cdots\wedge\theta_n.$$
Now
$$\eqalign{\theta_i\wedge \bar\theta_i &=\sum_j (\bar\beta_{ij}-\beta_{ij})
dx_i\wedge dy_j+\cdots\cr
&=-2idx_i\wedge\left(\sum_j \im\beta_{ij} dy_j\right)+\cdots.\cr}$$
Thus
$$\Omega\wedge\bar\Omega=V^2(-2)^n i^n\det(\im\beta)dx_1\wedge dx_2\wedge\dots
\wedge dy_n$$
so 
$$(-1)^{n(n-1)/2}(i/2)^n\Omega\wedge\bar\Omega=(-1)^{n(n-1)/2}V^2
\det(\im\beta)dx_1\wedge dx_2\wedge\dots
\wedge dy_n.$$
On the other hand, 
$$\omega^n/n!=(-1)^{n(n-1)/2}dx_1\wedge\cdots\wedge dy_n,$$
hence the result. $\bullet$

\proclaim Calculation 5.4. $\omega$ is a positive form of type $(1,1)$
if and only if $\beta$ is symmetric and 
$\im\beta$ is positive definite.

Proof. We first 
examine the condition that
$\omega$ is of type $(1,1)$,
i.e. we can write
$$\omega={i\over 2}\sum_{i,j} h_{ij} \theta_i\wedge\bar\theta_j.$$
Now
$$\eqalign{\theta_i\wedge\bar\theta_j=&dx_i\wedge dx_j+dx_i\wedge \left(\sum_k
\bar\beta_{jk} dy_k\right)-dx_j\wedge \left(\sum_k\beta_{ik}dy_k\right)\cr
&+\sum_{k,l}\beta_{ik}\bar\beta_{jl}dy_k\wedge dy_l\cr}$$
and
$$\eqalign{\sum h_{ij}\theta_i\wedge\bar\theta_j
=&\sum_{i<j} (h_{ij}-h_{ji})dx_i\wedge dx_j
+\sum_i dx_i\wedge\left(\sum_{j,k} h_{ij}\bar\beta_{jk}dy_k
-\sum_{j,k} h_{ji}\beta_{jk}dy_k\right)\cr
&+\sum_{i,j,k,l} h_{ij}\beta_{ik}\bar\beta_{jl} dy_k\wedge dy_l.\cr}$$

Thus in order for $\omega=\sum dx_i\wedge dy_i=
{i\over 2}\sum h_{ij}\theta_i\wedge\bar\theta_j$, we must
have, in particular, $h_{ij}=h_{ji}$, i.e. $h$ is symmetric.  On
the other hand, since $\omega$ is real, $h_{ij}=\bar h_{ji}$, and
thus the matrix $h$ is real. Also,
$$\eqalign{I&={i\over 2}(h\bar\beta-h\beta)\cr
&=h\im\beta.\cr}$$
Thus we have $h=(\im\beta)^{-1}$, and $\im\beta$ is symmetric.

To ensure the last term vanishes, we need
$$\sum_{i,j} h_{ij}\beta_{ik}\bar\beta_{jl}
=\sum_{i,j} h_{ij}\beta_{il}\bar\beta_{jk},$$
or equivalently, the matrix ${}^t\beta h\bar\beta$
is symmetric.

Now
$$\eqalign{{}^t\beta h\bar\beta&=(\re{}^t\beta+i\im\beta)(\im\beta)^{-1}
(\re\beta-i\im\beta)\cr
&=(\re{}^t\beta)(\im\beta)^{-1}(\re\beta)+i\re\beta-i\re{}^t\beta+\im\beta,\cr}
$$
while
$$\eqalign{{}^t\bar\beta h\beta&=(\re{}^t\beta-i\im\beta)(\im\beta)^{-1}
(\re\beta+i\im\beta)\cr
&=(\re{}^t\beta)(\im\beta)^{-1}(\re\beta)+i\re{}^t\beta-
i\re\beta+\im\beta,\cr}
$$
so symmetry of ${}^t\beta h\bar\beta$ is equivalent to 
$\re{}^t\beta=\re\beta$. 

Thus $\omega$ is of type $(1,1)$ if and only if $\beta$ is symmetric.
In addition, to ensure $\omega$ is a positive $(1,1)$ form,
$h=(\im\beta)^{-1}$ must be positive definite, so $\im\beta$ must
be positive definite.
$\bullet$

The real problem is understanding the condition $d\Omega=0$.
This is the heart of the difficulty, and we will return to this shortly.

We first connect $\Omega$ to the description of the choice of almost
complex structure given in the introduction, namely as a choice of
horizontal subspaces of an Ehresmann connection and the choice of a metric
on the fibres. 

\proclaim Proposition 5.5. The matrix $(\im\beta)^{-1}$ is the matrix
of the metric $(g_{ij})$ on the fibres of $f$. For a point $x\in X_b^{\#}$,
$J(\T_{X_b,x})$ is spanned by the tangent vectors
$$\{\partial/\partial y_j-\sum_i\re\beta_{ij}\partial/\partial x_i|1\le j\le n
\},$$ 
where $J:\T_{X^{\#}}\rightarrow \T_{X^{\#}}$ is the almost complex structure
induced by $\Omega$.

Proof. Since $\omega={i\over 2}\sum_{i,j} h_{ij}\theta_i\wedge\bar\theta_j$
with $h=(\im\beta)^{-1}$ is the K\"ahler form of the metric, the K\"ahler
metric itself is $g=\sum_{i,j} h_{ij} \theta_i\otimes\bar\theta_j$.
Thus $g_{ij}=g(\partial/\partial x_i,\partial/\partial x_j)
=h_{ij}$, giving the interpretation of $\im\beta$.

Next, let $J$
be the almost complex structure on $\T_{X^{\#}}$ induced by $\Omega$, and
${}^tJ$ the almost complex structure on $\T_{X^{\#}}^*$.
Since the space spanned by
$\theta_1,\ldots,\theta_n$ is the $+i$ eigenspace of ${}^tJ$ at a point
$x\in X^{\#}$, the cotangent space $\T_{X,x}^*$ decomposes
as $V_1\oplus V_2$ with 
$$\eqalign{V_1&=span(\re\theta_1,\ldots,\re\theta_n)\cr
V_2&=span(\im\theta_1,\ldots,\im\theta_n),\cr}$$
with ${}^tJ(V_1)=V_2$ and ${}^tJ(V_2)=V_1$. Note that $V_2=span(dy_1,\ldots,
dy_n)$ since $\im\beta$ is invertible. Thus also $\T_{X,x}=V_1^{\circ}
\oplus V_2^{\circ}$, $V_i^{\circ}$ the annihilator of $V_i$, with
$J$ interchanging $V_1^{\circ}$ and $V_2^{\circ}$. Now 
$$\eqalign{V_2^{\circ}&=span(\partial/\partial x_1,\ldots,\partial/\partial 
x_n)\cr
V_1&=span(\{dx_i+\sum_j\re\beta_{ij} dy_j\})\cr
V_1^{\circ}&=span(\{\partial/\partial y_j-\sum_i\re\beta_{ij}
\partial/\partial x_i\}).\cr}$$ 
Thus we see that $\re\beta$ determines $J(\T_{X_b,x})=V_1^{\circ}$ as 
claimed. $\bullet$

Summarizing, we now have

\proclaim Theorem 5.6. Specifying an $n$-form $\Omega$
on $X^{\#}$ satisfying
properties (2) and (3) of Theorem 5.1 and properties (2) and (3) of 
Theorem 5.2 and such that $f^{\#}$ is a special Lagrangian fibration 
with respect to $\Omega$ is equivalent to specifying
\item{(1)} A metric $(g_{ij})$ on each fibre of $f^{\#}:X^{\#}\rightarrow
B$.
\item{(2)} A splitting $\T_{X^{\#}}=\T_{X^{\#}/B}\oplus \F$, where
$\T_{X^{\#}/B}$ is the subbundle of $\T_{X^{\#}}$ with $\T_{X^{\#}/B,x}
=\T_{X_b^{\#},x}$, and $\F$ is a Lagrangian subbundle of $\T_{X^{\#}}$.

Proof. Proposition 5.5 shows that $\Omega$ specifies the metric on the
fibres and a splitting as desired with $\F=J(\T_{X^{\#}/B})$.
This is clearly a Lagrangian subbundle since $\T_{X^{\#}/B}$ is.

Conversely, giving a splitting of the exact sequence
$$0\mapright{}\T_{X^{\#}/B}\mapright{}
\T_{X^{\#}}\mapright{p}\F\mapright{}0$$ determines $\re\beta$ for us.
Indeed, at a point $x\in X^{\#}$, with local coordinates as usual,
such a splitting gives a map $s:\F_x\rightarrow 
\T_{X^{\#},x}$, and there is a matrix $(b_{ij})$
such that 
$$s(p(\partial/\partial y_j))=\partial/\partial y_j
-\sum b_{ij} \partial/\partial x_i.$$
We take $\re\beta_{ij}=b_{ij}$. Note that the symmetry of the
matrix $b_{ij}$ is equivalent to $s(\F_x)$
being Lagrangian. 
Thus specifying (1) and (2) is equivalent, in local
coordinates, to giving $\beta_{ij}=b_{ij}+ig^{ij}$, where $g_{ij}$ is
the metric on the fibre. Then we must have, with
$$\theta_i=dx_i+\sum_j\beta_{ij}dy_j,$$
$$\Omega=V \theta_1\wedge\cdots\wedge\theta_n$$
for some real function $V$.
Calculation 5.4 then tells us that $\omega$ is a positive $(1,1)$-form
in the corresponding almost complex structure, since by construction $\beta$
is symmetric and $\im\beta$ is positive definite. Calculation 5.3
shows that
$\omega^n/n!=(-1)^{n(n-1)/2}(i/2)^n\Omega\wedge\bar\Omega$ 
if and only if $V=\sqrt{\det(g_{ij})}$
since
$(\sqrt{\det(g_{ij})})^2\det(g^{ij})=1$. $\bullet$
\bigskip

We have now seen how the $n$-form $\Omega$ can be determined by choosing
a matrix $\beta=(\beta_{ij})=(b_{ij}+ig^{ij})$, so that $$\Omega
=\sqrt{\det(g_{ij})}\bigwedge_i(dx_i+\sum_j \beta_{ij} dy_j).$$
If $\beta$ is chosen to be symmetric and $\im\beta$ positive definite,
then
we have seen that $d\Omega=0$ implies both the integrability of the
almost complex structure induced by $\Omega$ and the Ricci-flatness of
the metric induced by this complex structure and the standard symplectic
form $\omega$. At first sight, the condition $d\Omega=0$ looks very
complicated if one proceeds by brute force and tries to compute
the exterior derivative of $\Omega$. In fact, this condition can be
simplified, and we wish to examine this here. We note that the calculation
below of $d\Omega$ is quite similar to calculations carried out
in [26] and [27] for an analagous situation in the study
of deformations of complex structures on Calabi-Yau manifolds. However,
we will introduce a formalism using differential operators to make
the calculation easier.

To begin, first note that the integrability of the almost complex
structure determined by $\Omega$ is a weaker condition than
$d\Omega=0$. Let us first understand this weaker condition.

\proclaim Theorem 5.7. The almost complex structure induced by $\Omega$ is
integrable if and only if 
$$\sum_i \left( {\partial \beta_{lk}\over\partial x_i}\beta_{ij}
-{\partial\beta_{lj}\over\partial x_i}\beta_{ik}\right)
={\partial \beta_{lk}\over\partial y_j}-{\partial\beta_{lj}\over
\partial y_k}.$$
for all $j,k$ and $l$.

Proof.
Writing as before
$$\Omega=V\theta_1\wedge\cdots\wedge\theta_n,$$
the almost complex structure is determined by the fact that 
$\theta_1,\ldots,\theta_n$ should span the space of $(1,0)$ forms.
To show that the almost complex structure induced by $\theta_1,\ldots,
\theta_n$ is integrable we need to show that $d\theta_i$ is of type
$(2,0)+(1,1)$ for all $i$. Since
$$\theta_i=dx_i+\sum_{j=1}^n\beta_{ij}dy_j,$$
$\theta_1,\ldots,\theta_n,dy_1,\ldots,dy_n$ form a basis for
the space of 1-forms, and thus we can write
$$d\theta_l=\sum_{i,j} {1\over 2}A_{ij}^l\theta_i\wedge\theta_j
+\sum_{ij} B_{ij}^l\theta_i\wedge dy_j
+\sum_{i,j} {1\over 2}C_{ij}^l dy_i\wedge dy_j.$$
The almost complex structure is integrable if and only if
$C_{ij}^l=0$ for all $i,j$ and $l$. Here $A^l$ and $C^l$ are skew-symmetric
matrices.
Note
$$d\theta_l=\sum_{i,j}{\partial \beta_{lj}\over\partial x_i}
dx_i\wedge dy_j+\sum_{i,j}{\partial \beta_{lj}\over\partial y_i}
dy_i\wedge dy_j.$$
Since $d\theta_l$ contains no $dx_i\wedge dx_j$ terms, we must have
$A^l_{ij}=0$, and then $B^l_{ij}=\partial \beta_{lj}/\partial x_i$. Then
the almost complex structure is integrable if and only if
$$\eqalign{d\theta_l&=\sum_{i,j} B^l_{ij} \theta_i\wedge dy_j\cr
&=\sum_{i,j} {\partial \beta_{lj}\over \partial x_i} dx_i\wedge dy_j
+\sum_{i,j,k} {\partial \beta_{lj}\over\partial x_i} \beta_{ik}
dy_k\wedge dy_j\cr}$$
which holds if and only
$$\sum_i \left( {\partial \beta_{lk}\over\partial x_i}\beta_{ij}
-{\partial\beta_{lj}\over\partial x_i}\beta_{ik}\right)
={\partial \beta_{lk}\over\partial y_j}-{\partial\beta_{lj}\over
\partial y_k}.$$ This is the desired condition. $\bullet$

We wish to rephrase this condition. 
We are going 
work locally for the moment, fixing coordinates $y_1,\ldots,y_n$
on an open subset $U\subseteq B$, and consider all forms as living
on $U\times{\bf R}^n$ with coordinates $y_1,\ldots,y_n,x_1,\ldots,x_n$.
We then write
$$\beta=\sum_{i,j} \beta_{ij} dy_j\otimes{\partial\over\partial x_i}
\in \Gamma(U\times {\bf R}^n,f^*\Omega^1_B\otimes\T_{X/B}).$$
(This is not a coordinate independent expression. The correct coordinate
independent expression would be
$$\sum_idx_i\otimes {\partial\over\partial x_i}+\beta\in \Gamma(f^{-1}(U),
\Omega^1_X\otimes \T_{X/B}),$$
but for practical purposes it is more convenient to work with the above
expression.)

\proclaim Definition 5.8.
For expressions of the type $v=\sum_j dy_j\otimes v_j$, $w=\sum_j dy_j\otimes 
w_j$ with $v_j$, $w_j$ vector fields, $v_j=\sum_i v_{ij}\partial/\partial x_i$,
$w_j=\sum_i w_{ij}\partial/\partial x_i$, we define
$$[v,w]=\sum_{l,m} [v_l,w_m] dy_l\wedge dy_m$$
where
$$[v_l,w_m]=\sum_{i,j} \left( v_{il}{\partial w_{jm}\over\partial x_i}
-w_{im}{\partial v_{jl}\over\partial x_i}\right){\partial\over
\partial x_j}$$
is the usual Lie bracket of vector fields. We also set
$$d_yv=\sum_{i,j,k} {\partial v_{ij}\over\partial y_k}dy_k\wedge dy_j
\otimes {\partial\over\partial x_i}.$$ 

We then obtain

\proclaim Theorem 5.9. The almost complex structure induced by $\Omega$ is
integrable if and only if
$$d_y\beta-{1\over 2}[\beta,\beta]=0.$$

Proof. $$d_y\beta=\sum_{l\atop i<j}\left({\partial \beta_{lj}\over
\partial y_i}-{\partial\beta_{li}\over\partial y_j}\right)dy_i\wedge dy_j
{\partial\over \partial x_l}$$
while
$$[\beta,\beta]=\sum_{i,j\atop l<m} 
2\left(\beta_{il}{\partial\beta_{jm}\over\partial x_i}
-\beta_{im}{\partial\beta_{jl}\over\partial x_i}\right) dy_l\wedge dy_m
{\partial\over\partial x_j}.$$
Comparing with the formula of Theorem 5.7, we see that the two formulae
are equivalent. $\bullet$

Our next goal is to prove

\proclaim Theorem 5.10. $d\Omega=0$ if and only if the almost complex structure
induced by $\Omega$ is integrable and $d\Omega\subseteq F^2\Omega_X^{n+1}$.

Thus, locally, one just needs to check that the coefficients of 
$dy_i\wedge dx_1\wedge\cdots\wedge dx_n$, $1\le i\le n$ in $d\Omega$ are
zero, and also check the integrability condition. 

To prove this theorem, we must introduce some additional algebraic structure
to accomplish the calculation. We continue to use a choice of local 
coordinates, and work on the space $U\times{\bf R}^n$, $f:U\times{\bf R}^n
\rightarrow U$ the projection. Let $\T_x$ be the subbundle of the tangent
bundle of $U\times{\bf R}^n$ generated by $\partial/\partial x_1,\ldots,
\partial/\partial x_n$, and let $\Omega_y$ be the subbundle of the cotangent
bundle generated by $dy_1,\ldots,dy_n$. For $q\ge 0,p\le 0$, set
$$\Omega^{p,q}_U=\Gamma(\bigwedge^q\Omega_y\otimes\bigwedge^{-p}\T_x)$$
and set $\Omega_0=dx_1\wedge\cdots\wedge dx_n$. Then there is an isomorphism
$$\bigoplus_{i=0}^p\Omega_U^{p-n-i,i}\mapright
{\cong}\Gamma(\Omega^p_{U\times{\bf R}^n})\leqno{(5.1)}$$
where $\Omega^p_{U\times {\bf R}^n}$ denotes the sheaf of $C^{\infty}$ $p$-forms
on $U\times{\bf R}^n$. This isomorphism sends $\theta\otimes v$ to 
$\theta\wedge \iota(v)\Omega_0$. In particular, $\Gamma(\Omega^n_{U\times
{\bf R}^n})\cong \bigoplus_{p=0}^n\Omega_U^{-p,p}$.

For $\alpha\otimes\beta\in\Omega^{p,q}_U$, $\alpha'\otimes\beta'\in
\Omega^{p',q'}_U$, we can define the product
$$(\alpha\otimes\beta)\cdot(\alpha'\otimes\beta'):=(\alpha\wedge \alpha')
\otimes (\beta\wedge\beta')\in\Omega^{p+p',q+q'}_U.$$
This satisfies the commutation relations
$$(\alpha\otimes\beta)\cdot(\alpha'\otimes\beta')=
(-1)^{pp'+qq'}(\alpha'\otimes\beta')\cdot(\alpha\otimes\beta).$$
This gives us a bigraded ring structure on $\bigoplus\Omega_U^{p,q}$.
Note that the subring $\bigoplus_{p=0}^n \Omega_U^{-p,p}\cong
\Gamma(\Omega^n_{U\times{\bf R}^n})$
is in fact a commutative ring with $1$, and $1\in \Omega_U^{0,0}$ corresponds
to $\Omega_0$ under this isomorphism.

\proclaim Lemma 5.11. Under the isomorphism (5.1),
$$\bigwedge_{i=1}^n\left( dx_i+\sum_j \beta_{ij}dy_j\right)
=exp(\beta):=\sum_{p=0}^{\infty} \beta^p/p!$$ 
where
$$\beta=\sum_{i,j} \beta_{ij} dy_j\otimes {\partial\over\partial x_i}$$
in the notation introduced above.

Proof. This is a straightforward though slightly tedious calculation.
Here is where all the signs must be dealt with correctly. For
this and subsequent calculations it is convenient to keep in mind that
$\iota(\partial/\partial x_I)\Omega_0
=(-1)^{M} dx_{I^*}$, where
$I^*=\{1,\ldots,n\}-I$ and $M=\#\{(i,j)| i\in I, j\in I^*,i>j\}$. $\bullet$

The next step is to turn $\Omega_X^{\cdot,\cdot}$ into a double complex.
We have exterior differentiation $$d:\Gamma(\Omega^i_{U\times{\bf R}^n})
\rightarrow\Gamma(\Omega^{i+1}_{U\times{\bf R}^n}),$$
and under the isomorphism (5.1), it is clear that $d(\Omega^{p,q}_U)
\subseteq \Omega_U^{p+1,q}\oplus\Omega_U^{p,q+1}$. Thus
we can write $d=d_x+d_y$, so that
$(\Omega_U^{\cdot,\cdot},d_x,d_y)$ defines a bicomplex, with
$$d_x:\Omega_U^{p,q}\rightarrow \Omega_U^{p+1,q},$$
$$d_y:\Omega_U^{p,q}\rightarrow \Omega_U^{p,q+1}.$$
One checks that 
$$d_x(dy_J\otimes \alpha)=(-1)^q\sum_{i=1}^n dy_J\otimes({\partial\alpha\over
\partial x_i}
\relbow dx_i).\leqno{(5.2)}$$
(Here, $(\partial/\partial x_I\wedge \partial/\partial x_i)\relbow dx_i
=\partial/\partial x_I$, and $\partial/\partial x_I\relbow dx_i=0$ if $i\not
\in I$.)

Let $D:\Omega_U\rightarrow\Omega_U$ be a graded endomorphism of $\Omega_U:=
\Omega_U^{\cdot,\cdot}$. We now recall what it means for $D$ to be
a differential operator of order $\le r$. Put on $\Omega_U\otimes\Omega_U$
the anti-commutative algebra structure given by
$$(a\otimes b)\cdot(a'\otimes b')=(-1)^{(\deg a')\cdot (\deg b)}
aa'\otimes bb'.$$ 
Here the dot is the standard dot product, keeping in mind $\Omega_U$
is bigraded. This turns $\Omega_U\otimes\Omega_U$ into a bigraded anti-commutative
algebra. Let $\lambda:\Omega_U\rightarrow\Omega_U\otimes\Omega_U$ be
given by
$\lambda(a)=a\otimes 1 - 1 \otimes a$. We define
$$\Phi_D^r:\Omega_U^r\rightarrow \Omega_U$$
by 
$$\Phi_D^r(a_1,\ldots,a_r)=m\circ (D\otimes id_{\Omega_U})(\prod_{i=1}^r
\lambda(a_i)).$$
Here $m(a\otimes b)=ab$. We say
$D$ is a differential operator of order $\le r$ if $\Phi^{r+1}_D$ is 
identically zero. Note that
$$\Phi^2_D(a,b)=D(ab)-D(a)b-(-1)^{(\deg a)\cdot(\deg b)}D(b)a
+D(1)ab$$
and
$$\Phi^3_D(a,b,c)=\Phi^2_D(a,bc)-\Phi^2_D(a,b)c-(-1)^{(\deg b)\cdot(\deg c)}
\Phi^2_D(a,c)b.$$
(See [19], \S 1.) As usual, the composition of differential operators
of orders $\le r$ and $s$ is order $\le r+s$.

\proclaim Definition 5.12. 
We set $d_x':\Omega_U\rightarrow\Omega_U$ to be the operator
acting on $\Omega_U^{p,q}$ by $(-1)^{p+q+1}d_x$.

\proclaim Lemma 5.13. 
$d_x$ is a differential operator of
order $\le 2$ and $d_y$ is a differential operator of order $\le 1$.

Proof. That $d_y$ is a differential operator of order $\le 1$ follows
immediately from the definition, while (5.2) shows that $d_x'$ can
be written as a sum of a composition of two operators: differentiation
in the direction $\partial/\partial x_i$ 
and $\alpha\otimes\beta\mapsto \alpha
\otimes(-1)^{p+1}\beta\relbow dx_i$. 
One checks easily that these are each first
order operators, and hence $d_x$ is second order. $\bullet$

Now define, for $\alpha,\beta\in\Omega_X$,
$$[\alpha,\beta]:=\Phi^2_{d'_x}(\alpha,\beta),$$
so
$$d'_x(\alpha\beta)=[\alpha,\beta]+d'_x(\alpha)\beta+(-1)^{(\deg
\alpha)\cdot(\deg\beta)} d'_x(\beta)\alpha.\leqno{(5.3)}$$
Since $\Phi^3_{d'_x}=0$, we obtain
$$[\alpha,\beta\gamma]=[\alpha,\beta]\gamma +(-1)^{(\deg\beta)\cdot(\deg\gamma)}
[\alpha,\gamma]\beta, \leqno{(5.4)}$$
and from Lemma 5.13,
$$d_y(\alpha\beta)=d_y(\alpha)\beta+(-1)^{(\deg\alpha)\cdot(\deg\beta)}d_y(\beta)
\alpha. \leqno{(5.5)}$$

We note that this definition of bracket is an extension of 
Definition 5.8:

\proclaim Proposition 5.14. If $\alpha,\beta\in
\Omega^{-1,1}$, then $[\alpha,\beta]$ as defined above agrees with
the earlier definition.

Proof. By linearity, we can assume that $\alpha=f dy_j\partial/\partial x_i$
and $\beta=g dy_l\partial/\partial x_k$. Then 
\def\dxi{{\partial\over\partial x_i}}
\def\dxk{{\partial\over\partial x_k}}
$$\eqalign{[\alpha,\beta]=&-d_x(\alpha\beta)+d_x(\alpha)\beta+
d_x(\beta)\alpha\cr
=&-d_x(fgdy_j\wedge dy_l\dxi\wedge\dxk)+d_x(fdy_j\dxi)gdy_l\dxk\cr
&+d_x(gdy_l\dxk)fdy_j\dxi\cr
=&-\left({\partial f\over\partial x_k}g+f{\partial g\over \partial x_k}\right)
dy_j\wedge dy_l\dxi+
\left({\partial f\over\partial x_i}g+f{\partial g\over \partial x_i}\right)
dy_j\wedge dy_l\dxk\cr
&-g{\partial f\over\partial x_i}dy_j\wedge dy_l\dxk
+f{\partial g\over\partial x_k} dy_j\wedge dy_l\dxi\cr
=&
f{\partial g\over\partial x_i}dy_j\wedge dy_l \dxk
-
{\partial f\over\partial x_k}g dy_j\wedge dy_l\dxi,
\cr}$$
proving the desired equality. $\bullet$
\bigskip

{\it Proof of Theorem 5.10.} 
First suppose that $d\Omega=0$. Then the almost complex structure is
integrable by Theorem 5.1, and obviously $d\Omega\subseteq F^2\Omega_X^{n+1}$.

Conversely, suppose $d\Omega\subseteq F^2\Omega_X^{n+1}$ and the almost complex
structure is integrable. We use the isomorphism (5.1) and Lemma 5.11 to write
$\Omega=V\exp(\beta)$ for some $\beta\in \Omega^{-1,1}$, $V\in\Omega^{0,0}$.
To show that $d\Omega=0$, we need to show that each graded piece of $d\Omega$
is zero, i.e.
$${1\over n!}d_y(V\beta^n)+{1\over (n+1)!}d_x(V\beta^{n+1})=0,\quad
\hbox{for all $n\ge 0$.} $$
This is equivalent to
$$d_y(V\beta^n)={1\over n+1} d'_x(V\beta^{n+1}). \leqno{(5.6)}$$
The fact that $d\Omega\subseteq
F^2\Omega_X^{n+1}$ is equivalent to (5.6) for $n=0$. Note that this
states
$$d_y(V)=d_x'(V\beta)=[\beta,V]+d_x'(\beta)V\leqno{(5.7)}$$
by (5.3), since $d_x'(V)=0$.
Observe that since $[\beta,\cdot]$ acts as an ordinary (non-graded) derivation
on the commutative ring $\bigoplus_p \Omega^{-p,p}$, we have $[\beta,\beta^n]
=n[\beta,\beta]\beta^{n-1}$.

We prove by induction that (5.7) implies
$$d_x'(V\beta^{n+1})=
(n+1)d_y(V)\beta^n+{n(n+1)\over 2}[\beta,\beta]V\beta^{n-1}.
\leqno{(5.8)}$$
Indeed this is true for $n=0$, by (5.7). Then
$$\eqalign{d_x'(V\beta^{n+1})&=d_x'(\beta\cdot V\beta^n)\cr
&=[\beta,V\beta^n]+d_x'(\beta)V\beta^n+d_x'(V\beta^n)\beta\cr
&=[\beta,V]\beta^n+n[\beta,\beta]V\beta^{n-1}+d_x'(\beta)V\beta^n
+d_x'(V\beta^n)\beta\cr
&=d_y(V)\beta^n+n[\beta,\beta]V\beta^{n-1}+d_x'(V\beta^n)\beta\cr}$$
by (5.7), so by induction the desired result holds.

Next, we note using the integrability condition $d_y\beta=[\beta,\beta]/2$
and (5.8) that
$$\eqalign{d_y(V\beta^n)&=d_y(V)\beta^n+d_y(\beta^n)V\cr
&=d_y(V)\beta^n+nd_y(\beta)V\beta^{n-1}\cr
&=d_y(V)\beta^n+{n\over 2}[\beta,\beta]V\beta^{n-1}\cr
&={1\over n+1}d_x'(V\beta^{n+1}).\cr}$$
This proves (5.6), and hence the theorem. $\bullet$

Next I would like to reinterpret the equations we've seen above
so that they may look more natural. As we have seen earlier,
$b=\re\beta$ defines an Ehresmann connection whose horizontal subspaces,
given by the subbundle $\F$, are Lagrangian. In local coordinates,
this connection is determined by 
$$b=\sum_{i,j} b_{ij}dy_j\otimes {\partial\over \partial x_i}=\re \beta.$$
It makes sense to define the covariant derivative with respect to this
connection. This will be an operator
$$\nabla_b:\Gamma(f^*\Omega^q_B\otimes\T_{X^{\#}/B})
\rightarrow \Gamma(f^*\Omega_B^{q+1}\otimes\T_{X^{\#}/B})$$
defined by
$$\nabla_b\alpha:=d_y\alpha-[b,\alpha].$$
It is easy to check that this definition is now independent of the choice
of coordinates. The curvature tensor of $\nabla_b$ is then
$F_b\in\Gamma(f^*\Omega^2_B\otimes\T_{X^{\#}/B})$ given by
$$F_b:=d_yb-{1\over 2}[b,b].$$
It is easy to check that $F_b=0$ if and only if the horizontal
distribution $\F$ is integrable. Of course, an Ehresmann connection
gives rise to parallel transport along a path contained in $B_0$;
we say a family of $p$-forms on the fibres of $f_0:X_0\rightarrow B_0$
is {\it parallel} if it is invariant under parallel transport.
If in local coordinates over $U\subseteq B_0$ this family of forms is 
written as $\alpha=\sum_I f_Idx_I$, $f_I$ a function on $f^{-1}(U)$,
$\alpha$ is parallel if
$$d_y\alpha-\L_b\alpha=0;$$
by this we mean
$${\partial\alpha\over\partial y_j}-\L_{\sum_i b_{ij}\partial/\partial x_i}
\alpha=0$$
for each $j$.
(For a similar treatment of Ehresmann connections, see [20]).

We can now rephrase the integrability conditions in a more invariant way.

\proclaim Corollary 5.15. Let $(\beta_{ij})=(b_{ij}+ig^{ij})$, so
we write $\beta=b+ig^{-1}$, $V=\sqrt{\det g}$. Then $d\Omega=0$ if and only if
$$F_b+{1\over 2}[g^{-1},g^{-1}]=0\leqno{(5.9)}$$
$$\nabla_bg^{-1}=0\leqno{(5.10)}$$
$$\hbox{$dx_1,\ldots,dx_n$ are harmonic forms on each fibre}\leqno{(5.11)}$$
$$\hbox{$Vdx_1\wedge \cdots\wedge dx_n$ is parallel.}\leqno{(5.12)}$$

Proof. The first two equations are the real and imaginary parts of
$d_y\beta-{1\over 2}[\beta,\beta]=0$. The last two are (5.7) broken
up again into its real and imaginary parts. $\bullet$

We remark here that in the study of Ricci curvature in the context
of Riemannian submersions, some similar structures arise. See [8], Chapter
9.

\proclaim Corollary 5.16. Suppose $d\Omega=0$ and $\nabla_b$ is
flat. Then the metric $g$ is flat along the fibres.

Proof. By (5.9) we have $[g^{-1},g^{-1}]=0$. We work on one fixed fibre
$X_b$ with coordinates $x_1,\ldots, x_n$. Let $g^j$ be the vector field
$\sum_i g^{ij}\partial/\partial x_i$, so $g^1,\ldots,g^n$ form a basis for
$\T_{X_b}$ at each point of $X_b$, and $[g^i,g^j]=0$. Let $\omega_1,\ldots,
\omega_n$ be the dual basis of one-forms. Then
$$\eqalign{
d\omega_i(g^j,g^k)&=g^j(\omega_i(g^k))-g^k(\omega_i(g^j))-\omega_i([g^j,g^k])
\cr
&=0\cr}$$
since $\omega_i(g^k)=\delta_{ik}$ is constant.
Thus $\omega_1,\ldots,\omega_n$ are closed 1-forms, and so $\omega_i=d\varphi_i$
for a function $\varphi_i$ on ${\bf R}^n$, the universal cover of $X_b$.
$\varphi_i$ will be a linear function plus a periodic function. Since
$\omega_i=\sum_i g_{ij}dx_j$, we see that $g_{ij}=\partial\varphi_i
/\partial x_j=g_{ji}=\partial \varphi_j/\partial x_i$, so there exists
a function $\varphi$ on ${\bf R}^n$ such that $\varphi_i=\partial\varphi/
\partial x_i$, and $g_{ij}=\partial^2\varphi/\partial x_i\partial x_j$.
$\varphi$ satisfies the real inhomogeneous Monge-Amp\`ere equation
$\det(\partial^2\varphi/\partial x_i\partial x_j)=V^2$.

Next we prove $V$ is constant. Note that $V^{-1}\omega_1\wedge\cdots
\wedge \omega_n=Vdx_1\wedge\cdots\wedge dx_n$, the volume form on 
$X_b$, so $*dx_i=\iota(g^i)(V^{-1}\omega_1\wedge\cdots\wedge \omega_n)
=\pm V^{-1}\omega_1\wedge\cdots\wedge \hat{\omega_i}\wedge\cdots\wedge
\omega_n$, and so $d(*dx_i)=\pm g^i(V^{-1})\omega_1\wedge\cdots\wedge\omega_n$.
Thus, by (5.11), we must have $g^i(V^{-1})=0$ for all $i$ so $V^{-1}$,
hence $V$, is constant.

Thus $\varphi$ is a solution to the equation
$$\det\left({\partial^2\varphi\over \partial x_i\partial x_j}\right)=C,$$
$C$ a constant, and $\varphi$ is a function on ${\bf R}^n$, with
$\varphi=\varphi_{quad}+\varphi_{lin}+\varphi_{per}$, the decomposition
into a quadratic, linear, and periodic part. We can of course assume
that $\varphi_{lin}=0$ and $\int_{X_b}\varphi_{per} 
dx_1\wedge\cdots\wedge dx_n=0$. Then I claim $\varphi=\varphi_{quad}$,
and so $g_{ij}=\partial^2\varphi/\partial x_i\partial x_j$ is constant.
To see this, one applies a standard technique for non-linear elliptic
partial differential equations. Let $\varphi_t=t\varphi+(1-t)\varphi_{quad}$.
Let $m_{ij}(x_{kl})$ denote the $ij$th cofactor of the matrix $(x_{kl})$,
so in particular
$$m_{ij}(x_{kl})={\partial \det(x_{kl})\over \partial x_{ij}}.$$
Put 
$$a_{ij}=\int_0^1 m_{ij}(\partial^2\varphi_t/\partial x_k\partial x_l) dt.$$
Now $(\partial^2\varphi/\partial x_i\partial x_j)$ is positive definite.
So in fact is $(\partial^2\varphi_{quad}/\partial x_i\partial x_j)$. 
To see this, note that if we
put $h_{ij}=h(\partial/\partial y_i,\partial
/\partial y_j)$, we calculate
$$\eqalign{h_{ij}&=-\int_{X_b} \iota(\partial/\partial y_i)\omega
\wedge \iota(\partial/\partial y_j)\im\Omega\cr
&=\int_{X_b} dx_i\wedge \left( V\sum_k (-1)^{k-1} g^{kj} dx_1
\wedge\cdots\wedge\hat{dx_k}\wedge\cdots\wedge dx_n\right)\cr
&=\int_{X_b} Vg^{ij} dx_1\wedge\cdots\wedge dx_n.\cr}$$
Since $V$ is constant, $h_{ij}=V{\partial^2\varphi_{quad}\over\partial x_i
\partial x_j}\int_{X_b}dx_1\wedge\cdots\wedge dx_n$. 
Thus since $(h_{ij})$ is positive definite,
so is $(\partial^2\varphi_{quad}/\partial x_i\partial x_j)$.
Thus $(\partial^2\varphi_t/\partial x_i\partial x_j)$ is positive
definite for all $0<t<1$, so the matrix $(a_{ij})$ is also positive
definite.

{\it Claim:} 
$$\sum_{i,j} a_{ij}{\partial^2(\varphi-\varphi_{quad})\over \partial x_i
\partial x_j}=constant.$$

Proof.
$$\eqalign{{\partial\over\partial t}\left(\det\left({\partial^2
\varphi_t\over\partial x_k\partial x_l}\right)\right)&=
\sum_{i,j} m_{ij}\left({\partial^2\varphi_t\over
\partial x_k\partial x_l}\right){\partial^3\varphi_t\over
\partial x_i\partial x_j\partial t}\cr
&=\sum_{i,j} m_{ij}\left({\partial^2\varphi_t\over\partial x_k
\partial x_l}\right) {\partial^2(\varphi-\varphi_{quad})\over
\partial x_i\partial x_j}.\cr}$$
Integrating with respect to $t$ gives the desired result. $\bullet$

Now $\varphi-\varphi_{quad}$ is a periodic function, so applying the 
maximum principal (or minimum principal, depending on
the sign of the constant), $\varphi-\varphi_{quad}$ is
constant. Hence $\varphi=\varphi_{quad}$. $\bullet$

This explains why in [17], the integrable complex structures constructed
on torus fibrations had to have flat metric on the fibres
given that $\re\beta$ was taken to be zero in that paper.
\bigskip

{\hd \S 6. The complex structure on the mirror.}

Having understood, at least to some extent, how one describes
complex structures on torus fibrations, we now wish to explain
how one should put a complex structure on the $D$-brane moduli space.
We cannot solve this problem at present due to the complexity of
the equation $d\Omega=0$; however, here we will give guidelines
as to where to look for the correct solutions.

We continue with an integral
special Lagrangian torus fibration $f:X\rightarrow B$
with a Lagrangian section, along with forms $\omega$ and $\Omega$. In \S 4, we
have seen how to put a symplectic form $\check\omega$ on 
$\check f:\check X_0\rightarrow B$ which has the property that
$[\check\omega]$ and $[\im\Omega_n]$ agree in $\Gamma(B_0,
\Omega^1_{B_0}\otimes R^1\check f_{0*}{\bf R})
\cong \Gamma(B_0, \Omega^1_{B_0}\otimes R^{n-1}f_{0*}{\bf R})$.
To specify the complex structure on $\check X$, we need to construct
the form $\check\Omega$. This should, according to the appropriate
conjectures, be determined by the $B$-field, i.e. something like an element
$\b\in H^2(X,{\bf R}/\boldz)$, and $\omega$ the K\"ahler form on $X$.
Certainly the first requirement for $\check\Omega$ should be that
$[\omega]$ and $[\im\check\Omega_n]$ should agree on
$\Gamma(B_0,\Omega^1_{B_0}\otimes R^1f_{0*}{\bf R})\cong
\Gamma(B_0,\Omega^1_{B_0}
\otimes R^{n-1}\check f_{0*}{\bf R})$, so that the double dual
brings us back to $X$. The second requirement should involve
$\re\check\Omega_n$ and is much less precise at this point.
We can only be guided by item (5) of the introduction, but will try to be
more precise later.

\proclaim Proposition 6.1. If $\check\Omega_n$ is a normalised holomorphic
$n$-form on $\check X_0$ making $\check f_0:\check X_0\rightarrow B$
special Lagrangian
and if $\check h_n$ is the induced normalised Riemannian metric on the base,
then
$[\omega]=[\im\check\Omega_n]$ in
$\Gamma(B_0,\Omega^1_{B_0}\otimes R^1f_{0*}{\bf R})
\cong\Gamma(B_0,\Omega^1_{B_0}\otimes R^{n-1}\check f_{0*}{\bf R})$
if and only if $h_n=\check h_n$.

Proof. Let $y_1,\ldots,y_n$ be action coordinates for $f$. (We
now have two different special Lagrangian fibrations, $f$ and $\check f$,
hence two different sets of action coordinates.)
Fixing $b\in B$, $X_b=\T^*_{B,b}/\Lambda$, $\Lambda$ is generated
by $e_1,\ldots,e_n$ with $e_i=dy_i$ and $\dual{\Lambda}\subseteq
\T_{B,b}$ is generated by $e_1^*,\ldots,e_n^*$, $e_i^*=\partial/\partial y_i$.
Now $\check X_b$ is identified with $\T_{B,b}^*/\dual{\Lambda}$, where
$e_i^*$ is identified with the 1-form $-h_n(e_i^*,\cdot)
=-\sum_j (h_n)_{ij} dy_j$, and thus $e_i\in\Lambda$ is identified
with $-\sum_j h_n^{ij}\partial/\partial y_j\in \T_{B,b}$.
Thus
$$\eqalign{\check h_n(\sum_j h_n^{ij} \partial/\partial y_j,
\partial/\partial y_k)&=
-\int_{\check X_b} \iota(\sum_j h_n^{ij}\partial/\partial y_j)\check
\omega\wedge \iota(\partial/\partial y_k)\im\check\Omega_n\cr
&=-\int_{e_i}\iota(\partial/\partial y_k)\im\check\Omega_n\cr}$$
by Proposition 4.6, where $e_i\in\Lambda=H_{n-1}(\check X_b,\boldz)
\cong H_1(X_b,\boldz)$. If $[\omega]=[\im\check\Omega_n]$,
this latter integral coincides with
$$-\int_{e_i}\iota(\partial/\partial y_k)\omega=\delta_{ik}.$$
Thus $\check h_n(\partial/\partial y_i,\partial/\partial y_j)
=(h_n)_{ij}$, so $\check h_n=h_n$. The argument reverses
to prove the converse. $\bullet$

\proclaim Moral 6.2. $\check\Omega_n$ must be chosen on $\check X$
so that $\check h_n=h_n$.

{\it Remark 6.3.} While this moral was deduced by beginning with a 
special Lagrangian torus fibration and applying the principal that
double dualising should bring one back to the initial fibration,
there is no reason this can't then be generalised to provide a guide
for putting complex structures on more general $D$-brane moduli
spaces. In the situation of \S 4, given a family $\U\rightarrow B$,
one has the metric $h_n$ on $B$. Then a holomorphic $n$-form should
be chosen on $\M$ so that $\M\rightarrow B$ is special Lagrangian
and the induced metric on $B$ is $h_n$.
\bigskip

It is more difficult to say exactly what role the $B$-field plays. According
to the conjecture originally stated in [15], and restated in
the introduction, $\check\Omega_n$ should be chosen so that
$[\check\Omega_n]-[\sigma_0]\in H^1(B,R^{n-1}\check f_*{\bf R})$
should coincide with the choice of the $B$-field $\b\in H^1(B,R^1f_*{\bf R})$.
This provides little guidance, but we will see an example of this 
below which may point in the correct direction for interpreting the
$B$-field. 

{\it Example 6.4.} (The Hitchin solution.) Hitchin [17] gives a choice
of $\check\Omega_n$ which under certain assumptions about the metric
$h_n$ satisfies $d\check\Omega_n=0$. He expresses it
locally in terms of action-angle coordinates, but it can be written
down in arbitrary coordinates on the base in a natural way.
One takes for $\check \nabla_b$ the Gauss-Manin connection. This is a linear
connection $\check\nabla_{GM}$ 
on $\T_{B_0}^*=(R^{n-1}\check f_{0*}{\bf R})\otimes
C^{\infty}(B_0)$, whose flat sections are the sections of $R^{n-1}f_{0*}
{\bf R}$. Taking the horizontal subspaces of this connection, it
is easy to see that these descend to give a flat Ehresmann connection on
$\check X_0$, which one takes to define $\check \nabla_b$.
Note that in action coordinates $v_1,\ldots,v_n$ for $\check f_0$,
the Gauss-Manin connection is trivial, so in these coordinates
one takes $b=0$. Since now $\check \nabla_b$ is flat, we must have
$\check g$ constant along fibres by Corollary 5.16. Thus
$$\eqalign{\check h_n(\partial/\partial y_i,\partial/\partial y_j)
&=\int_{\check X_b} \check V\check g^{ij} dx_1\wedge\cdots\wedge dx_n
/\int_{\check X_b} \check V dx_1\wedge\cdots\wedge dx_n\cr
&=\check g^{ij},\cr}$$
and since we want $\check h_n=h_n$, we have no choice but to take $\check g^{ij}
=(h_n)_{ij}$, giving rise to a choice of $\check\Omega_n$. 
To check to see if $d\check\Omega_n=0$, one checks conditions
(5.9)-(5.12). (5.9) and (5.11) are immediate, while (5.10) can be checked. (5.12)
is then equivalent, in this case, to $\int_{\check X_b}\check\Omega_n$
being independent of $b$. But
$$\eqalign{\int_{\check X_b}\check\Omega_n
&=\int_{\check X_b} \sqrt{\det \check g_{ij}}
dx_1\wedge\cdots\wedge dx_n\cr
&={1\over\sqrt{\det(\check h_n)}} \int_{\check X_b} dx_1\wedge\cdots\wedge
dx_n.\cr}$$ 
Hence this quantity must be independent of $b$ for (5.12) to be satisfied. 
In particular,
if $u_1,\ldots,u_n$ are action coordinates for $f_0$ and
$v_1,\ldots,v_n$ are action 
coordinates on the same open subset for $\check f_0$,
and $v_1,\ldots,v_n,x_1,\ldots,x_n$ are canonical coordinates,
then $\int_{\check X_b} dx_1\wedge\cdots\wedge dx_n$ is a constant
independent of $b$, so in these coordinates, $d\check\Omega_n=0$
is equivalent to
$$\det(\check h_n)=constant.$$
A simple calculation now shows that if this is the case, then
$\check\Omega_n$ coincides with Hitchin's $\tilde\Omega^c$ (up
to a constant factor) in [17], \S 6, where
$$\tilde\Omega^c=\bigwedge_i(dx_i+\sqrt{-1}du_i).$$
Thus one recovers [17], Proposition 5. Of course,
$\check\omega=\sum dx_i\wedge dv_i$, and Hitchin views mirror
symmetry as an exchanging of the roles of the two sets of action
coordinates $\{u_i\}$ and $\{v_i\}$. By [17], Proposition 3,
the condition $\det(\check h_n)=constant$ (in coordinates $v_1,\ldots,v_n$)
is equivalent to the condition $\det(h_n)=constant$ (in coordinates
$u_1,\ldots,u_n$). This holds in particular if the metric $g_{ij}$
is constant on fibres.

{\it Example 6.5.} (The Hitchin solution twisted by the $B$-field.) 
Continuing with the above example, assume $d\check\Omega_n=0$.
Now choose
a symmetric cohomology class $\b\in H^1(B,R^{n-1}\check f_*{\bf R})$,
with a symmetric representative $b\in \Gamma(B_0,\Omega^1_{B_0}\otimes
\Omega^1_{B_0})$; in action-angle coordinates for $\check f$,
this will be of the form $\sum b_{ij} dv_i\otimes dv_j$. Now take, in this
coordinate system, the $n$-form $\check\Omega_{n,b}$ to be given by the 
matrix $(\beta_{ij})=(b_{ij}+\sqrt{-1}(h_n)_{ij})$. This in fact gives
a well-defined $n$-form on all of
$\check X_0$. Note that since $(b_{ij})$ is symmetric, so is $\beta$,
and thus we just need to show that $d\check\Omega_{n,b}=0$ in order to
show that $\check\Omega_{n,b}$ induces a complex structure with a Ricci-flat
metric. This closedness can be seen to be true as follows. If
$U\subseteq B_0$ is a sufficiently small open set with action-angle
coordinates $v_1,\ldots,v_n,x_1,\ldots,x_n$ for $\check f$, 
we can find an element
$a\in\Gamma(U,\Omega^0_U\otimes\Omega^1_U)$ such that $\check\nabla_{GM}(a)=b$,
since $\check\nabla_{GM}(b)=0$. Here
$a=\sum_i a_i 1\otimes dv_i$, with $\partial a_i/\partial v_j=b_{ij}$.
Thus by symmetry of $b_{ij}$, $(v_1,\ldots,v_n)\mapsto (v_1,\ldots,v_n,a_1,
\ldots,a_n)$ gives a Lagrangian section $\sigma$ of $\check f^{-1}(U)\rightarrow
U$,
and it is easy to see that
$$\check\Omega_{n,b}=T_{\sigma}^*\check\Omega_n,$$
where $\check\Omega_n$ is the Hitchin solution of Example 6.4. Since
$d\check\Omega_n=0$, $d\check\Omega_{n,b}=0$ also.

Note also that if $b'=\sum b_{ij}'dy_i\otimes dy_j$ is a different symmetric
representative for $\b$, then
there exists an $a\in\Gamma(B_0,\Omega^0_{B_0}\otimes \Omega^1_{B_0})$
such that $\nabla_{GM}(a)=b'-b$. As before, $a$ gives a Lagrangian 
section $\sigma$ of $\check f_0$, and
$$T_{\sigma}^*\check\Omega_{n,b}=\check\Omega_{n,b'}.$$
Finally, one sees that $\check\Omega_{n,b}-\check\Omega_n$ represents
the class $\b\in H^1(B,R^{n-1}\check f_{*}{\bf R})$, and  $\check
\Omega_{n,b}$ and $\check\Omega_n$ yield the same complex structure
if $\b\in H^1(B, R^{n-1}\check f_{*}\boldz)$, for then there is a global
section $\sigma$ of $f:X\rightarrow B$
with $T_{\sigma}^*\check\Omega_n=\check\Omega_{n,b}$.

One can in fact go further. We have observed that for small open sets
$U\subseteq B_0$, the complex and K\"ahler structures on $\check X_0$
induced by $\check\Omega_n$ and $\check\Omega_{n,b}$ on $\check f^{-1}(U)$
are isomorphic, so one should think of the K\"ahler structure induced
by $\check\Omega_{n,b}$ as a torseur over that induced by $\check\Omega_n$.
Specifically, fixing the complex structure $\check\Omega_n$ on
$\check X_0$, note that the sheaf ${\cal A}$ on $B_0$ defined by
$${\cal A}(U)=\{\hbox{sections $\sigma:U\rightarrow \check f^{-1}(U)$ such
that $T_{\sigma}^*\check\omega=\check\omega$ and $T_{\sigma}^*\check\Omega_n
=\check\Omega_n$}\}$$ 
coincides with $R^{n-1}f_*{\bf R}/\boldz$. In fact, writing these
conditions in the coordinates $u_i$ and $x_i$ of Example 6.4, one
sees that the condition $T_{\sigma}^*\check\Omega_n=\check\Omega_n$
guarantees that the section $\sigma$ is constant with respect
to the Gauss-Manin connection. Thus the set
of all special Lagrangian fibrations over $B_0$ obtained from
$\check f:\check X_0\rightarrow B_0$ by regluing using these
translations is $H^1(B_0,R^{n-1}\check f_{0*}{\bf R}/\boldz)$.
Because $\check\omega$ and $\check\Omega_n$ are preserved by these translations,
they glue to give forms on the twisted fibrations. Thus each
element $\b\in H^1(B_0, R^{n-1}\check f_{0*}{\bf R}/\boldz)$ gives
rise to a fibration $\check f_{\b}:\check X_{0,\b}\rightarrow B_0$
with symplectic form $\check\omega_{\b}$ and holomorphic $n$-form 
$\check\Omega_{n,\b}$. This is a potentially wider class of examples than
were constructed above using symmetric cohomology classes in $H^1(B,
R^{n-1}\check f_*{\bf R})$; if $\b$ does not come from a symmetric class, then 
$\check f_{\b}$ will not possess a Lagrangian section, and may not even
possess a topological section.

Note also that if $f$ and $\check f$ are ${\bf R}/\boldz$-simple, then
$H^1(B,R^{n-1}\check f_{*}{\bf R}/\boldz)
\cong H^1(B,R^1f_{*}{\bf R}/\boldz)$. 
This leads us to conjecture that
the correct group for the $B$-field to live in is $H^1(B,R^1f_*
{\bf R}/\boldz)$. This new 
proposed definition for the $B$-field is dependent not just on $X$
but on the topology of the fibration, and
even in the threefold case does not necessarily coincide
with $H^2(X,{\bf R}/\boldz)$, as we saw in Example 3.11. 
Nevertheless, I believe this is
the correct interpretation of the $B$-field.
\bigskip

This now leads us to a refined mirror symmetry conjecture.

\proclaim Conjecture 6.6. Let $X$ be a compact Calabi-Yau $n$-fold with
$\omega$, $\Omega$ the symplectic form and holomorphic $n$-form respectively,
and suppose $f:X\rightarrow B$ is an integral special Lagrangian
fibration. Then for each $\b\in H^1(B,R^1f_*{\bf R}/\boldz)$, there
is a Calabi-Yau $n$-fold $\check X$ with symplectic form $\check\omega$
and holomorphic normalized $n$-form $\check \Omega_n$ along with
an integral special Lagrangian fibrations $\check f:\check X\rightarrow B$
such that
\item{(1)} for each open set $U\subseteq B_0=B-\Delta$ on which
both $f$ and $\check f$ have sections, $f^{-1}(U)\rightarrow U$
and $\check f^{-1}(U)\rightarrow U$ are topologically dual fibrations.
\item{(2)} For $b\in B_0$,  $\gamma\in H_1(\check X_b,\boldz)\cong
H_{n-1}(X_b,\boldz)$ and $v\in\T_{B,b}$,
$$\int_{\gamma} \iota(v)\check\omega=\int_{\gamma} \iota(v)\im\Omega_n.$$
\item{(3)} For $\gamma\in H_{n-1}(\check X_b,\boldz)\cong H_1(X_b,\boldz)$,
$$\int_{\gamma}\iota(v)\omega=\int_{\gamma}\iota(v)\im\check\Omega_n.$$
\item{(4)} $\check f$ possesses a topological section if
$$\b\in H^1(B,R^1f_*{\bf R})/H^1(B,R^1f_*\boldz)\subseteq
H^1(B,R^1f_*{\bf R}/\boldz).$$
In this case, if $\check\sigma$ is a topological section,
then $[\re\check\Omega_n]-[\check\sigma]$ defines a class in
$H^1(B,R^{n-1}\check f_*{\bf R})$ which is well-defined modulo
$H^1(B,R^{n-1}\check f_*\boldz)$, and agrees with $\b$ in
$$H^1(B,R^1f_*{\bf R})/H^1(B,R^1f_*\boldz).$$
\item{(5)} If $\check J$ is the Jacobian of $\check X$, then
$\check X$ is obtained from $\check J$ as a symplectic manifold
via the image of $\b$ under the composed map $H^1(B,
R^1f_*{\bf R}/\boldz)\cong H^1(B,R^{n-1}\check f_*{\bf R}/\boldz)
\rightarrow H^1(B,\Lambda(\check J^{\#}))$ of Remark 3.15. 
\item{(6)} Once $\check X$ and $\check\omega$ are fixed,
$\check\Omega_n$ is unique up to translation by
a Lagrangian section of $\check J$ acting on $\check X$.

I do not however suggest that $\check f:\check X\rightarrow B$ is
obtained as a {\it K\"ahler} manifold as a torseur over some
basic $\check J\rightarrow B$. While this occurred in Example 6.5,
there is no reason to suspect this works when the metric on the fibres
is not flat. We simply don't expect there to be 
isometries given by translation by a section. However, in some sense
this might provide an initial approximation to the correct answer.

{\it Remark 6.7.} We can show a local form of the conjectured 
uniqueness. Suppose $\Omega_t$ is a family of holomorphic $n$-forms
on $X$ with respect to which a fixed symplectic form
$\omega$ induces a Ricci-flat metric
and $f:X\rightarrow B$ is special Lagrangian. In addition assume
$[\Omega_t]\in H^n(X,{\bf C})$ is a fixed cohomology class.
Then by local Torelli
all $\Omega_t$ induce the same complex structure, so there exists
diffeomorphisms $\phi_t:X\rightarrow X$ such that $\phi_t^*\Omega_t
=\Omega_0$, $\phi_0=id$. Now $\phi_t^*\omega$ is a symplectic form
on $X$ inducing a Ricci-flat metric in the complex structure induced by
$\phi_t^*\Omega_t=\Omega_0$, and represents the same cohomology class
as $\omega$, so by uniqueness of Ricci-flat metrics, $\phi_t^*\omega
=\omega$. Thus $\phi_t$ is a family of symplectomorphisms. Assuming
$H^1(X,{\bf R})=0$, differentiating this family of diffeomorphisms at
$t=0$ yields a Hamiltonian vector field $v$ induced by a Hamiltonian 
function $H$ on $X$: $\iota(v)\omega=dH$.
Then $\im\Omega_t|_{X_b}=0$ for all $b$ implies that
$(\L_{v}\im\Omega_0)|_{X_b}=0$ for all $b$. But 
$\L_v\im\Omega_0=d(\iota(v)\im\Omega_0)$, and if $\Omega_0$ is
given as usual by a matrix $\beta=(\beta_{ij})$, $\beta_{ij}=b_{ij}+ig^{ij}$,
then a simple calculation shows that
$$(d(\iota(v)\im\Omega_0))|_{X_b}=-\sum_{i,j} Vg^{ij} {\partial^2 H
\over \partial x_i\partial x_j}dx_1\wedge\cdots\wedge dx_n.$$
Here we have used the fact that $\iota(\partial/\partial y_i)\im\Omega_0$
is closed. Thus we see that on each fibre, $H$ satisfies the second
order elliptic partial differential equation
$$\sum_{i,j} g^{ij}{\partial^2 H\over\partial x_i\partial x_j}=0.$$
By the maximum principal, $H$ cannot have a local maximum on each non-singular
fibre unless
$H$ is constant on the fibre. Since the set of non-singular fibres is dense,
we conclude that $H$ is the pullback of a function on $B$.
In particular, it follows from \S 2 that
$\phi_t$ must be translation by a Lagrangian section.

{\it Remark 6.8.} One natural question is to determine
the relationship between $Vol(X_b)=\int_{X_b} \Omega$ and
$Vol(\check X_b)=\int_{\check X_b}\check\Omega$, since
knowledge of the latter allows us to reconstruct $\check \Omega$
from $\check\Omega_n$. We can describe this relationship if
the metric is constant
along the fibres of $f:X\rightarrow B$ and $\check f:\check X\rightarrow
B$. Let $y_1,\ldots,y_n$ be action coordinates for $f$. As in 
Example 6.4, 
$$h_n^{ij}=g_{ij}.$$
Then
$$\eqalign{Vol(X_b)&=\int_{X_b}\sqrt{\det(g_{ij})} dx_1\wedge\cdots\wedge
dx_n\cr
&=\sqrt{\det(h_n^{ij}) }.\cr}$$
On the other hand, for a fixed $b\in B$, $\check X_b$ can be written
canonically as $\T_{B,b}^*/\dual{\Lambda}$ as in Proposition 4.7,
where $\dual{\Lambda}$ is generated by the one-forms $h_n(\partial/\partial y_i,
\cdot)$, i.e. the one-forms $\sum_j (h_n)_{ij} dy_j$. The  metric on the fibre
$\check X_b$ is still given by $\check g_{ij}=(h_n)^{ij}$, since $h_n
=\check h_n$. Thus
$$\eqalign{Vol(\check X_b)&=\int_{\check X_b} \sqrt{\det (h_n^{ij})}
dx_1\wedge\cdots\wedge dx_n\cr
&=\sqrt{\det (h_n^{ij})}\det ((h_n)_{ij})\cr
&={1\over Vol(X_b)}.\cr}$$
This is the familiar ``$R\mapsto 1/R$'' relationship of $T$-duality.
If the metric is not flat, we expect some corrections to the volume,
and this may affect this relationship. However, as we shall see in \S 7,
this relationship continues to hold for K3 surfaces.
\bigskip

We end this section with a brief discussion of the Yukawa coupling. Mirror
symmetry instructs us that given a Calabi-Yau $X$, the Yukawa
coupling on $\check X$ contains information about the genus 0 Gromov-Witten
invariants on $X$, and in particular in the three dimensional
case, these Gromov-Witten invariants can be completely 
recovered from the Yukawa coupling on $\check X$. How do we see the
Yukawa coupling in the context of special Lagrangian fibrations?

Suppose that Conjecture 6.6 holds. Assume for simplicity that we only
consider values of the $B$-field $\b\in H^1(B,R^{n-1}\check f_*
{\bf R}/\boldz)$ which come from symmetric classes in $H^1(B,R^{n-1}\check
f_*{\bf R})$, so that although $\b$ varies, we can fix the underlying manifold
$\check X$ and symplectic form $\check \omega$ and
simply let $\check\Omega_n$ vary; we denote the dependence on $\b$
by writing $\check\Omega_{n,\b}$. Of course, we will not have
$\check\Omega_{n,\b+\alpha}=\check\Omega_{n,\b}$ for
$\alpha\in H^1(B,R^{n-1}\check f_*\boldz)$, but merely expect that there
exists a Lagrangian section $\sigma$ of $\check f:\check X\rightarrow
B$ with $T_{\sigma}^*\check\Omega_{n,\b}
=\check\Omega_{n,\b+\alpha}$.
The $(1,n-1)$-Yukawa
coupling of interest is then, for tangent directions $\partial/
\partial b_1,\ldots,\partial/\partial b_n\in H^1(B,R^{n-1}\check f_*
{\bf R})$, the tangent space of the
torus $H^1(B,R^{n-1}\check f_*{\bf R}/\boldz)$,
$$\left\langle {\partial\over\partial b_1},\ldots,{\partial\over\partial
b_n}\right\rangle=
\int_{\check X} \check\Omega_{n,\b}\wedge {\partial^n\over\partial b_1\ldots
\partial b_n}
\check\Omega_{n,\b}.$$
In local coordinates, we write
$$\check\Omega_{n,\b}=V_{n,\b}\theta_1(\b)\wedge\cdots\wedge\theta_n(\b).$$
Note that in taking the $n$ derivatives of $\check\Omega_{n,\b}$
by using the product rule, all terms still containing any undifferentiated
$\theta_i$ will disappear after we wedge with $\check\Omega_{n,\b}$.
Thus
$$\int_{\check X} \check\Omega_{n,\b}\wedge {\partial^n\over\partial b_1
\ldots\partial b_n}
\check\Omega_{n,\b}=\sum_{\sigma\in S_n}
\int_{\check X} V_{n,\b}^2\theta_1\wedge\dots\wedge\theta_n
\wedge {\partial \theta_1\over\partial b_{\sigma(1)}
}\wedge\cdots\wedge{\partial \theta_n
\over\partial b_{\sigma(n)}}.$$
Writing $\theta_i=dx_i+\sum_{j} \check\beta_{ij}(\b)dy_j$, 
$${\partial\theta_i\over\partial b_{\sigma(i)}}=
\sum{\partial\check\beta_{ij}\over
\partial b_{\sigma(i)}}dy_j,$$
so the above integral is
$$\int_{\check X} V_{n,\b}^2\sum_{\sigma\in S_n}
\det\left({\partial\check\beta_{ij}\over\partial b_{\sigma(i)}}
\right) dx_1\wedge\cdots\wedge dx_n\wedge dy_1\wedge\cdots\wedge
dy_n.$$
So far, this is not particularly illuminating in the general case.
However, for the twisted Hitchin solutions, it is
an elementary calculation to show this integral can be evaluated
in terms of the topological coupling on $X$, as expected. This comes from
observing that if $\partial\theta_j/\partial b_i=\sum_k b_{jk}^idy_k$,
then in action-angle coordinates, where $V_n=1$, the integrand
above coincides with
$$(-1)^{n(n-1)/2}\bigwedge_{i=1}^n\left(\sum_{j,k} b_{jk}^i dx_j\wedge dy_k
\right).$$
\bigskip

{\hd \S 7. K3 Surfaces.}

Mirror symmetry for K3 surfaces has been completely understood using
Torelli theorems for K3 surfaces. We will now show that the
previous material of this paper gives us a differential geometric
construction of mirror symmetry for K3 surfaces, and in doing so, we will
show Conjecture 6.6 holds in two dimensions. In other words, given a special
Lagrangian $T^2$-fibration on a K3 surface, and a choice of $B$-field,
we will construct the mirror in the sense made explicit in Conjecture
6.6. This will prove to be a variant of the mirror symmetry
for K3 surfaces described in [6] and [11].

To begin, let $S$ be a K3 surface with holomorphic 2-form $\Omega$
and K\"ahler form $\omega$ corresponding to a Ricci-flat metric. We
insist on the usual normalisation and this implies in particular
that $(\re\Omega)^2=(\im\Omega)^2=\omega^2>0$ and $(\re\Omega)\wedge
(\im\Omega)=\omega\wedge(\re\Omega)=\omega\wedge(\im\Omega)=0$. In order
for $S$ to possess a special Lagrangian fibration there must be
a cohomology class $E\in H^2(S,\boldz)$ such that $E^2=0$ and
$\omega.E=0$. We assume such a class exists, and we fix it. 
We take $E$ to be primitive. Now one 
constructs a special Lagrangian fibration on $S$ by the usual
hyperk\"ahler trick, as originally suggested in [25]. First, multiply
$\Omega$ by a phase $e^{i\theta}$ to ensure $\im\Omega.E=0$.
Following the notation of [15], there is a complex structure $K$
compatible with the Ricci-flat metric in which
$\Omega_K=\im\Omega+i\omega$ and $\omega_K=\re\Omega$. Then special Lagrangian
submanifolds on $S$ are complex submanifolds in the $K$ complex structure.
Denote the K3 surface in the $K$ complex structure by $S_K$. Then
by construction $E.\Omega_K=0$ so $E\in Pic(S_K)$. It is then
standard that $S_K$ possesses an elliptic fibration. However, the class
of the fibre need not be $E$; $E$ might be represented by a fibre
plus a sum of $-2$ divisors. In any event, replace $E$ by the
class of the fibre, positively oriented. $E$ will now remain fixed,
and we have an elliptic fibration $f:S_K\rightarrow\Pone$, which we
identify with a special Lagrangian $T^2$ fibration $f:S\rightarrow B=S^2$.
We will assume $f$ is integral, which is true if and only if there
is no $\delta\in Pic(S_K)$ with $\delta^2=-2$ and $\delta.E=0$. This will
certainly hold for general $S$.

We now consider the spectral sequence for $f$ over $\boldz$. Since there
exists a class $\sigma$ such that $\sigma.E=1$, the sequence
in fact degenerates and takes the form
$$\matrix{H^0(B,\boldz)&0&H^2(B,\boldz)\cr
0&H^1(B,R^1f_*\boldz)&0\cr
H^0(B,\boldz)&0&H^2(B,\boldz)\cr}$$
It is also clear that $H^1(B,R^1f_*\boldz)$ is canonically isomorphic
to $E^{\perp}/E$.

We now wish to construct a mirror fibration given the data
$$\vbox{\narrower\narrower
$f:S\rightarrow B$ as above and a choice of $B$-field $\b\in
E^{\perp}/E\otimes{\bf R}/\boldz\cong H^1(B,R^1f_*{\bf R}/\boldz)$.

}$$

We first recall some facts about elliptic fibrations. See [7]
for general facts about the analytic theory of elliptic surfaces.

The fibration $f:S_K\rightarrow\Pone$ in general does not possess a
holomorphic section; in fact for general choice of $S$, $\Pic S_K
=\boldz E$. However, there is a Jacobian fibration $j:J_K\rightarrow\Pone$
of $f$ which is locally isomorphic to $f$, and which does possess a 
holomorphic section.

\proclaim Proposition 7.1. There is a diffeomorphism $\phi:J_K\rightarrow S_K$
over $\Pone$ which is holomorphic when restricted to each fibre.
Furthermore, if $U\subseteq\Pone$ is an open subset on which there exists
a biholomorphic map $\xi:j^{-1}(U)\rightarrow f^{-1}(U)$ over
$U$, then $\phi^{-1}\circ\xi:j^{-1}(U)\rightarrow j^{-1}(U)$ is
given by translation by a (not necessarily holomorphic) section of 
$j^{-1}(U)$.

Proof. This follows easily from the fact that the
``$C^{\infty}$ Tate-Shafarevich group,'' i.e. the first cohomology
group of the sheaf of $C^{\infty}$ sections of $j:J_K^{\#}\rightarrow
\Pone$, is zero. Thus $f:S_K\rightarrow\Pone$ possesses a $C^{\infty}$
section, and $\phi$ can be taken to identify this $C^{\infty}$
section of $f$ with a holomorphic section of $j$, such
that $\phi$ is holomorphic on each fibre.
$\bullet$

We fix one holomorphic section $\sigma_0$ of $j:J\rightarrow\Pone$,
and identify $\sigma_0$ with the topological section $\phi(\sigma_0)$ of
$f:S\rightarrow B$.
Having chosen this section, we can take it to be the zero section
of $j:J_K\rightarrow\Pone$ and obtain a standard exact sequence
$$0\rightarrow R^1f_*\boldz\rightarrow R^1f_*\O_{S_K} \mapright{\psi}
J_K^{\#}\rightarrow 0$$
where $J_K^{\#}$ denotes the sheaf of holomorphic sections of $j:J_K\rightarrow
\Pone$.
Here $R^1f_*\O_{S_K}$ can be
identified with the normal bundle of the zero section, and the map
$\psi$ is just the fibre-wise exponential map. For K3 surfaces, $R^1f_*\O_{S_K}
\cong\omega_{\Pone}$. The underlying real bundle is $\T_{S^2}^*$.
This also gives a map $\pi:\T_{S^2}^*\rightarrow S_K^{\#}$ with $\pi=\phi\circ
\psi$.

Just as in the real case, the total space of $\omega_{\Pone}$, the holomorphic
cotangent bundle of $\Pone$, has a canonical holomorphic symplectic
form $\Omega_c$. In local coordinates, if $z$ is a coordinate on $\Pone$
and $w$ the canonical coordinate on the cotangent bundle, then 
$\Omega_c=dw\wedge dz$. Furthermore, any holomorphic symplectic form on the
cotangent bundle of $\Pone$ is proportional to $\Omega_c$. 

\proclaim Proposition 7.2. There is a map $\chi:\T_{S^2}^*\rightarrow\T_{S^2}^*$
given by fibrewise
multiplication by a complex constant so that, for $\pi'=\pi\circ\chi$,
$\pi'^*(\Omega_K)=\Omega_c+f^*\alpha$, where $\alpha$ is a 2-form on $S^2$.

Proof.
First note that there are two different complex structures on the total
space of $\T^*_{S^2}$ in this picture: one is the standard complex
structure coming from being the holomorphic line bundle $\O_{\Pone}(-2)$,
while the other is induced by $\pi^*\Omega_K$. To distinguish
between these two complex structures, let $\tilde J$ be the
total space of $\T^*_{S^2}$ with the standard complex structure,
and let $\tilde S$ denote the total space of $\T_{S^2}^*$ with the complex
structure induced by $\pi^*\Omega_K$. Let $\tilde j:\tilde J\rightarrow
S^2$, $\tilde f:\tilde S\rightarrow S^2$ be the projections.

Now let $U\subseteq S^2$ be a sufficiently small open set so that there 
exists a biholomorphic map $\xi:j^{-1}(U)\rightarrow f^{-1}(U)$ over $U$.
By Proposition 7.1, $\phi^{-1}\circ\xi=T_{\sigma}$ for some section
$\sigma$ of $j^{-1}(U)\rightarrow U$, and if $U$ is small enough,
$\sigma$ can be lifted to a section $\tilde\sigma$ of $\tilde j^{-1}(U)
\rightarrow U$. We let $\tilde\xi$ denote translation by the
section $\tilde\sigma$ (with the zero-section of $\T_{S^2}^*$ taken
to be the origin). We then have a commutative diagram
$$\matrix{\tilde j^{-1}(U)&\mapright{\tilde\xi}&\tilde f^{-1}(U)\cr
\mapdown{\psi}&&\mapdown{\pi}\cr
j^{-1}(U)&\mapright{\xi}&f^{-1}(U)\cr}$$
since $\pi\circ\tilde\xi=\phi\circ\psi\circ T_{\tilde\sigma}
=\phi\circ T_{\sigma}\circ \psi=\xi\circ\psi$.
In particular, since $\Omega_K$ is a holomorphic 2-form on $\tilde 
f^{-1}(U)$, $\tilde\xi^*\pi^*\Omega_K$ is a holomorphic 2-form on
$\tilde j^{-1}(U)$, and hence can be written as $g dw\wedge dz$ in local
coordinates, for some holomorphic function $g$. This function $g$ must
be constant along the fibres of $\tilde j$ since $g$ must
descend to a holomorphic function on the compact fibres of $j$. Thus on
$\tilde f^{-1}(U)$,
$$\eqalign{\pi^*\Omega_K&=(\tilde\xi^{-1})^*(g dw\wedge dz)\cr
&=T_{-\tilde\sigma}^*(gdw \wedge dz)\cr
&=(g\circ T_{-\tilde\sigma}) dw\wedge dz+h dz\wedge d\bar z\cr}$$
for some function $h$. Of course, $g\circ T_{-\tilde\sigma}=g$
since $g$ is constant on fibres.

Let $i:\tilde J\rightarrow \tilde S$ be the identity map; this is 
of course non-holomorphic. The above equation shows that the $(2,0)$
part of $i^*\pi^*\Omega_K$, being locally of the form $gdz\wedge dw$,
is in fact a holomorphic 2-form. In addition, this holomorphic
2-form is nowhere vanishing: if $g$ vanishes then $\Omega_K\wedge
\bar\Omega_K=0$ at that point. Thus the $(2,0)$ part of $i^*\pi^*\Omega_K$
is proportional to $\Omega_c$, say $C\Omega_c$. In addition, we then
see from $d\Omega_K=0$ that $h$ must be constant along fibres and
hence $\pi^*\Omega_K-C\Omega_c$ is the pullback of a $(1,1)$ form $\alpha$
on $\Pone$, i.e.
$$\pi^*\Omega_K=C\Omega_c+\tilde f^*\alpha.$$
Now let $\chi:\T^*_{S^2}\rightarrow \T^*_{S^2}$ be given by 
$w\mapsto C^{-1}w$. Then $\chi^*\pi^*\Omega_K=\Omega_c+
\tilde f^*\alpha$ as desired. $\bullet$

To sum up, replacing $\pi$ by $\pi'$, we now have a map
$\pi:\T_{S^2}^*\rightarrow S_K^{\#}$ with kernel $R^1f_*\boldz$ and such
that 
$\pi^*\Omega_K=\Omega_c+\tilde f^*\alpha$.
This gives, identifying the underlying topological spaces $S_K^{\#}$
and $S^{\#}$, $\pi^*\omega=\im(\Omega_c+\tilde f^*\alpha)$
and $\pi^*(\im\Omega)=\re(\Omega_c+\tilde f^*\alpha)$. 

A local description of these forms are as follows: given complex
canonical coordinates $z,w$ on $\T_U^*$, $z$ a coordinate on $U$,
write $z=y_1-iy_2$, $w=x_1+ix_2$.
The signs are chosen so that with real coordinates $y_1,y_2$ on the
base, $y_1,y_2,x_1,x_2$ are canonical coordinates on $\T_{S^2}^*$. Then
$$\eqalign{\pi^*\omega&=\im(dw\wedge dz+\tilde f^*\alpha)\cr
&=dy_2\wedge dx_1+dx_2\wedge dy_1+\tilde f^*\im\alpha\cr}\leqno{(7.1)}$$
and
$$\eqalign{
\pi^*\im\Omega&=\re(dw\wedge dz+\tilde f^*\alpha)\cr
&=
dx_1\wedge dy_1+dx_2\wedge dy_2+\tilde f^*\re\alpha.\cr}
\leqno{(7.2)}$$
This completes our first goal of finding coordinates on $S^{\#}$ in which
$\omega$ and $\im\Omega$ have simple forms. Note that our map
$\pi:\T^*_{S^2}\rightarrow S^{\#}$ is not the same as defined in \S 2
because the symplectic form is not the standard one. We will see why
we have made this choice in the proof of the following theorem.

\proclaim Theorem 7.3. Conjecture 6.6 holds for the general integral
special Lagrangian fibration $f:S\rightarrow B=S^2$.

Proof. Choose a $B$-field $\b\in (E^{\perp}/E)\otimes{\bf R}/\boldz$.
Lift this to a representative $\b\in E^{\perp}/E\otimes {\bf R}$. At times,
we will also further lift $\b$ to an element $\b\in E^{\perp}\otimes
{\bf R}$ chosen so that $\b.[\sigma_0]=0$, where $\sigma_0$
is the fixed topological section of $f:S\rightarrow B$
chosen previously.
We are now trying to construct a special Lagrangian fibration
$\check f:\check S\rightarrow B$ satisfying the properties of Conjecture
6.6. Because this fibration may not have a Lagrangian section,
we first construct the Jacobian 
$\check j:\check J\rightarrow B$ (in the sense of \S 2) as a symplectic
manifold. This Jacobian
should be the dual fibration with a symplectic form $\check\omega_{\check J}$
as constructed in \S 4.

To do so, we reembed
$R^1f_*\boldz\hookrightarrow \T_{S^2}^*$ using the periods given
by $\im\Omega_n$. Now $\Omega_n={1\over Vol(S_b)}\Omega$, and this
embedding takes a cycle $\gamma\in H_1(S_b,\boldz)\cong H^1(S_b,\boldz)
\subseteq \T_{S^2}^*$ to the one-form
$$\eqalign{v&\mapsto -\int_{\gamma} \iota(v)\im\Omega_n\cr
&=-{1\over Vol(S_b)}\int_{\gamma} \iota(v) (dx_1\wedge dy_1+dx_2\wedge dy_2)
\cr}$$
which yields the 1-form $\gamma/Vol(S_b)$.

The moral is:
{\it The dual lattice is simply the original lattice scaled by a factor of
$1/Vol(S_b)$.}

We have in fact chosen the map $\pi$ so that this would happen and so
make it transparent that dualising does not change the topology of
the fibration.

Instead of rescaling the lattice, it is easier to identify $\check J$
with $S$ topologically, and rescale the symplectic form. Since
in local coordinates we want $\check\omega_{\check J}
=dx_1\wedge dy_1+dx_2\wedge dy_2$, by (7.2) we take
on $\check J=S$
$$\check\omega_{\check J}=(\im\Omega-f^*\re\alpha)/Vol(S_b).$$
How do we obtain $\check S$ as a symplectic manifold?
Item (5) of Conjecture 6.6 instructs us to proceed as follows.
$\b\in H^1(B,R^1f_*{\bf R}/\boldz)\cong H^1(B,R^1\check j_*{\bf R}/\boldz)$
maps to an element of $H^1(B,\Lambda(\check J^{\#}))$. Having
lifted $\b$ to an element of $H^1(B,R^1f_*{\bf R})$, we obtain
in this way an element $-\b\wedge[\check\omega_{\check J}]\in H^2(B,{\bf R})$, 
via Remark 
3.15, which then maps to the appropriate element of $H^1(B,\Lambda(\check
J^{\#}))$. As remarked before Example 2.7, this means
$\check\omega=\check\omega_{\check J}+\check j^*\alpha_1$, where
$\alpha_1$ is a form on $B$ such that
$$\int_{B}\alpha_1=-\b.[\check\omega_{\check J}].$$
Now $[\check\omega_{\check J}]=[\im\Omega_n]$ by construction, so it
would appear that we take $\int_B\alpha_1=-\b.[\im\Omega_n]$. There is 
a slight subtlety in this: here we are representing $\b\in H^1(B,R^1f_*{\bf R})$
as an element of $E^{\perp}/E\otimes{\bf R}$, but that does not mean that
if we reinterpret $\b$ as a class in $H^1(B,R^1\check j_*{\bf R})$,
this class will coincide with the original $\b$ in $S$ under our identification
of $S$ with $\check J$. In fact, the correct choice is
$$\int_{B}\alpha_1=\b.[\im\Omega_n].$$
We will be vague about this here, but will see this sign change more
explicitly shortly.
We then set
$$\check\omega=\check\omega_{\check J}+\check j^*\alpha_1.$$
Note that the choice of $\alpha_1$ is not important; any
two choices representing the same cohomology class
can be identified by translation by a section.
Let $\check S$ be the symplectic manifold obtained with 
underlying manifold $S$ and symplectic form $\check\omega$, 
and let $\check f:\check S\rightarrow B$
be the same map as $f:S\rightarrow B$. We note that the cohomology class
of $\check\omega$ satisfies the relation
$$[\check\omega]=[\im\Omega_n]+(\im\Omega_n.(\b-[\sigma_0]))E.$$

Next we construct the form $\im\check
\Omega_n$. The first observation is that $\iota(v)\im\check\Omega_n$
must be harmonic for any $v\in\T_{B,b}$. On the other
hand, the same is true of $\iota(v)\check\omega$,
and $\iota(v)\check\omega=adx_1+bdx_2$ for $a$ and $b$ constant.
Thus we already know, in the 2-dimensional case, the harmonic
$n-1$ forms. This is a crucial point in dimension 2 which fails
in higher dimensions. Applying item (3) of Conjecture 6.6, in the form
given in Proposition 6.1, we can now determine $\check\Omega_n$ as follows.
First
$$h_n(\partial/\partial y_i,\partial/\partial y_j)={\delta_{ij}\over Vol(S_b)}
\int_{S_b} dx_1\wedge dx_2$$
at a point $b\in B$, as is easily computed from the definition and
(7.1), (7.2). Then
in order for $\check h_n=h_n$, we must have
$$-\int_{S_b}\iota(\partial/\partial y_i)\check\omega
\wedge\iota(\partial/\partial y_j)\im\check\Omega_n=
{\delta_{ij}\over Vol(S_b)}\int_{S_b} dx_1\wedge dx_2,$$
and a quick calculation shows this implies 
we locally can write
$$\im\check\Omega_n=-\omega+hdy_1\wedge dy_2,$$
where $h$ is a function. But the condition that $d\im\check\Omega_n=0$
implies $h$ is constant on fibres, so 
$$\im\check\Omega_n=-\omega+f^*\alpha_2$$
for some form $\alpha_2$ on the base. 

Here we see the sign reversal explicitly.
It might be a bit surprising that we have obtained $-\omega$ instead
of $\omega$. But this is the fault of the identification we have chosen.
We want $[\im\check\Omega_n]=[\omega]$ as classes in 
$H^1(B,R^1\check f_*{\bf R})\cong H^1(B,R^1f_*{\bf R})$; it is only
an accident of dimension that we have been able to identify $S$ and
$\check S$ as manifolds and then compare cohomology classes directly.
In fact, this sign change {\it must} occur if we want to identify $S$ and
$\check S$ without changing the orientation of the fibres.

Condition (3) of Conjecture 6.6
does not tell us what $\alpha_2$ must be; $\alpha_2$ is not
determined until one knows something about $\re\check\Omega_n$.
In fact, condition (4) tells us that we require
$$[\re\check\Omega_n]=[\sigma_0]-\b\mod E,$$
where again we are making use of the sign reversal observed
above.
Now our form $\im\check\Omega_n$ constructed above satisfies
$$[\im\check\Omega_n]=[-\omega]\mod E,$$
which tells us that in order for $\check\Omega_n^2=0$, we must have
$[\check\Omega_n]$ satisfying 
$$[\check\Omega_n]=[\sigma_0]-(\b+i\omega)+(1-(\b+i\omega)^2/2
+i(\omega.\sigma_0))E.\leqno{(7.3)}$$
Here we have chosen a representative of $\b\in E^{\perp}$ such that
$\b.\sigma_0=0$. Thus, in particular, we need
to choose $\alpha_2$ so that
$$\int_B\alpha_2=-\b.\omega+\sigma_0.\omega.$$
Again, we need to ask how much freedom we have to choose $\alpha_2$,
given that we have fixed $\alpha_1$. We had seen that $\alpha_1$
could be chosen to be any representative of its
cohomology class, as any choice could be obtained from any other
by translating $\check\omega$ by a section of $\T_{S^2}^*$. Once we have fixed
the form $\alpha_1$, however, we can only translate by
sections corresponding to 1-forms $\sigma$ with $d\sigma=0$.
Let $\sigma=\sigma_1dy_1+\sigma_2dy_2$. Then
$$T_{\sigma}^*(dy_2\wedge dx_1+dx_2\wedge dy_1)=
dy_2\wedge dx_1+dx_2\wedge dy_1-\left({\partial \sigma_1\over\partial y_1}
+{\partial\sigma_2\over\partial y_2}\right)dy_1\wedge dy_2.$$
Thus $T_{\sigma}^*(\omega)-\omega=-f^*(d*\sigma)$, where $*$ denotes
the Hodge $*$-operator in, say, the Fubini-Study metric on $B=\Pone$.
Thus if $\alpha_2,\alpha_2'$ are two 2-forms on $B$ representing
the same cohomology class, we just need to find a 
1-form $\sigma$ on $B$ such that $d\sigma=0$ and $d*\sigma=\alpha_2-\alpha_2'$,
and then $T_{\sigma}^*(\omega+f^*\alpha_2)=\omega+f^*\alpha_2'$. By the
Hodge theorem, such a $\sigma$ can always be found. Thus we have
complete freedom to choose $\alpha_2$, and any two choices
are related by translation by a Lagrangian section.

This will be the last remaining choice in the construction which is not forced
on us by any items in Conjecture 6.6; thus the uniqueness of item (6)
of Conjecture 6.6 will hold, given that the lack of uniqueness in the lifting
of $\b$ can be rectified by changing the choice of the zero section $\sigma_0$.

Finally, we set 
$$\im\check\Omega=(\im\check\Omega_n)/Vol(S_b).$$
It now follows immediately from (7.1) and (7.2) that as forms,
$$(\im\check\Omega)\wedge (\im\check\Omega)=\check\omega\wedge\check\omega
>0,$$
and $$(\im\check\Omega)\wedge\check\omega=0.$$
Thus $\check\Omega_K=\im\check\Omega+i\check\omega$ is a 2-form
which satisfies the conditions of Theorem 5.1 and hence determines
a complex structure on $\check S$. We call $\check S$ with this
complex structure $\check S_K$.
As observed above, $[\check\Omega_n]$
must satisfy (7.3), so we need to look for a form $\re\check\Omega$
such that
$$[\re\check\Omega]={1\over Vol(S_b)}([\sigma_0]-\b-(\b^2-\omega^2-2)E/2).$$
If $[\re\check\Omega]$ is a K\"ahler class on $\check S_K$, 
then by Yau's theorem,
there exists a unique K\"ahler form $\re\check\Omega$ whose metric
is Ricci-flat, so we only need to ensure $[\re\check\Omega]$ is
a K\"ahler class. We first note that $[\re\check\Omega].[\im\check\Omega]
=[\re\check\Omega].[\check\omega]=0$, so $[\re\check\Omega]$ is a $(1,1)$
class.
Next we observe it is a positive class. Indeed,
$\check f:\check S_K\rightarrow B$ is still a holomorphic elliptic
fibration, and the complex structure on each fibre $\check S_{K,b}$
is the same as that of $S_{K,b}$. Since $[\re\check\Omega].E>0$,
this shows then that $[\re\check\Omega]$ is in fact positive
on $E$. As long as $Pic(\check S_K)$ contains no $-2$ classes, this
shows that $[\re\check\Omega]$ is a K\"ahler class. Hence we obtain
a K\"ahler form $\re\check\Omega$ as desired, and
set $\check\Omega=\re\check\Omega+i\im\check\Omega$.
$\bullet$

We make a few closing comments. First, in the above proof we can write
$$[\Omega]=Vol(S_b)([\sigma_0]+\check\b+i[\check\omega])\mod E$$
and 
$$[\check\Omega]={1\over Vol(S_b)}([\sigma_0]-(\b+i[\omega]))\mod E.$$
Thus we see that $Vol(\check S_b)=1/Vol(S_b)$, conforming with Remark
6.8. However this does not quite look like mirror symmetry: if we
repeat the process we appear to get
$$[\check{\check\Omega}]=Vol(S_b)([\sigma_0]-(\check\b+i[\check\omega]))
\mod E.$$
This again is the fault of our identifications, essentially having
put in a $90^{\circ}$ twist in $\T_{B}^*$ during
each dualising. This is rectified on the double mirror by pulling back all
forms by the fibrewise negation map on $S$. This acts trivially on 
$H^0(B,R^2f_*{\bf R})$ and $H^2(B,f_*{\bf R})$, but by negation
on $H^1(B,R^1f_*{\bf R})$.

Finally we note that the construction of the mirror K3 surface
given in this proof is of a different nature from previous constructions
of mirror symmetry for K3 surfaces.
Normally one appeals to the Torelli theorem to construct the mirror.
Here, once we have produced a special Lagrangian fibration on $S$,
we produce the mirror without an appeal to Torelli. Instead,
we are essentially
applying Yau's Theorem to solve the equations of Corollary 5.15.
However we are still aided by some key points which don't hold in
higher dimensions. These are:
\item{(1)} We know the harmonic $n-1$-forms on fibres, since $n-1=1$.
\item{(2)} We know the cohomology class of a holomorphic 2-form $\Omega$
if we know its class modulo $E$; this is completely determined by
the requirement $[\Omega]^2=0$.
\item{(3)} We can use the hyperk\"ahler trick.
\bigskip

{\hd Bibliography}
\item{[1]} Almgren, F. J., Jr. 
``$Q$ Valued Functions Minimizing Dirichlet's
Integral and the Regularity of Area Minimizing Rectifiable Currents up to
Codimension Two,'' {\it Bull. Amer. Math. Soc. (N.S.)}, {\bf 8}
(1983), 327--328.
\item{[2]} Almgren, F. J., Jr.,
``$Q$ Valued Functions Minimizing Dirichlet's
Integral and the Regularity of Area Minimizing Rectifiable Currents up to
Codimension Two,'' preprint.
\item{[3]} Arnold, V. I., {\it Mathematical Methods of Classical
Mechanics}, 2nd edition, Springer-Verlag, 1989.
\item{[4]} Aspinwall, P., ``An $N=2$ Dual Pair and a Phase Transition,''
{\it Nucl. Phys. B}, {\bf 460} (1996), 57--76.
\item{[5]} Aspinwall, P., and Morrison, D., ``Stable Singularities
in String Theory,'' with an appendix by M. Gross,
{\it Comm. Math. Phys.} {\bf 178}, (1996) 115--134.
\item{[6]} Aspinwall, P., and Morrison, D., ``String Theory on K3
surfaces,'' in {\it Essays on Mirror Manifolds II},
Greene, B.R., Yau, S.-T. (eds.) Hong Kong, International Press 1996, 703--716.
\item{[7]} Barth, W., Peters, C., and van de Ven, A., {\it
Compact Complex Surfaces}, Springer-Verlag, 1984.
\item{[8]} Besse, A., {\it Einstein Manifolds}, Springer-Verlag, 1987.
\item{[9]} Borel, A. et al, {\it Intersection Cohomology}, Birkh\"auser,
1984.
\item{[10]} Bredon, G., {\it Sheaf Theory}, 2nd edition, Springer-Verlag,
1997.
\item{[11]} Dolgachev, I., ``Mirror Symmetry for Lattice Polarized K3
surfaces,'' Algebraic Geometry, {\bf 4}, {\it J. Math. Sci.} {\bf 81},
(1996) 2599--2630. 
\item{[12]} Duistermaat, J., ``On Global Action-Angle Coordinates,''
{\it Comm. Pure Appl. Math.,} {\bf 33} (1980) 687--706.
\item{[13]} Federer, H., {\it Geometric Measure Theory}, Springer-Verlag,
1969.
%\item{[8]} Fukaya, K., ``Morse Homotopy, $A^{\infty}$-category,
%and Floer Homologies,'' in {\it Proceedings of GARC Worskshop on Geometry
%and Topology `93 (Seoul 1993)}, Lecture Notes Ser., {\bf 18}, Seoul Nat.
%Univ., Seoul, 1993, 1--102.
\item{[14]} Gross, M., ``Special Lagrangian Fibrations I: Topology,'' 
to appear in the Proceedings of the Taniguchi Symposium on Integrable
Systems and Algebraic Geometry.
\item{[15]} Gross, M., and Wilson, P.M.H., ``Mirror Symmetry via 3-tori for
a Class of Calabi-Yau Threefolds,'' 
{\it Math. Ann.}, {\bf 309}, (1997) 505--531.
\item{[16]} Harvey, R., and Lawson, H.B. Jr.,  ``Calibrated Geometries,'' {\it
 Acta
Math.} {\bf 148}, 47-157 (1982).
\item{[17]}
Hitchin, N., ``The Moduli Space of Special Lagrangian Submanifolds,''
preprint, dg-ga/9711002.
\item{[18]} King, J., ``The Currents Defined by Analytic Varieties,''
{\it Acta Math.}, {\bf 127} (1971), 185--220.
\item{[19]} Koszul, J.-L., ``Crochet de Schouten-Nijenhuis et Cohomologie,''
in {\it The Mathematical Heritage of Elie Cartan (Lyon, 1984)},
{\it Ast\'erisque}, Numero Hors Serie, (1985) 257--271.
%\item{[14]} Kontsevich, M., ``Homological Algebra of Mirror Symmetry,''
%in {\it Proceedings of the International Congress of Mathematicians, (Z\"urich,
%1994)}, Birkh\"auser, Basel, 1995, 120--139.
\item{[20]} Mangiarotti, L., and Modugno, M., ``Graded Lie Algebras and
Connections on a Fibered Space,'' {\it J. Math. Pures. et Appl.,}
{\bf 63}, (1984), 111--120.
\item{[21]} McLean, R.C., `` Deformations of Calibrated Submanifolds,'' 
Duke University Preprint, 1996.
\item{[22]} Morgan, F., {\it Geometric Measure Theory. A Beginner's Guide,}
2nd Edition,
{\it Academic Press, Inc.}, 1995. 
%\item{[17]} Morrison, D., ``Compactifications of Moduli Spaces Inspired
%by Mirror Symmetry,'' 
%In {\it Journ\'ees de G\'eometrie Alg\'ebrique d'Orsay, Juillet 1992,} 
%Asterisque
%{\bf 218}, 243-271 (1993).
\item{[23]} Morrison, D.R., ``The Geometry Underlying Mirror Symmetry,''
Preprint 1996.
\item{[24]} Simon, L.,  {\it Lectures on Geometric Measure Theory,}
Proceedings of the Centre for Mathematical Analysis, Australian National
University, {\bf 3}, {\it Australian National University, Centre
for Mathematical Analysis, Canberra}, 1983.
\item{[25]} Strominger, A., Yau, S.-T., and Zaslow, E.,  ``Mirror Symmetry is
T-Duality,'' {\it Nucl. Phys.} {\bf B479}, (1996) 243--259.
\item{[26]} Tian, G., ``Smoothness of the Universal Deformation Space of
Compact Calabi-Yau Manifolds and its Petersson-Weil Metric,'' in
{\it Mathematical Aspects of String Theory, (San Diego, California, 1986),}
629--646, {\it Adv. Ser. Math. Phys.} {\bf 1}, World Scientific
Publishing.
\item{[27]} Todorov, A., ``The Weil-Petersson Geometry of the Moduli Space of
$SU(n\ge 3)$
(Calabi-Yau) Manifolds I,'' {\it Commun. Math. Phys} {\bf 126}, (1989) 325--346.
\item{[28]} Weil, A., ``Final Report on Contract AF 18(603)-57,'' 
{\it Oeuvres Scientifiques}, Vol. II, 390--395.
\item{[29]} Zharkov, I., ``Torus Fibrations of Calabi-Yau Hypersurfaces in
Toric Varieties and Mirror Symmetry,'' preprint, alg-geom/9806091.
\end